\documentclass[12pt]{article}
\usepackage[dvips]{graphicx}
\usepackage{amsmath}
\usepackage{amssymb}
\usepackage{latexsym}
\usepackage{french}
\newcommand{\spa}{\hspace{5mm}}
\newcommand{\g}{\mathfrak{g}}
\newcommand{\nH}{{\cal H}^{\cal D}_{R}}
\newcommand{\Hr}{{\cal H}^{\cal D}_{P,R}}
\newcommand{\sHr}{{\cal H}_{P,R}}
\newcommand{\Ar}{{\cal A}^{\cal D}_{P,R}}
\newcommand{\esp}{\hspace{2mm}}
\newcommand{\strees}{{\cal T}_{P,R}}
\newcommand{\sforets}{{\cal F}_{P,R}}
\newcommand{\trees}{{\cal T}_{P,R}^{\cal D}}
\newcommand{\forets}{{\cal F}_{P,R}^{\cal D}}
\newcommand{\ntrees}{{\cal T}_{R}^{\cal D}}
\newcommand{\nforets}{{\cal F}_{R}^{\cal D}}
\newcommand{\ad}{{\cal A}d}
\newcommand{\cocy}{Z_{\varepsilon}^1}
\newcommand{\tdelta}{\tilde{\Delta}}
\newcommand{\ttop}{{\top}_{\hspace{-1mm} C}}
\newcommand{\otop}{\overline{\top}}
\newcommand{\mtop}{\tilde{\top}}
\newcommand{\treesb}{{\cal T}_{b}}
\newcommand{\Hfr}{{\cal H}_{Fr}}
\title{Les alg\`ebres de Hopf des arbres enracin\'es d\'ecor\'es}
\date{}
\author{L. Foissy \\
\\
{\small{\it Laboratoire de Math\'ematiques - UMR6056, Universit\'e de Reims}}\\
\small{{\it Moulin de la Housse - BP 1039 - 51687 REIMS Cedex 2, France}}\\
\small{e-mail: loic.foissy@univ-reims.fr}}

\begin{document}

\newtheorem{lemme}{Lemme}
\newtheorem{prop}[lemme]{Proposition}
\newtheorem{theo}[lemme]{Th\'eor\`eme}
\newtheorem{cor}[lemme]{Corollaire}
\newtheorem{defi}[lemme]{D\'efinition}

\maketitle
\section{Introduction}
\spa Dans \cite{Kreimer1,Connes}, une alg\`ebre de Hopf des arbres enracin\'es ${\cal H}_R$ est introduite dans le but d'\'etudier un probl\`eme de renormalisation. Cette alg\`ebre de Hopf est commutative, non 
cocommutative, et v\'erifie une propri\'et\'e universelle en cohomologie de Hochschild. On montre que le dual gradu\'e de cette alg\`ebre de Hopf est l'alg\`ebre enveloppante de l'alg\`ebre de Lie des arbres enracin\'es ; 
enfin, on construit dans \cite{Kreimer}  une bigraduation de ${\cal H}_R$ qui permet de calculer les dimensions des composantes homog\`enes de l'espace des \'el\'ements primitifs de ${\cal H}_R$, sans toutefois donner une expression pour ces \'el\'ements.

Dans ce papier, comme il \' etait sugg\'er\'e dans \cite{Connes},  nous introduisons une nouvelle alg\`ebre de Hopf $\Hr$, ni commutative, ni cocommutative, g\'en\'eralisant la construction pr\'ec\'edente. 
Apr\`es quelques r\'esultats pr\'eliminaires sur les alg\`ebres de Hopf gradu\'ees, nous montrons que cette alg\`ebre v\'erifie une propri\'et\'e universelle en cohomologie de Hochschild. On montre de plus que son dual gradu\'e $(\Hr)^{*g}$ est isomorphe \`a $\Hr$ ;
ce r\'esultat permet de construire un couplage de Hopf non d\'eg\'en\'er\'e $(\hspace{1mm},)$ entre $\Hr$ et elle-m\^eme. La base duale $(e_F)$ de la base des for\^ets est parti\-culi\`erement bien adapt\'ee au calcul du coproduit ; en particulier,
elle permet de trouver facilement les \'el\'ements primitifs de $\Hr$.  
De plus, nous donnons un sens combinatoire au couplage $(\hspace{1mm},)$, ce qui permet de donner une expression combinatoire directe des \'el\'ements primitifs. 
Ces r\'esultats permettent de montrer que $\sHr$ est isomorphe \`a l'alg\`ebre de Hopf des arbres plans binaires introduite dans \cite{Frabetti}.

Enfin, on d\'emontre plusieurs r\'esultats sur les cog\`ebres tensorielles ; appliqu\'es \`a ${\cal H}_R$, ils permettent de retrouver les r\'esultats de \cite{Connes}. On termine en montrant que les
alg\`ebres de Hopf des arbres plans d\'ecor\'es par un ensemble $\cal D$ fini ou d\'enombrable (voir \cite{Kreimer2}) sont toutes isomorphes.

\section{Dual gradu\'e}
\subsection{Cas d'un espace vectoriel}

\spa Soit $V$ un espace vectoriel sur un corps commutatif $K$.
On suppose $V$ muni d'une graduation $(V_n)_{n \in \mathbb{N}}$, telle que $dim(V_n)$ soit finie pour tout $n$.
Pour tout $x \in V$, $x \neq 0$, on pose $poids(x)=\min\{n / x \in V_0\oplus \ldots \oplus V_n\}$.\\

On identifie $V_n^*$ avec $\{f\in V^* / f(V_k)=(0) \mbox{ si } k \neq n\} \subseteq V^*$, et on pose $V^{*g}=\bigoplus V_n^*=\{f \in V^*/ \exists \esp n_0, f(V_n)=(0) \mbox{ si } n \geq n_0\}$.   
$V^{*g}$ est un  espace gradu\'e, avec $(V^{*g})_n=V_n^*$.

\begin{lemme}
\label{lemme1}
Soit $(e_i)_{i \in I}$ une base de $V$ form\'ee d'\'el\'ements homog\`enes. Pour $i \in I$, on d\'efinit $f_i \in V^*$ par $f_i(e_j)=\delta_{i,j}$, o\`u $\delta$ est le symbole de Kronecker.
Alors $(f_i)_{i \in I}$ est une base de $V^{*g}$.
\end{lemme}
{\it Preuve : }si $e_i \in V_k$, on a $f_i \in V_k^*$. Soit $J_k=\{j \in I/ e_j \in V_k\}$. Il suffit alors de montrer que $(f_j)_{j \in J_k}$ est une base de $V_k^*$. Comme $V_k$ est de dimension finie, c'est imm\'ediat. $\Box$

\begin{lemme}
\label{lemme2}
Soient $V$  et $W$ deux espaces gradu\'es et soit $\gamma :V \longmapsto W$, homog\`ene de degr\'e $k \in \mathbb{Z}$. 
Alors il existe une unique application $\gamma^{*g} :W^{*g} \longmapsto V^{*g}$, telle que :
$$\gamma^{*g}(f)(x) = f( \gamma(x)) \esp \forall f \in W^{*g}, \esp \forall x \in V.$$
De plus, $\gamma^{*g}$ est homog\`ene de degr\'e $-k$.
\end{lemme}
{\it Preuve  :}

Unicit\'e : soit $\gamma^* :W^* \longmapsto V^*$ la transpos\'ee de $\gamma$. On a alors $\gamma^{*g}=\gamma^*_{\mid W^{*g}}$.

Existence : il faut montrer que $\gamma^*(W^{*g}) \subseteq V^{*g}$. Soit $f\in W_n^*$. Soit $x \in V_i$, $i\neq n-k$.\\
$\gamma^*(f)(x)=f(\gamma(x))=0$ car $\gamma(x) \in W_{i+k}$, $i+k \neq n$. Donc $\gamma^*(W^*_n) \subseteq V^*_{n-k}$. $\Box$\\

On munit $V \otimes V$ d'une graduation donn\'ee par :
\begin{eqnarray*}
(V \otimes V)_n&=&\sum_{k+l=n} V_k \otimes V_l.
\end{eqnarray*}

\begin{lemme}
\label{lemme3}
On consid\`ere l'application suivante :
\begin{eqnarray*}
\theta_V :V^{*g}\otimes V^{*g}&\longmapsto & (V \otimes V)^{*g}\\
f\otimes g & \longmapsto &\left\{ \begin{array}{rcl}
V\otimes V &\longmapsto & K\\
x \otimes y & \longmapsto & f(x)g(y).
\end{array} \right.
\end{eqnarray*}
Alors $\theta_V$ est un isomorphisme d'espaces gradu\'es.
\end{lemme}
{\it Preuve :} classiquement, $\theta_V$ est injectif (voir \cite{Kassel}). Soit $f\in V_m^*$, $g \in V_n^*$, $x \in V_k$, $y \in V_l$, avec $k+l \neq m+n$.
$\theta_V(f\otimes g)(x \otimes y)=f(x)g(y)=0$ car $ k \neq m$ ou $l \neq n$. Donc $\theta_V$ est homog\`ene de degr\'e 0.
De plus,
\begin{eqnarray*}
dim((V^{*g}\otimes V^{*g})_n)&=& \sum_{k+l=n}dim(V_k^*)dim(V_l^*)\\
&=&  \sum_{k+l=n}dim(V_k)dim(V_l)\\
&=& dim((V\otimes V)_n^*).
\end{eqnarray*}
Donc $\theta_V$ est \'egalement surjectif. $\Box$
\\

Soit $V$ un espace gradu\'e, et $W$ un sous-espace de $V$. On dira que $W$ est un sous-espace gradu\'e de $V$  si $W=\bigoplus (W \cap V_n)$.

\begin{lemme}
\label{prop4}
Soit $W$ un sous-espace gradu\'e de $V$.
\begin{enumerate}
\item Soit $W_n = W \cap V_n$ et soit $W_n^{\perp_n}$ l'orthogonal de $W_n$ dans la dualit\'e entre $V_n$ et $V_n^*$.
 Alors dans la dualit\'e entre $V$ et $V^{*g}$, on a :
$$ W^\perp =\bigoplus_{n=0}^{+\infty} W^{\perp_n}.$$
\item De plus, $W^{\perp \perp}=W$.
\end{enumerate}
\end{lemme}
{\it Preuve : } 
On a $W^{\perp}= (\bigoplus W_n)^{\perp}=\bigcap (W_n^{\perp})$.
Or $W_n^{\perp}=\bigoplus_{i \neq n} V_i^* \oplus W_n^{\perp_n}$, donc $W^{\perp}=\bigoplus W_n^{\perp_n}$.
Comme $W_n^{\perp_n \perp_n}=W_n$, on a les r\'esultats annonc\'es. $\Box$
 
\subsection{Cas d'une alg\`ebre de Hopf}

\spa Soit $(A,m,\eta,\Delta,\varepsilon,S)$ une alg\`ebre de Hopf gradu\'ee sur un corps commutatif $K$, c'est-\`a-dire 
qu'il existe une graduation $(A_n)_{n \in \mathbb{N}}$ de l'espace vectoriel  ${A}$, avec :
\begin{eqnarray}
m(A_n \otimes A_m) &\subseteq & A_{n+m}, \esp \forall n,m \in \mathbb{N},\\
\Delta(A_n) &\subseteq & \sum_{k+l=n} A_k \otimes A_l, \esp \forall n \in \mathbb{N}.
\end{eqnarray}
(C'est-\`a-dire que $m$ et $\Delta$ sont homog\`enes de degr\'e 0)\\
On suppose de plus :\begin{description}
\item[$(C_1)$]  $dim(A_0)=1$ ;
\item[$(C_2)$]  les $A_n$ sont de dimension finie. 
\end{description}
On a alors $A_0=(1)$.
D'apr\`es \cite{Bou1}, proposition III.3.5,
 on a $Ker(\varepsilon)=\bigoplus_{n \geq 1} A_n.$

\begin{theo}
$A^{*g}$ est muni d'une structure d'alg\`ebre de Hopf gradu\'ee donn\'ee par :
\begin{enumerate}
\item $\forall f,g \in A^{*g}$, $\forall x \in A$, $(fg)(x)=(f \otimes g) ( \Delta(x))$ ;
\item $1_{A^{*g}}=\varepsilon$ ; 
\item $\forall f \in A^{*g}$, $\forall x,y \in A$, $\Delta(f)(x \otimes y)=f(xy)$ ;
\item $\forall f \in A^{*g}$, $\varepsilon(f)=f(1)$ ;
\item $\forall f \in A^{*g}$, $\forall x \in A$, $\left( S(f) \right)(x)=f\left( S(x) \right)$ ;
\item $(A^{*g})_n=A_n^*$.
\end{enumerate}
\end{theo}  
{\it Preuve :}
$m,$ $\eta$, $\Delta$, $\varepsilon$, $S$ sont homog\`enes de degr\'e 0
($K$ \'etant muni de la graduation triviale). 
En utilisant les lemmes \ref{lemme2} et \ref{lemme3}, on consid\`ere :
\begin{eqnarray*}
m^{*g}& :& A^{*g} \longmapsto (A\otimes A)^{*g} \stackrel{\theta_A^{-1}}{\longmapsto} A^{*g} \otimes A^{*g} \esp ;\\
\eta^{*g}& :&A^{*g} \longmapsto K \esp ;\\
\Delta^{*g}& :&A^{*g}\otimes A^{*g} \stackrel{\theta_A}{\longmapsto} (A\otimes A)^{*g} \longmapsto A^{*g} \esp ;\\
\varepsilon^{*g}& :& K \longmapsto A^{*g} \esp ;\\
S^{*g}& :& A^{*g} \longmapsto  A^{*g}.
\end{eqnarray*}
Classiquement, $(A^{*g},\Delta^{*g},\varepsilon^{*g},m^{*g},\eta^{*g},S^{*g})$ est une alg\`ebre de Hopf v\'erifiant $1-5$ (voir \cite{Kassel}).
Comme $ \Delta^{*g}$ et $m^{*g}$ sont homog\`enes de degr\'e 0, $6$ est v\'erifi\'ee. $\Box$\\

\begin{prop}
\label{prop6}
\begin{enumerate}
\item $(A^{*g})^{*g}$ et $A$ sont isomorphes comme alg\`ebres de Hopf gradu\'ees.
\item Soit $M$ l'id\'eal d'augmentation de $A$, c'est-\`a-dire $M=Ker(\varepsilon)$. 
Soit $Prim(A^{*g})=\{f \in A^{*g}/ \Delta(f)=1\otimes f+f \otimes1\}$.
Alors dans la dualit\'e entre $A$ et $A^{*g}$,
\begin{eqnarray*}
Prim(A^{*g})^{\perp}&=& (1) \oplus M^2,\\
\left( (1) \oplus M^2\right)^{\perp}&=&Prim(A^{*g}).
\end{eqnarray*}
\end{enumerate}
\end{prop}
{\it Preuve : }

$\begin{array}{rcl}
                   \mbox{1. Soit }i_n : \esp A_n & \longmapsto &(A_n^*)^*\\
                     x &  \longmapsto & \left\{  \begin{array}{rcl}
                                                 A_n^* &\longmapsto& K\\
                                                    f & \longmapsto & f(x).
                                                 \end{array}
                                                \right.
                                          \end{array}$ \\
Comme $A_n$ est de dimension finie, $i_n$ est un isomorphisme d'espaces vectoriels ; 
par suite, $i :A \longmapsto (A^{*g})^{*g}$ d\'efini par $i_{\mid A_n}=i_n$ est un isomorphisme d'espaces vectoriels gradu\'es.
On montre facilement qu'il s'agit d'un isomorphisme d'alg\`ebres de Hopf (voir \cite{Kassel}).\\

2. $(1) \subseteq Prim(A^{*g})^{\perp}$ : soit $p \in Prim(A^{*g})$. Alors $p(1)=\varepsilon(p)=0$.\\
$M^2 \subseteq Prim(A^{*g})^{\perp}$ : soit $m \in M^2$. On peut supposer $m=m_1m_2,$ $\varepsilon(m_1)=\varepsilon(m_2)=0$.
Soit $p\in Prim(A^{*g})$. On a :
\begin{eqnarray*}
p(m_1m_2)&=&\Delta(p)(m_1 \otimes m_2)\\
&=&(1\otimes p + p \otimes 1)(m_1 \otimes m_2)\\
&=& \varepsilon(m_1)p(m_2) + p(m_1) \varepsilon(m_2)\\
&=&0.\\
\end{eqnarray*}
$\left((1) \oplus M^2\right)^{\perp} \subseteq Prim(A^{*g})$ : soit $f \in \left((1) \oplus M^2\right)^{\perp}$.
Il s'agit de montrer  : $\forall x,y \in A$, $\Delta(f)(x \otimes y)=(f\otimes 1 + 1 \otimes f)(x \otimes y)$, 
c'est-\`a-dire : $f(xy)=f(x)\varepsilon(y)+\varepsilon(x)f(y)$. Comme $A=(1)\oplus Ker(\varepsilon)$, il suffit de consid\'erer les 4 cas suivants :
\begin{enumerate}
\item $x=y=1$ : il faut montrer $f(1)=2 f(1)$, ce qui est vrai car $f \in (1)^{\perp}$ ;
\item $x=1$, $\varepsilon(y)=0$ : \'evident ;
\item $\varepsilon(x)=0$, $y=1$ : \'evident ;
\item $\varepsilon(x)=\varepsilon(y)=0$ : alors $xy \in M^2$, et donc $f(xy)=0.$
\end{enumerate}
On a montr\'e $(1) \oplus M^2 \subseteq Prim(A^{*g})^{\perp}$ et  $\left((1) \oplus M^2\right)^{\perp} \subseteq Prim(A^{*g})$. Comme $A$ et $A^{*g}$ sont des alg\`ebres de Hopf gradu\'ees, $Prim(A^{*g})$ et $M^2$ sont des sous-espaces gradu\'es. On obtient donc les inclusions r\'eciproques par passage \`a l'orthogonal en utilisant le lemme \ref{prop4}.  $\Box$

\subsection{Cas des alg\`ebres enveloppantes}

\spa On suppose  que $\mathfrak{g}$ est une alg\`ebre de Lie gradu\'ee, avec $\g_0=(0)$, et les $\g_n$ de dimension finie.
Alors ${\cal U}(\g)$ est une alg\`ebre de Hopf gradu\'ee v\'erifiant $(C_1)$ et $(C_2)$.
On suppose de plus que $K$ est de caract\'eristique nulle. 
\begin{prop}
\label{prop7}
Soit $M=Ker(\varepsilon_{{\cal U}(\g)^{*g}} )=\bigoplus_{n \geq 1}{\cal U}(\g)^*_n$. Soit $V$ un suppl\'ementaire gradu\'e de $M^2$ dans $M$. Soit $S(V)$ l'alg\`ebre sym\'etrique de $V$.
On consid\`ere le morphisme d'alg\`ebres suivant :
\begin{eqnarray*}
\xi_V :S(V)&\longmapsto & {\cal U}(\g)^{*g}\\
      v \in V & \longmapsto & v.
\end{eqnarray*}
Alors $\xi_V$ est un isomorphisme d'alg\`ebres gradu\'ees.
\end{prop}
{\it Preuve :} remarquons d'abord que ${\cal U}(\g)^{*g}$ est une alg\`ebre commutative, et donc $\xi_V$ est bien d\'efini.

Soit $(e_i)_{i \in I}$ une base de $\mathfrak{g}$ form\'ee d'\'el\'ements homog\`enes, o\`u $I=\{1, \ldots,n\}$ ou $\mathbb{N}^*$ suivant la dimension de $\mathfrak{g}$.
On pose $S_I$ l'ensemble des \'el\'ements $(\nu_1,\ldots,\nu_k,\ldots)$ de $\mathbb{N}^I$ tels que les $\nu_j$ soient presque tous nuls.
Pour $\nu=(\nu_1,\ldots,\nu_k,\ldots)\in S_I$, avec $k$ tel $\nu_j=0$ si $j>k$, 
on pose : $$e_{\nu}=\frac{e_1^{\nu_1}\ldots e_k^{\nu_k}}{\nu_1!\ldots \nu_k!}.$$
Suivant \cite{Dixmier}, pages 91 et suivantes,
pour tout $\nu \in S_I$ :
$$ \Delta(e_{\nu})=\sum_{\lambda+\mu=\nu} e_{\lambda} \otimes e_{\mu}.$$
D'apr\`es le th\'eor\`eme de Poincar\'e-Birkhoff-Witt, $(e_{\nu})_{\nu \in S_I}$ est une base de ${\cal U}(\g)$, form\'ee
d'\'el\'ements homog\`enes. Soit $(f_{\nu})_{\nu \in S_I}$ la base de ${\cal U}(\g)^{*g}$, d\'efinie par  $f_{\nu}(e_{\mu})=\delta_{\nu,\mu}.$
\begin{eqnarray*}
f_{\nu} f_{\mu} (e_{\lambda})&=&(f_{\nu}\otimes  f_{\mu}) \Delta(e_{\lambda})\\
&=&(f_{\nu}\otimes  f_{\mu})(\sum_{\alpha+\beta=\lambda} e_{\alpha} \otimes e_{\beta})\\
&=&\delta_{\nu+\mu,\lambda}.
\end{eqnarray*}
Par suite, $f_{\nu} f_{\mu}=f_{\nu +\mu}$. On consid\`ere alors  $V_0$ l'espace engendr\'e par les $f_{\nu}$, $\sum \nu_i=1$, ainsi que :
\begin{eqnarray*}
\xi_{V_0} :S(V_0)& \longmapsto & {\cal U}(\g)^{*g}\\
f_{\nu} \in V_0 & \longmapsto & f_{\nu}.
\end{eqnarray*}
On d\'eduit alors imm\'ediatement du calcul pr\'ec\'edent que $\xi_{V_0}$ est un isomorphisme d'alg\`ebres.
De plus, comme $Prim({\cal U}(\g))=\g$ (car $K$ est de caract\'eristique nulle, voir \cite{Bou2}),
on a ${\cal U}(\g)^{*g}=V_0 \oplus Prim({\cal U}(\g))^{\perp}$, donc $V_0$ est un suppl\'ementaire gradu\'e de 
$M$ dans $M^2$ d'apr\`es la proposition \ref{prop6}-2. Il est \'evident que $\xi_{V_0}$ est homog\`ene de degr\'e z\'ero.
\\

Consid\'erons $V$ un suppl\'ementaire gradu\'e de $M^2$ dans $M$. Montrons que $V$ g\'en\`ere ${\cal U}(\g)^{*g}$. Soit $f  \in {\cal U}(\g)^{*g}$, de poids $n$. Si $n=0$, alors $f \in <V>$, la sous-alg\`ebre engendr\'ee par $V$. Supposons que tous les \'el\'ements de poids strictement inf\'erieur \`a $n$ soient dans $<V>$. On peut alors supposer $f$
homog\`ene ; alors $f \in M=M^2\oplus V$ : on peut donc supposer $f=f_1f_2$, $f_i \in M$, et donc $poids(f_i)<n$ pour $i=1,2$. Par suite les $f_i$ sont dans $<V>$, et donc $f\in <V>$.

 Donc $V$ g\'en\`ere ${\cal U}(\g)^{*g}$. Par suite, $\xi_{V}$ est surjectif, homog\`ene de degr\'e z\'ero. De plus, $V$ et $V_0$ sont des espaces gradu\'es isomorphes, donc $S(V)$ et $S(V_0)$ sont des alg\`ebres gradu\'ees isomorphes.
Pour tout $n \in \mathbb{N}$, on a donc $dim(S(V)_n)=dim(S(V_0)_n)=dim({\cal U}(\g)^{*g}_n)$ ;  on en d\'eduit que $\xi_{V}$ est un isomorphime. $\Box$

\section{R\'esultats sur les cog\`ebres et les alg\`ebres de Hopf gradu\'ees}

\subsection{Filtration par $deg_p$}
\label{partB}

\begin{lemme}
\label{lemme9}
Soient $(C,\Delta, \varepsilon)$ une cog\`ebre, et $e \in C$ tel que $\Delta(e)=e \otimes e$. On pose :
$$\tilde{\Delta}(x)=\Delta(x)-e\otimes x -x \otimes e, \esp \forall x \in C.$$
Alors $\tdelta$ est coassociatif, c'est-\`a-dire : pour tout $x\in C,$ $(\tilde{\Delta}\otimes Id) \circ \tilde{\Delta}(x)=(Id \otimes \tilde{\Delta}) \circ \tilde{\Delta}(x)$.
\end{lemme}
{\it Preuve :}
pour tout $y \in C$, on pose $\tilde{\Delta}(y)=\sum y' \otimes y''$.
\begin{eqnarray*}
(\Delta \otimes Id) \circ {\Delta}(x)&=& e \otimes e \otimes x + e \otimes x \otimes e + x \otimes e \otimes e + \sum x' \otimes x'' \otimes e\\
&& +\sum x' \otimes e \otimes x'' +\sum e\otimes x' \otimes x'' + \sum \sum {(x')}' \otimes {(x')}'' \otimes x'' \esp ;\\
(Id \otimes \Delta) \circ {\Delta}(x)&=& e \otimes x \otimes e + e \otimes e \otimes x + \sum e \otimes x' \otimes x''+ x \otimes e \otimes e \\ 
&&+ \sum x' \otimes x'' \otimes e + \sum x' \otimes e \otimes x'' + \sum \sum x' \otimes {(x'')}^{'} \otimes {(x'')}^{''}.
\end{eqnarray*} 
Comme $\Delta$ est coassociatif, les deux membres de droite sont \'egaux. On en d\'eduit alors le r\'esultat voulu. $\Box$\\

Soit $C$ une cog\`ebre gradu\'ee v\'erifiant la condition $(C_1)$. 
D'apr\`es \cite{Bou1}, il existe un unique $x\in C$ non nul tel que $\Delta(x)=x \otimes x$. Cet \'el\'ement sera not\'e $1$. Il est homog\`ene de poids z\'ero.
On a de plus $M=Ker(\varepsilon)=\bigoplus_{n \geq 1} C_n$, et $\forall x \in M$, $\Delta(x)=x\otimes 1+1 \otimes x +M \otimes M$. Enfin, $C^{*g}$ est muni d'une structure d'alg\`ebre donn\'ee par :
$$(fg)(x)=(f\otimes g)\left(\Delta(x)\right) \esp \forall f,g \in C^{*g}, \esp \forall x \in C.$$
On d\'efinit $\tilde{\Delta}(x)=\tilde{\Delta}^1(x)=\Delta(x) -1\otimes x - x \otimes 1$ pour tout $x \in C$, et par r\'ecurrence
$\tilde{\Delta}^k=(\tilde{\Delta}^{k-1} \otimes Id) \circ \tilde{\Delta}$.

\begin{lemme}
\label{lemme10}
Soit $\rho$ la projection sur $\bigoplus_{n \geq 1} C_n$ parall\`element \`a $C_0$. Alors pour tout $k \geq 1$ :
$$ \tilde{\Delta}^k \circ \rho=\rho^{\otimes (k+1)}\circ \Delta^k.$$
\end{lemme}
{\it Preuve :}
on remarque que $\rho(x)=x - \varepsilon(x)1$, $\forall x \in C$.
Montrons par r\'ecurrence que $\tilde{\Delta}^k \circ \rho=\rho^{\otimes k+1} \circ \Delta^k$.
Pour $k=1$, c'est imm\'ediat si $x=1$, et d\'ecoule des remarques pr\'ec\'edentes si $ \varepsilon(x)=0$. Supposons l'hypoth\`ese de r\'ecurrence vraie
au rang $k$ : 
\begin{eqnarray*}
\tilde{\Delta}^{k+1} \circ \rho&=&(\tilde{\Delta} \otimes Id^{\otimes k}) \circ \tilde{\Delta}^k \circ \rho\\
&=&(\tilde{\Delta} \otimes Id^{\otimes k}) \circ {\rho}^{\otimes k+1} \circ \Delta^k\\
&=&\rho^{\otimes k+2} \circ (\Delta \otimes Id^{\otimes k}) \circ \Delta^k\\
&=&\rho^{\otimes k+2} \circ \Delta^{k+1}.
\end{eqnarray*}
(On a utilis\'e l'hypoth\`ese de r\'ecurrence pour la deuxi\`eme \'egalit\'e, et le r\'esultat avec $k=1$ pour la troisi\`eme.) $\Box$

\begin{lemme}
\label{lemme11}
Soit $x\in C$, tel que  $\tilde{\Delta}^{n}(x)=0$. Alors $\tilde{\Delta}^{n-1}(x) \in Prim(C)^{\otimes n}$.
\end{lemme}
{\it Preuve :} comme $\tilde{\Delta} $ est coassociatif,
 $\tilde{\Delta}^{n-1}(x) \in Ker(Id^{\otimes (i-1)} \otimes \tilde{\Delta} \otimes Id^{\otimes(n-i)})$ $\forall i \in \{1,\ldots ,n\}$. Donc
$\tilde{\Delta}^{n-1}(x) \in \bigcap \left(C^{\otimes (i-1)} \otimes Prim(C) \otimes C^{\otimes (n-i)}\right)= Prim(C)^{\otimes n}$. $\Box$

\begin{lemme}
\label{lemme12}
Pour tout $x \in M$, on a $\tilde{\Delta}^{poids(x)}(x)=0$.
\end{lemme}
{\it Preuve :} par r\'ecurrence sur $poids(x)$. Si $poids(x)=1$, alors  $x$ est n\'ecessairement primitif, et donc $\tilde{\Delta}(x)=0$.
Supposons l'hypoth\`ese de r\'ecurrence vraie au rang $n$. Soit $x$ de poids $n$ ; on pose $\tilde{\Delta}(x)=\sum x' \otimes x'',$ $poids(x')<n$. 
Par coassociativit\'e de $\tilde{\Delta}$, on a :
$$ \tilde{\Delta}^{n}(x)=\sum \tilde{\Delta}^{n-1}(x') \otimes x''=0. \esp \Box$$

On pose $C_{deg_p \leq n}=(1) \oplus Ker(\tilde{\Delta}^n)=Ker(\tilde{\Delta}^n \circ \rho)$. D'apr\`es le lemme pr\'ec\'edent, c'est une filtration de l'espace $C$.
Pour $x \in C$, $x \neq 0$, on pose $deg_p(x)=\min\{n/ \tilde{\Delta}^n\circ \rho (x)=0\}$, de sorte que $C_{deg_p\leq n}=\{x \in C/deg_p(x)\leq n\}$.
On a $deg_p(x) \leq poids(x)$ $\forall x\neq 0$.

\begin{prop}
\label{proF}
Soit $M_*=(1)^{\perp}\subset C^{*g}$. Alors dans la dualit\'e entre $C$ et $C^{*g}$, on a :
$$(C_{deg_p \leq n})^{\perp}=M_*^{n+1}.$$
\end{prop}
{\it Preuve :}
c'est \'evident si $n=0$. Supposons $n \geq 1$. Posons $\rho_*$ la projection sur $\bigoplus_{n \geq 1} C^*_n=M_*$ parall\`element \`a $C_0^*$.
On a facilement $\rho_*=\rho^{*g}$. Soient $m_1, \ldots, m_{n+1} \in M_*$, $x \in C_{deg_p \leq n}$.  
\begin{eqnarray*}
(m_1 \ldots m_{n+1},x)&=&(m_1 \otimes \ldots \otimes m_{n+1}, \Delta^n(x))\\
&=&(\rho_*^{\otimes(n+1)}(m_1 \otimes \ldots \otimes m_{n+1}), \Delta^n(x))\\
&=&(m_1 \otimes \ldots \otimes m_{n+1}, \rho^{\otimes(n+1)}\circ\Delta^n(x))\\
&=&(m_1 \otimes \ldots \otimes m_{n+1}, \tilde{\Delta}^n(\rho(x)))\\
&=&(m_1 \otimes \ldots \otimes m_{n+1}, 0)\\
&=&0.
\end{eqnarray*}
(On a utilis\'e le lemme \ref{lemme10} pour la quatri\`eme \'egalit\'e.)\\
Donc $M_*^{n+1}\subseteq (C_{\deg_p\leq n})^{\perp}.$\\

Pour montrer l'inclusion r\'eciproque, il suffit de montrer que $(M_*^{n+1})^{\perp} \subseteq C_{deg_p \leq n}$.
Soit $x \in (M_*^{n+1})^{\perp}$, $m_1, \ldots , m_{n+1} \in M_*$.
\begin{eqnarray*}
(m_1 \otimes \ldots \otimes m_{n+1}, \tilde{\Delta}^n\circ \rho (x)))&=&(m_1 \ldots m_{n+1},x)\\
&=&0.
\end{eqnarray*}
Donc $\tilde{\Delta}^n\circ \rho (x)\in M^{\otimes (n+1)}\cap {(M_*^{\otimes (n+1)})}^\perp=(0)$.
Par suite, $\tilde{\Delta}^n\circ \rho(x)=0$, et donc $x \in C_{deg_p\leq n}$. $\Box$

\begin{cor}
$(C_{deg_p\leq n})_{n \in \mathbb{N}}$ est une filtration de cog\`ebre, c'est-\`a-dire :
$$\Delta(C_{deg_p\leq n}) \subseteq \sum_{k+l=n} C_{deg_p \leq k} \otimes C_{deg_p \leq l}.$$
\end{cor}
\begin{lemme}
 Soit $A$ une alg\`ebre et $M$ un id\'eal  de $A$. Pour tout $n \in \mathbb{N}$, on a :
\begin{eqnarray*}
\bigcap_{k+l=n} \left( M^{k+1} \otimes A + A\otimes  M^{l+1} \right) &=& \sum_{i+j>n} M^i \otimes M^j.
\end{eqnarray*}
\end{lemme}
 {\it Preuve :}  \\
$\supseteq $ : soit $x =x_i \otimes x_j \in M^i \otimes M^j$, $i+j>n$. Soient $k,l \in \mathbb{N}$, $k+l=n$. Si $k<i$, alors $x_i \in M^i \subseteq M^{k+1}$, et donc $x \in M^{k+1}\otimes A$. Sinon, $l<j$ car $k+l\leq n<i+j$. Donc $x_j \in M^{l+1}$, et $x \in A \otimes M^{l+1}$.\\
$\subseteq$ : pour tout $i\in \mathbb{N}$, soit $W_i$ un suppl\'ementaire de $M^{i+1}$ dans $M^i$, et soit $W_\infty = \cap M^i$. On a alors :
\begin{eqnarray*}
A&=& \bigoplus_{i=0}^{\infty} W_k,\\
M^k&=&\bigoplus_{i=k}^{\infty} W_k \esp ;\\ 
A\otimes A &=& \bigoplus_{i,j=0}^{\infty} W_k\otimes W_j.
\end{eqnarray*}
Pour $i,j \in  \mathbb{N}\cup \{\infty\}$, soit $p_{i,j}:A\otimes A \longmapsto W_i\otimes W_j$ la projection sur $W_i \otimes W_j$ dans cette somme directe.
Soit $x \in \bigcap (M^{k+1} \otimes A+A\otimes M^{l+1})$.  Soit $k \in \mathbb{N}$.  Comme $x \in M^{k+1} \otimes A+A\otimes M^{n-k+1}$, 
si $i \leq k$ et $j \leq n-k$, alors $p_{i,j}(x)=0$. Par suite, si $p_{i,j}(x) \neq 0$, alors pour tout $k \in \mathbb{N}$, $i >k$ ou $j>n-k$. En particulier pour $k=i$,
on a $j>n-i$. Donc :
$$x  \in  \bigoplus_{i+j>n} W_i \otimes W_j \subseteq \ \sum_{i+j>n}  M^i \otimes M^j.\esp \Box$$

\noindent {\it Preuve du corollaire :}
soit $x \in C_{deg_p \leq n}$ ; montrons que $\Delta(x) \in \sum_{k+l=n} C_{deg_p \leq k} \otimes C_{deg_p \leq l}.$ On a :
\begin{eqnarray}
\nonumber \sum_{k+l=n} C_{deg_p \leq k} \otimes C_{deg_p \leq l}&=& \sum_{k+l=n} (M^{k+1}_*)^{\perp} \otimes (M^{l+1}_*)^{\perp}\\
\nonumber &=&\sum_{k+l=n} (M^{k+1}_* \otimes C^{*g}+C^{*g}\otimes M^{l+1}_*)^{\perp}\\
\nonumber &=&\left( \bigcap_{k+l=n}(M^{k+1}_* \otimes C^{*g}+C^{*g} \otimes M_*^{l+1})\right)^{\perp}\\
\label{eqnH} &=&\left( \sum_{i+j >n} M^i_* \otimes M^j_*\right)^{\perp}.
\end{eqnarray}
(On a utilis\'e le lemme pour la derni\`ere \'egalit\'e, avec $A=C^{*g}$ et $M=M_*$.)\\
Soit $f_1 \in M^i_*$, $f_2 \in M^j_*$, $i+j>n$. Alors $f_1 f_2 \in M_*^{n+1}=(C_{deg_p \leq n})^{\perp}$, donc :
\begin{eqnarray*}
(\Delta(x),f_1 \otimes f_2)&=&(x,f_1f_2)\\
&=&0.
\end{eqnarray*}
D'apr\`es (\ref{eqnH}), $\Delta(x) \in  \sum_{k+l=n} C_{deg_p \leq k} \otimes C_{deg_p \leq l}.$ $\Box$

\subsection{Cas d'une alg\`ebre de Hopf gradu\'ee}
\spa Soit $A$ une alg\`ebre de Hopf gradu\'ee v\'erifiant $(C_1)$. On suppose de plus que $K$ est de caract\'eristique nulle.\\

$S_n$ agit sur $A^{\otimes n}$ par $\sigma.(x_1 \otimes \ldots \otimes x_n)=x_{\sigma(1)} \otimes \ldots \otimes x_{\sigma(n)}$.
On rappelle qu'un $(p,q)$-$battage$ est un \'el\'ement $\sigma$ de $S_{p+q}$, croissant sur $\{1, \ldots,p\}$ et sur
$\{p+1, \ldots, p+q\}$.  On note $bat(p,q)$ l'ensemble des $(p,q)$-battages.

\begin{lemme}
\label{lemI}
\begin{enumerate}
\item Soient $x,y \in A-\{0\}$, $deg_p(x)=p$, $deg_p(y)=q$. On suppose que $\varepsilon(x)=\varepsilon(y)=0.$
 On pose $\tilde{\Delta}^{p-1}(x)=\sum_i x^{(1)}_i \otimes \ldots \otimes x^{(p)}_i$ et
$\tilde{\Delta}^{q-1}(y)=\sum_j y^{(1)}_j \otimes \ldots \otimes y^{(q)}_j$, les $x^{(k)}_i$ et les $y^{(l)}_j$ \'etant primitifs (lemme \ref{lemme11}).
Alors :
$$\tilde{\Delta}^{p+q-1}(xy)=\sum_{i,j} \sum_{\sigma \in bat(p,q)} \sigma.(x^{(1)}_i \otimes \ldots \otimes x^{(p)}_i \otimes y^{(1)}_j \otimes \ldots \otimes y^{(q)}_j).$$ 
\item Soient $p_1,\ldots, p_n \in Prim(A)$. On a :
$$\tilde{\Delta}^{n-1}(p_1 \ldots p_n)=\sum_{\sigma \in S_n} p_{\sigma(1)} \otimes \ldots \otimes p_{\sigma(n)}.$$ 
\end{enumerate}
\end{lemme}
{\it Preuve :} 

1. D'apr\`es le lemme \ref{lemme10} :
\begin{eqnarray*}
\tilde{\Delta}^{p+q-1}(xy)&=&\tilde{\Delta}^{p+q-1}(\rho(xy))\\
&=& \rho^{\otimes(p+q)} \left( \Delta^{p+q-1}(xy)\right)\\
&=& \rho^{\otimes(p+q)} \left( \Delta^{p+q-1}(x) \Delta^{p+q-1}(y)\right).
\end{eqnarray*}
Le r\'esultat est alors imm\'ediat.\\

2. R\'ecurrence sur $n$. La formule est vraie pour $n=2$. Supposons l'hypoth\`ese de r\'ecurrence vraie au rang $n-1$.
D'apr\`es 1, on a :
$$\tilde{\Delta}^{n-1}(p_1 \ldots p_n)= \sum_{\tau \in bat(n-1,1)}\esp \sum_{\sigma \in S_{n-1}} \tau. (p_{\sigma(1)} \otimes \ldots \otimes p_{\sigma(n-1)} \otimes p_n).$$
Or il y a exactement $n$ $(n-1,1)$-battages : $\tau \in bat(n-1,1)$ est caract\'eris\'e par $\tau(n)$. Alors :
\begin{eqnarray*}
\tilde{\Delta}^{n-1}(p_1 \ldots p_n)
&=& \sum_{i=1}^{n} \sum_{\sigma \in S_{n-1}} (p_{\sigma(1)} \ldots \otimes p_{\sigma(i-1)} \otimes p_n \otimes p_{\sigma(i)} \otimes \ldots \otimes p_{\sigma(n-1)})\\
&=&\sum_{\sigma \in S_n} p_{\sigma(1)} \otimes \ldots \otimes p_{\sigma(n)}. \esp \Box
\end{eqnarray*}
\begin{prop}
\label{proJ}
$(A_{deg_p\leq n})_{n \in \mathbb{N}}$ est une filtration de l'alg\`ebre de Hopf $A$.
\end{prop}
{\it Preuve :}
on a d\'ej\`a  vu qu'il s'agissait d'une filtration de cog\`ebre.
Soit $x,y\in A$, $deg_p(x) = p$, $deg_p(y)= q$. Alors d'apr\`es 
le lemme \ref{lemI}, $\tilde{\Delta}^{p+q-1}(\rho(xy)) \in Prim(A)^{\otimes (p+q)}$, et donc $\tilde{\Delta}^{p+q}(\rho(xy))=0$. Donc $deg_p(xy)\leq p+q$.
On a donc bien une filtration d'alg\`ebre. $\Box$

\subsection{Alg\`ebres de Hopf gradu\'ees cocommutatives ou commutatives}
\spa On suppose toujours $K$ de caract\'eristique nulle.

Les consid\'erations pr\'ec\'edentes permettent de donner la preuve suivante du th\'eor\`eme de Cartier-Milnor-Moore-Quillen :

\begin{theo}[Cartier-Milnor-Moore-Quillen]
\label{theoK}
Soit $A$ une alg\`ebre de Hopf cocommutative gradu\'ee v\'erifiant $(C_1)$. Alors $A$ est isomorphe \`a
${\cal U}(Prim(A))$ comme alg\`ebre de Hopf gradu\'ee.
\end{theo}
{\it Preuve :} montrons d'abord que $Prim(A)$ g\'en\`ere l'alg\`ebre $A$. Soit $\cal P$ la sous alg\`ebre de $A$ engendr\'ee par $Prim(A)$.
Soit $x\in A$, montrons que $x \in \cal P$. On peut se ramener \`a $\varepsilon(x)=0$. Proc\'edons par r\'ecurrence sur $deg_p(x)$.
Si $deg_p(x)=1$, alors $x$ est primitif, et donc $x \in \cal P$. Supposons l'hypoth\`ese de r\'ecurrence vraie au rang $n-1$, et supposons $deg_p(x)=n$.
D'apr\`es le lemme \ref{lemme11}, on peut poser $\tilde{\Delta}^{n-1}(x)=\sum_i x_i^{(1)} \otimes \ldots \otimes x_i^{(n)}$, les $x_i^{(j)}$ \'etant primitifs.
Comme $A$ est cocommutative, on a :
\begin{eqnarray*}
\tilde{\Delta}^{n-1}(x)&=&\sum_i x_i^{(1)} \otimes \ldots \otimes x_i^{(n)}\\
 &=&\frac{1}{n!}  \sum_i \sum_{\sigma \in S_n} x_i^{(\sigma(1))} \otimes \ldots \otimes x_i^{(\sigma(n))}\\
&=&\frac{1}{n!}\sum_i \tilde{\Delta}^{n-1}(x_i^{(1)} \ldots x_i^{(n)}).
\end{eqnarray*}
(On a utilis\'e le lemme \ref{lemI}-2 pour la derni\`ere \'egalit\'e).\\
Donc $deg_p(x-\frac{1}{n!}\sum x_i^{(1)} \ldots x_i^{(n)})<n$, et donc $ x-\frac{1}{n!}\sum x_i^{(1)} \ldots x_i^{(n)}\in \cal P$ ;
on en d\'eduit que $x\in \cal P$.

Par suite, on a un morphisme surjectif d'alg\`ebres de Hopf, homog\`ene de poids z\'ero :
\begin{eqnarray*}
\xi :{\cal U}(Prim(A))& \longmapsto & A \\
 p \in Prim(A)&\longmapsto & p.
\end{eqnarray*}
Comme $A$ v\'erifie la condition $(C_1)$, la composante homog\`ene de poids 0 de $Prim(A)$ est nulle, et donc ${\cal U}(prim(A))$ v\'erifie aussi la condition $(C_1)$.
Par suite ${\cal U}(Prim(A))$ est filtr\'ee par $deg_p$ (proposition \ref{proJ}). \\

Supposons $\xi$ non injectif, et soit $x$ non nul tel que $\xi(x)=0$. On choisit $x$ de $deg_p$ minimal. Comme $\xi(p)=p$ pour tout $p$ primitif,
n\'ecessairement $deg_p(x) >1$. Posons $\Delta(x) =x \otimes 1 +1 \otimes x +\sum x' \otimes x''$, $deg_p(x')<deg_p(x)$, $deg_p(x'')<deg_p(x)$.
$\xi$ est un morphisme de cog\`ebres, donc :
\begin{eqnarray*}
\Delta(\xi(x))&=&\xi \otimes \xi ( \Delta(x))\\
&=&\xi \otimes \xi (x \otimes 1 +1 \otimes x +\sum x' \otimes x'')\\
&=& \sum \xi(x') \otimes \xi(x'')\\
&=&0.
\end{eqnarray*}
Par choix de $x$, $\xi_{\mid A_{deg_p\leq deg_p(x)-1}}$ est injectif, d'o\`u $\sum x' \otimes x''=0$. Donc $x$ est primitif, et donc $deg_p(x) =1$ : contradiction.
Donc $\xi$ est injectif. $\Box$

\begin{prop}
\label{proL}
Soit $A$ une alg\`ebre de Hopf gradu\'ee commutative v\'erifiant $(C_1)$ et $(C_2)$. Posons $M=Ker(\varepsilon)$.
Alors si $G$ est un suppl\'ementaire gradu\'e de $M^2$ dans $M$, on a un isomorphisme d'alg\`ebres gradu\'ees :
\begin{eqnarray*}
\xi_G :S(G)&\longmapsto & A\\
g \in G & \longmapsto& g.
\end{eqnarray*}
\end{prop}
{\it Preuve :} $A^{*g}$ est une alg\`ebre cocommutative. En appliquant le th\'eor\`eme pr\'ec\'edent et la proposition \ref{prop7}, on obtient
imm\'ediatement le r\'esultat, car $A$ et $(A^{*g})^{*g}$ sont isomorphes d'apr\`es la proposition \ref{prop6}-1. $\Box$

\begin{cor}
\label{corM}
Soit $A$ une alg\`ebre de Hopf gradu\'ee v\'erifiant $(C_1)$ et $(C_2)$. Alors $\forall x,y \in A$, on a :
\begin{eqnarray*}
 deg_p(xy)&=&deg_p(x)+deg_p(y) \esp ;\\
 poids(xy)&=&poids(x)+poids(y).
\end{eqnarray*}
\end{cor}
{\it Preuve :} il s'agit de montrer que l'alg\`ebre gradu\'ee associ\'ee \`a la fitration par $deg_p$ est int\`egre. En utilisant le r\'esultat pr\'ec\'edent, il suffit de montrer que 
c'est une alg\`ebre commutative. Soient $x,y \in A$, $\varepsilon(x)=\varepsilon(y)=0$. On utilise les notations du lemme \ref{lemI}-1. 
Soit $\sigma \in bat(p,q)$. Soit $\tilde{\sigma} \in S_{p+q}$  d\'efini par : 
\begin{eqnarray*}
\tilde{\sigma}(i) &=&\sigma( i+p) \mbox{ si $i \leq q$}\\
&=&\sigma(i-q) \mbox{ si $i >q$.}
\end{eqnarray*}

On remarque imm\'ediatement que $\sigma \longmapsto \tilde{\sigma}$ est une bijection de $bat(p,q)$ vers $bat(q,p)$.\\
 De plus,
$\sigma.(x_1 \otimes \ldots \otimes x_p \otimes y_1 \otimes \ldots \otimes y_q)=
\tilde{\sigma}.(y_1 \otimes \ldots \otimes y_q \otimes x_1 \otimes \ldots \otimes x_p)$ $\forall x_i,y_j \in A$. Donc :
\begin{eqnarray*}
\tilde{\Delta}^{p+q-1}(xy)&=&\sum_{i,j} \sum_{\sigma \in bat(p,q)} \sigma.(x^{(1)}_i \otimes \ldots \otimes x^{(p)}_i \otimes y^{(1)}_j \otimes \ldots \otimes y^{(q)}_j)\\ 
&=&\sum_{i,j} \sum_{\tilde{\sigma} \in bat(q,p)} \tilde{\sigma}.(y^{(1)}_j \otimes \ldots \otimes y^{(q)}_j\otimes x^{(1)}_i \otimes \ldots \otimes x^{(p)}_i )\\ 
&=&\tilde{\Delta}^{p+q-1}(yx).
\end{eqnarray*}
Par suite, $deg_p(xy-yx)<p+q=deg_p(x) +deg_p(y)$, et donc l'alg\`ebre gradu\'ee associ\'ee est commutative, donc int\`egre d'apr\`es le r\'esultat pr\'ec\'edent.
On en d\'eduit que $A$ est elle-m\^eme int\`egre, et donc $poids(xy)=poids(x)+poids(y)$, $\forall x,y \in A$. $\Box$\\

{\it Remarque :} cette propri\'et\'e est fausse si $K$ est de caract\'eristique $p$ non nulle. En effet, soit $x$ un \'el\'ement primitif non nul de $A$ ; alors $deg_p(x)=1.$ De plus, $x^p$ est aussi primitif,
donc $deg_p(x^p)=1 \neq p$.  
 
\section{Alg\`ebre de Hopf des arbres enracin\'es plans d\'ecor\'es}

\subsection{Construction}
\label{partie2}
\begin{defi}
\textnormal{ Un {\it arbre plan enracin\'e} $t$ est la donn\'ee d'un graphe fini orient\'e connexe et sans boucles, muni d'un plongement dans le plan ; on suppose que l'un des sommets de ce graphe n'est l'arriv\'ee d'aucune ar\^ete ; ce sommet est appel\'e {\it racine} de $t$. 
Les arbres plans enracin\'es seront dessin\'es avec la racine en bas.
Le $poids$ de $t$ est le nombre de ses sommets. L'ensemble des arbres plans enracin\'es est not\'e ${\cal T}_{P,R}$}.

\textnormal{ Soit ${\cal D}$ un ensemble non vide. Un {\it arbre plan enracin\'e d\'ecor\'e par} ${\cal D}$ est un arbre plan enracin\'e $t$ muni d'une application $d_t$ de l'ensemble de ses sommets vers $\cal D$.
L'image d'un sommet  $s$ par cette application est appel\'ee {\it d\'ecoration} de $s$. L'ensemble des arbres plans enracin\'es d\'ecor\'es par $\cal D$ sera not\'e $\trees$.} 
\textnormal{Pour tout $d \in {\cal D}$, on notera $\bullet_d$ l'\'el\'ement de $\trees$ form\'e d'un seul sommet d\'ecor\'e par $d$.}
\end{defi}  

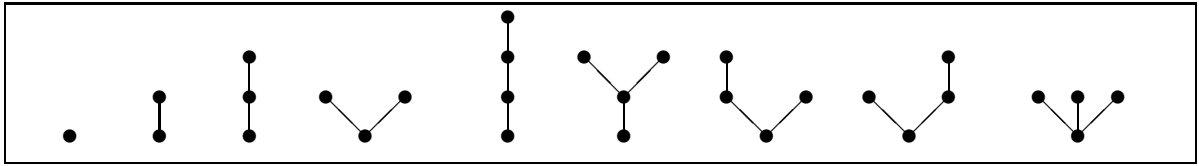
\begin{figure}[h]
\framebox(450,60){
\begin{picture}(10,20)(140,-10)
\circle*{5}
\end{picture}
\begin{picture}(10,40)(120,-10)
\put(0,0){\circle*{5}}
\put(0,0){\line(0,1){15}}
\put(0,15){\circle*{5}}
\end{picture}
\begin{picture}(10,60)(100,-10)
\put(0,0){\circle*{5}}
\put(0,0){\line(0,1){15}}
\put(0,15){\circle*{5}}
\put(0,15){\line(0,1){15}}
\put(0,30){\circle*{5}}
\end{picture}
\begin{picture}(10,60)(70,-10)
\put(0,0){\circle*{5}}
\put(0,0){\line(1,1){15}}
\put(0,0){\line(-1,1){15}}
\put(-15,15){\circle*{5}}
\put(15,15){\circle*{5}}
\end{picture}
\begin{picture}(10,60)(30,-10)
\put(0,0){\circle*{5}}
\put(0,0){\line(0,1){15}}
\put(0,15){\circle*{5}}
\put(0,15){\line(0,1){15}}
\put(0,30){\circle*{5}}
\put(0,45){\circle*{5}}
\put(0,30){\line(0,1){15}}
\end{picture}
\begin{picture}(10,60)(0,-10)
\put(0,0){\circle*{5}}
\put(0,0){\line(0,1){15}}
\put(0,15){\circle*{5}}
\put(0,15){\line(1,1){15}}
\put(0,15){\line(-1,1){15}}
\put(-15,30){\circle*{5}}
\put(15,30){\circle*{5}}
\end{picture}
\begin{picture}(10,60)(-40,-10)
\put(0,0){\circle*{5}}
\put(0,0){\line(1,1){15}}
\put(0,0){\line(-1,1){15}}
\put(-15,15){\circle*{5}}
\put(15,15){\circle*{5}}
\put(-15,15){\line(0,1){15}}
\put(-15,30){\circle*{5}}
\end{picture}
\begin{picture}(10,60)(-80,-10)
\put(0,0){\circle*{5}}
\put(0,0){\line(1,1){15}}
\put(0,0){\line(-1,1){15}}
\put(-15,15){\circle*{5}}
\put(15,15){\circle*{5}}
\put(15,15){\line(0,1){15}}
\put(15,30){\circle*{5}}
\end{picture}
\begin{picture}(10,60)(-130,-10)
\put(0,0){\circle*{5}}
\put(0,0){\line(1,1){15}}
\put(0,0){\line(-1,1){15}}
\put(0,15){\circle*{5}}
\put(-15,15){\circle*{5}}
\put(0,0){\line(0,1){15}}
\put(15,15){\circle*{5}}
\end{picture} }
\caption{\it les arbres plans enracin\'es de poids inf\'erieur o\`u \'egal \`a 4.}
\end{figure}

Soit $\cal D$ un ensemble non vide, et soit $\Hr$ l'alg\`ebre libre engendr\'ee sur $\mathbb{Q} $ par les \'el\'ements de $\trees$. 
Les mon\^omes en les arbres plans enracin\'es de cette alg\`ebre sont appel\'es {\it for\^ets planes enracin\'ees d\'ecor\'ees} ; il sera souvent utile de consid\'erer $1$ comme la for\^et vide.
L'ensemble des for\^ets planes enracin\'ees d\'ecor\'ees est not\'e $\forets$.
Le poids d'une for\^et $F=t_1 \ldots t_n$ est par d\'efinition $poids(t_1)+\ldots +poids(t_n)$.
\\

On va munir $\Hr$ d'une structure d'alg\`ebre de Hopf.
Soit $t \in \trees$. Une {\it coupe \'el\'ementaire} de $t$ est une coupe sur une seule ar\^ete de $t$.
Une {\it coupe admissible} de $t$ est une coupe non vide telle que tout trajet d'un sommet de $t$ vers un autre ne rencontre au plus qu'une seule coupe \'el\'ementaire. 
L'ensemble des coupes admissibles de $t$ est not\'e ${\cal A}d_*(t)$.
Une coupe admissible $c$ envoie $t$ vers un couple $(P^c(t),R^c(t)) \in \forets \times \trees$, tel que $R^c(t)$ est la composante connexe de la racine de $t$ apr\`es la coupe, et $P^c(t)$ est la for\^et plane form\'ee par les autres composantes connexes.

D'autre part, si $c_v$ est la coupe vide de $t$, on pose $P^{c_v}(t)=1$ et $R^{c_v}(t)=t$. 
On d\'efinit la {\it coupe totale } de $t$ comme une coupe $c_t$ telle  $P^{c_t}(t)=t$ et $R^{c_t}(t)=1$.
L'ensemble form\'e des coupes admissibles de $t$, de la coupe vide et de la coupe totale de $t$ est not\'e ${\cal A}d(t)$.\\

Soit maintenant $F \in \forets$, $F \neq 1$. Il existe $t_1, \ldots, t_n \in \trees$, tels que $F=t_1\ldots t_n$.
Une {\it coupe admissible} de $F$ est un $n$-uplet $(c_1,\ldots,c_n)$ tel que $c_i \in {\cal A}d(t_i) \esp \forall i$.
Si toutes les $c_i$ sont vides (respectivement totales), $c$ est appel\'ee la coupe vide de $F$ (respectivement la coupe totale de $F$). L'ensemble des coupes admissibles non vides et non totales de $F$ est not\'e ${\cal A}d_*(F)$.
L'ensemble de toutes les coupes admissibles de $F$ est not\'e ${\cal A}d(F)$.
Pour $c=(c_1,\ldots,c_n) \in {\cal A}d(F)$, on pose $P^c(F)= P^{c_1}(t_1) \ldots P^{c_n}(t_n)$ et  $R^c(F)= R^{c_1}(t_1) \ldots R^{c_n}(t_n)$.\\

On d\'efinit $\Delta : \Hr \longmapsto \Hr \otimes \Hr$ par :
\begin{eqnarray}
\nonumber
\Delta(1)&=&1 \otimes 1, \\
\nonumber
\Delta(F) &=& \sum_{c \in \ad(F)} P^c(F) \otimes R^c(F) \\
\label{defDelta}
&=& 1 \otimes F + F \otimes 1 + \sum_{c \in \ad_*(F)} P^c(F) \otimes R^c(F),
\mbox{ pour $F \in \forets$, $F\neq 1$.}
\end{eqnarray}

\begin{figure}[h]
{\it coupes admissibles : }
\begin{picture}(27,30)(-47,0)
\put(-3,-1){\circle*{3}}
\put(-3,-1){\line(0,1){12}}
\put(-3,11){\circle*{3}}
\put(-3,11){\line(-1,1){10}}
\put(-3,11){\line(1,1){10}}
\put(-13,21){\circle*{3}}
\put(7,21){\circle*{3}}
\put(-15,12){\line(1,1){10}}
\put(1,-4){$a$}
\put(1,6){$b$}
\put(5,24){$d$}
\put(-15,24){$c$}
\end{picture}
\begin{picture}(27,30)(-72,0)
\put(-3,-1){\circle*{3}}
\put(-3,-1){\line(0,1){12}}
\put(-3,11){\circle*{3}}
\put(-3,11){\line(-1,1){10}}
\put(-3,11){\line(1,1){10}}
\put(-13,21){\circle*{3}}
\put(7,21){\circle*{3}}
\put(-2,21){\line(1,-1){10}}
\put(1,-4){$a$}
\put(1,6){$b$}
\put(5,24){$d$}
\put(-15,24){$c$}
\end{picture}
\begin{picture}(27,30)(-170,0)
\put(-3,-1){\circle*{3}}
\put(-3,-1){\line(0,1){12}}
\put(-3,11){\circle*{3}}
\put(-3,11){\line(-1,1){10}}
\put(-3,11){\line(1,1){10}}
\put(-13,21){\circle*{3}}
\put(7,21){\circle*{3}}
\put(-9,5){\line(1,0){10}}
\put(1,-4){$a$}
\put(1,6){$b$}
\put(5,24){$d$}
\put(-15,24){$c$}
\end{picture}
\begin{picture}(27,30)(-70,0)
\put(-3,-1){\circle*{3}}
\put(-3,-1){\line(0,1){12}}
\put(-3,11){\circle*{3}}
\put(-3,11){\line(-1,1){10}}
\put(-3,11){\line(1,1){10}}
\put(-13,21){\circle*{3}}
\put(7,21){\circle*{3}}
\put(-15,12){\line(1,1){10}}
\put(-2,21){\line(1,-1){10}}
\put(1,-4){$a$}
\put(1,6){$b$}
\put(5,24){$d$}
\put(-15,24){$c$}
\end{picture}

 $\Delta($
\begin{picture}(27,30)(-12,11)
\put(-3,-1){\circle*{3}}
\put(-3,-1){\line(0,1){12}}
\put(-3,11){\circle*{3}}
\put(-3,11){\line(-1,1){10}}
\put(-3,11){\line(1,1){10}}
\put(-13,21){\circle*{3}}
\put(7,21){\circle*{3}}
\put(1,-4){$a$}
\put(1,6){$b$}
\put(5,24){$d$}
\put(-15,24){$c$}
\end{picture}
)=  
\begin{picture}(27,30)(-12,11)
\circle*{3}
\put(-3,-1){\line(0,1){12}}
\put(-3,11){\circle*{3}}
\put(-3,11){\line(-1,1){10}}
\put(-3,11){\line(1,1){10}}
\put(-13,21){\circle*{3}}
\put(7,21){\circle*{3}}
\put(1,-4){$a$}
\put(1,6){$b$}
\put(5,24){$d$}
\put(-15,24){$c$}
\end{picture}
$\otimes$ 1 +
\begin{picture}(10,30)(-7,0)
\circle*{3}
\put(-5,4){$c$}
\end{picture}
$\otimes$
\begin{picture}(10,40)(-2,5)
\put(-2,-1){\circle*{3}}
\put(-2,-1){\line(0,1){12}}
\put(-2,11){\circle*{3}}
\put(-2,11){\line(0,1){12}}
\put(-2,23){\circle*{3}}
\put(1,-4){$a$}
\put(1,6){$b$}
\put(1,20){$d$}
\end{picture}
+
\begin{picture}(10,30)(-7,0)
\circle*{3}
\put(-5,4){$d$}
\end{picture}
$\otimes$
\begin{picture}(12,40)(-2,5)
\put(-2,-1){\circle*{3}}
\put(-2,-1){\line(0,1){12}}
\put(-2,11){\circle*{3}}
\put(-2,11){\line(0,1){12}}
\put(-2,23){\circle*{3}}
\put(1,-4){$a$}
\put(1,6){$b$}
\put(1,20){$c$}
\end{picture}
+ 
\begin{picture}(14,30)
\circle*{3}
\put(7,0){\circle*{3}}
\put(5,4){$d$}
\put(-5,4){$c$}
\end{picture}
$\otimes$ 
\begin{picture}(10,30)(-3,2)
\circle*{3}
\put(-3,0){\line(0,1){12}}
\put(-3,12){\circle*{3}}
\put(0,-3){$a$}
\put(0,8){$b$}
\end{picture}
+
\begin{picture}(27,30)(-18,11)
\put(-3,11){\circle*{3}}
\put(-3,11){\line(-1,1){10}}
\put(-3,11){\line(1,1){10}}
\put(-13,21){\circle*{3}}
\put(7,21){\circle*{3}}
\put(1,6){$b$}
\put(5,24){$d$}
\put(-15,24){$c$}
\end{picture}
$\otimes$
\begin{picture}(10,30)(-2,0)
\circle*{3}
\put(0,-3){$a$}
\end{picture}
+ 1 $\otimes$
\begin{picture}(27,30)(-12,11)
\put(-3,-1){\circle*{3}}
\put(-3,-1){\line(0,1){12}}
\put(-3,11){\circle*{3}}
\put(-3,11){\line(-1,1){10}}
\put(-3,11){\line(1,1){10}}
\put(-13,21){\circle*{3}}
\put(7,21){\circle*{3}}
\put(1,-4){$a$}
\put(1,6){$b$}
\put(5,24){$d$}
\put(-15,24){$c$}
\end{picture}

\begin{picture}(-510,-200)
\put(-5,-5){\line(0,1){100}}
\put(-5,95){\line(1,0){450}}
\put(-5,-5){\line(1,0){450}}
\put(445,-5){\line(0,1){100}}
\end{picture}
\caption{\it calcul d'un coproduit dans $\Hr$, avec ${\cal D}=\{a,b,c,d\}$. }
\end{figure}

Soit $F=t_1\ldots t_n \in \forets$. On d\'eduit imm\'ediatement de la d\'efinition des coupes admissibles de $F$ que $\Delta(F)=\Delta(t_1) \ldots \Delta(t_n)$, et donc $\Delta$ est un morphisme d'alg\`ebres.\\

Soit $\varepsilon : \left\{ \begin{array}{rcl}
                     \Hr &\longmapsto &\mathbb{Q}\\
                       1 & \longmapsto & 1\\
                       F & \longmapsto & 0 \mbox{ si $F\in \forets$, $F\neq 1$.}
\end{array} \right. $\\
Il est imm\'ediat que $\varepsilon$ est une counit\'e pour $\Delta$.\\

Pour montrer que $(\Hr,m,\eta,\Delta,\varepsilon)$ est une big\`ebre, il ne reste plus qu'\`a montrer que $\Delta$ est coassociatif.
Pour cela, on introduit les op\'erateurs suivants : pour $d \in {\cal D}$, on pose $B_d^+ :\forets \longmapsto \trees$, tel que $B_d^+(t_1 \ldots t_n)$ est l'arbre plan enracin\'e d\'ecor\'e obtenu en reliant les racines de $t_1,\ldots, t_n$ (dans cet ordre) \`a une racine commune d\'ecor\'ee par $d$.
En particulier, $B_d^+(1)=\bullet_d$. On prolonge $B_d^+$ en un op\'erateur de $\Hr$.

On d\'efinit \'egalement $B^- :\trees \longmapsto \forets$, qui envoie un arbre enracin\'e plan d\'ecor\'e $t$ sur la for\^et obtenue en \^otant la racine de $t$. On a $B^- \circ B^+_d(F)=F$ $ \forall F \in \forets$, $\forall d \in \cal D$.  

\begin{prop}
\label{pro7}
$\forall x \in \Hr, \esp \Delta\left(B_d^+(x)\right)=B_d^+(x) \otimes 1 + (Id \otimes B_d^+)\circ \Delta(x)$.
\end{prop}
{\it Preuve :} il suffit de le montrer pour $F=t_1\ldots t_n \in \forets$.

Dans ce cas, on a une bijection $\alpha :\ad(F) \longmapsto \ad(B_d^+(F))-\{\mbox{coupe totale de $B_d^+(F)$}\}$, qui envoie la coupe $c$ de $F$ sur la coupe $\alpha(c)$ de $B_d^+(F)$ 
telle que  $P^{\alpha(c)}(B_d^+(F))=P^c(F)$ et $R^{\alpha(c)}(B_d^+(F))=B_d^+(P^c(F))$. Le r\'esultat d\'ecoule alors imm\'ediatement de (\ref{defDelta}). $\Box$\\

Montrons que $\Delta$ est coassociatif. 

Il suffit de montrer que $(\Delta \otimes Id)\circ \Delta(F) =(Id \otimes \Delta)\circ \Delta(F) $, $\forall F \in \forets$.
Proc\'edons par r\'ecurrence sur $poids(F)$.
Si $poids(F)=0$ : alors $F=1$, et le r\'esultat est \'evident.
Supposons la propri\'et\'e vraie pour toute for\^et de poids inf\'erieur ou \'egal \`a $n-1$, et soit $F$ une for\^et de poids $n$. Deux cas se pr\'esentent :
\begin{enumerate}
\item $F \notin \trees$ : alors il existe deux for\^ets non vides $F_1$ et $F_2$ telles que $F=F_1F_2$. On a $poids(F_1)< n$ et $poids(F_2)< n$, donc la propri\'et\'e est vraie pour $F_1$ et $F_2$,
et comme $\Delta$ est un morphisme d'alg\`ebres, elle est vraie pour $F$.
\item $F \in \trees$ : alors il existe un unique $d \in {\cal D}$ et une unique $F_1 \in \forets$, $F=B_d^+(F_1)$. De plus $poids(F_1)=poids(F)-1$, donc la propri\'et\'e est vraie pour $F_1$. On a :
\end{enumerate}
\begin{eqnarray*}
\Delta(F)&=&\Delta\left(B_d^+(F_1)\right)= F\otimes 1 + (Id \otimes B_d^+)\circ\Delta(F_1).
\end{eqnarray*}
D'o\`u :
\begin{eqnarray*}
 \left(Id\otimes \Delta \right)\circ\Delta(F)&=& F \otimes 1 \otimes 1 + \left(Id \otimes (\Delta \circ B_d^+)\right)\circ \Delta(F_1)\\
&=& F \otimes 1 \otimes 1 + \left(Id \otimes Id \otimes B_d^+\right)\circ \left( Id \otimes \Delta\right)\circ\Delta(F_1)\\
&&+ \left((Id \otimes B_d^+)\circ \Delta(F_1)\right) \otimes 1\\
&=& F \otimes 1 \otimes 1 + \left(Id \otimes Id \otimes B_d^+\right)\circ\left( \Delta \otimes Id\right)\circ\Delta(F_1)\\
&&+ \left((Id \otimes B_d^+)\circ\Delta(F_1)\right) \otimes 1\\
&=& F \otimes 1 \otimes 1 + \left(\Delta \otimes B_d^+\right)\circ\Delta(F_1)\\
&&+ \left((Id \otimes B_d^+)\circ\Delta(F_1)\right) \otimes 1\\
&=&(\Delta \otimes Id)\left( F \otimes 1 + (Id \otimes B_d^+)\circ\Delta(F_1)\right)\\
&=&\left(\Delta \otimes Id \right)\circ\Delta(F).
\end{eqnarray*}
(On a utilis\'e l'hypoth\`ese de r\'ecurrence pour la troisi\`eme \'egalit\'e, et la proposition \ref{pro7} pour la premi\`ere, la deuxi\`eme, la cinqui\`eme et la sixi\`eme \'egalit\'e.)\\

Montrons que $\Hr$ poss\`ede une antipode. 

Pour $t \in \trees$, posons $\Delta(t)=t\otimes 1+1\otimes t +\sum t' \otimes t''$, $poids(t')<poids(t),$ $poids(t'')<poids(t)$. Soient $S_1$ et $S_2$ les antimorphismes d'alg\`ebres d\'efinis par r\'ecurrence sur le poids de la mani\`ere suivante :
$$S_1(t)=-t -\sum S_1(t')t'', \spa S_2(t)=-t- \sum t'S_2(t'').$$
Alors pour tout $t \in\trees$, on a :
\begin{eqnarray*}
m\circ (S_ 1 \otimes Id)\circ \Delta(t)&=&0\\
&=& \eta \circ \varepsilon(t) \esp ;\\
m\circ (Id \otimes S_2)\circ \Delta(t)&=&0\\
&=& \eta \circ \varepsilon(t).
\end{eqnarray*}
Donc dans ${\cal L}(\Hr)$ muni du produit de convolution $\star$, on a :
$$ S_1 \star Id=Id \star S_2= 1_{{\cal L}(\Hr)}.$$
Par suite, on a :
\begin{eqnarray*}
(S_1 \star Id)\star S_2&=&1 \star S_2\\
&=&S_2\\
&=&S_1 \star (Id \star S_2)\\
&=&S_1 \star 1\\
&=&S_1.
\end{eqnarray*}
Donc $S=S_1=S_2$ est une antipode de $\Hr$.

\begin{theo}
$(\Hr,m,\eta,\Delta,\varepsilon,S)$ est une alg\`ebre de Hopf.
\end{theo}
On donnera dans le th\'eor\`eme \ref{theo37} une expression directe de l'antipode.\\ 

{\it Remarques :} 
\label{rem}
\begin{enumerate}
\item Si $\cal D$ et $\cal D'$ ont le m\^eme cardinal, il est \'evident que $\Hr$ et ${\cal H}_{P,R}^{\cal D'}$ sont isomorphes comme alg\`ebres de Hopf.
\item On peut \'egalement effectuer cette construction \`a partir de $ {\cal T}_{P,R}$, l'ensemble des arbres enracin\'es plans (non d\'ecor\'es) ;
l'alg\`ebre de Hopf obtenue est not\'ee ${\cal H}_{P,R}$. 
Quand le cardinal de $\cal D$ est \'egal \`a 1, on a un isomorphisme d'alg\`ebres de Hopf entre $\Hr$ et ${\cal H}_{P,R}$ qui envoie l'arbre d\'ecor\'e $(t,d_t)$ sur $t$.
\end{enumerate}

\subsection{Graduation}
\spa Il est clair que le poids d\'efinit une graduation de l'alg\`ebre de Hopf $\Hr$. De plus, si $\cal D$ est fini, de cardinal $D$,
les composantes homog\`enes sont de dimension finie. On se propose de calculer ces dimensions. 
On note ${\cal H}_n$ la composante homog\`ene de poids $n$, et $r_n$ sa dimension. Soit $R(X)=\sum r_n X^n$ la s\'erie g\'en\'eratrice des $r_n$.

Soit $\tau_k$ le nombre d'arbres plans enracin\'es de poids $k$. Le nombre d'arbres plans enracin\'es d\'ecor\'es par $\cal D$ de poids $k$ est $D^k \tau_k$. On note $T(X)=\sum \tau_k X^k$ la s\'erie g\'en\'eratrice des $\tau_k$.\\
Une base de ${\cal H}_n$ est $(t_1 \ldots t_m)_{t_1, \ldots, t_m \in \trees, \esp poids(t_1)+\ldots +poids(t_m)=n}$. D'o\`u :
\begin{eqnarray*}
r_n&=&\sum_{m=1}^n\esp \sum_{\stackrel{ a_1 + \ldots + a_m=n}{\mbox{\scriptsize{les $a_i$ tous non nuls}}}}
    D^{a_1+\ldots +a_m} \tau_{a_1} \ldots \tau_{a_m} \\
&=& D^n \sum_{b_1+2b_2+\ldots+nb_n=n} \frac{(b_1+\ldots+b_n)!}{b_1!\ldots b_n!} \tau_1^{b_1} \ldots \tau_n^{b_n}.\\
\end{eqnarray*}
($b_i$ est le nombre de $a_j$ \'egal \`a $i$.)\\
On a donc :
\begin{eqnarray*}
R(X)&=&\sum_{n=0}^{+\infty} \esp 
\sum_{b_1+2b_2+\ldots+nb_n=n} \frac{(b_1+\ldots+b_n)!}{b_1!\ldots b_n!} (\tau_1D^1X^1)^{b_1} \ldots (\tau_nD^nX^n)^{b_n}\\
&=&\sum_{n=0}^{+\infty} \left(\sum_{k=0}^{+\infty} \tau_k (DX)^k \right)^n.
\end{eqnarray*}
D'o\`u :
\begin{eqnarray}
\label{eqnR}
R(X)=\frac{1}{1-T(DX)}.
\end{eqnarray}
De plus, $\trees= \bigcup_{d \in \cal D} B_d^+(\forets)$ (union disjointe).
On en d\'eduit : $D^k \tau_k=Dr_{k-1}$, et donc :
\begin{eqnarray}
R(X)=\frac{1}{DX}T(DX).
\end{eqnarray}
Puis en combinant ces deux r\'esultats, et en prenant la valeur particuli\`ere $D=1$ :
\begin{eqnarray}
T(X)^2-T(X)-X=0.
\end{eqnarray}
On en d\'eduit imm\'ediatement le th\'eor\`eme suivant :
\begin{theo}
\label{the9}
\begin{enumerate}
\item $$ T(X)=\frac{1-\sqrt{1-4X}}{2} \esp ; $$
$$ \tau_k = \frac{(2k-2)}{k!(k-1)!}.$$
\item $$ R(X)=\frac{1-\sqrt{1-4DX}}{2DX} \esp ;$$
$$ r_n=\frac{(2n)!}{(n+1)!n!}D^n.$$
\end{enumerate}
\end{theo}   
Voir la section \ref{part81} pour les valeurs de $\tau_k$, $k \leq 24$.\\

{\it Remarque :} $\tau_k$ est donc \'egal au $k$-i\`eme nombre de Catalan ; on a donc la relation :
\begin{eqnarray}
  \tau_k&=&\sum_{i=1}^{k-1} \tau_i\tau_{k-i}\esp \forall k \geq 2.
\end{eqnarray}
On peut donner une preuve directe de ce r\'esultat en consid\'erant la bijection :    
\begin{eqnarray*}
\trees &\longmapsto &\trees \times \trees\\
B^+(t_1\ldots t_n)&\longmapsto &(B^+(t_1\ldots t_{n-1}),t_n).
\end{eqnarray*}

\subsection{Propri\'et\'e universelle}
\spa Soit $A$ une big\`ebre. Suivant \cite{Connes}, On appelera 1-cocycle de $A$ toute application lin\'eaire $L$ de $A$ dans $A$ v\'erifiant :
\begin{eqnarray}
\label{eqn9}
\Delta_A\left(L(x)\right)&=&L(x) \otimes 1 + (Id \otimes L)\circ \Delta_A(x),  \esp \forall x \in A.
\end{eqnarray}
Par exemple, pour tout $d\in \cal D$, $B_d^+$ est un 1-cocycle de $\Hr$. On notera $\cocy(A^*)$ l'espace des 1-cocycles de $A$.

\begin{theo}
\label{theo10}
\begin{enumerate}
\item Soit $A$ une big\`ebre, $\cal D$ un ensemble non vide, et $L_d \in \cocy(A^*)$ pour tout $d \in \cal D$.
Alors il existe un unique morphisme de big\`ebres  $\varphi$ de $\Hr$ dans $A$ v\'erifiant :
$$ \varphi \circ B_d^+=L_d \circ \varphi \esp \forall d \in \cal D.$$
Si de plus $A$ est une alg\`ebre de Hopf, alors $\varphi$ est un morphisme d'alg\`ebres de Hopf.
\item Cette propri\'et\'e caract\'erise $\Hr$, c'est-\`a-dire : si ${\cal H}$ est une alg\`ebre de Hopf, $b^+_d \in \cocy({\cal H}^*)$  v\'erifiant 1, alors il existe un isomorphisme d'alg\`ebres de Hopf $\Psi$ de $\Hr$ dans $\cal H$ 
tel que $\Psi \circ B_d^+=b^+_d \circ \Psi$, $ \forall d \in \cal D$.
\end{enumerate}
\end{theo}
{\it Preuve} :

1. Unicit\'e : on doit avoir $\varphi(1)=1$ ; de plus, pour tout $t\in \trees$, dont la racine est d\'ecor\'ee par $d$ :
\begin{eqnarray}
\nonumber \varphi(t)&=& \varphi\left(B_d^+\circ B^-(t)\right)\\
                   &=& L_d \circ \varphi(B^-(t)).
\label{equphi}
\end{eqnarray}
Comme $poids(B^-(t))=poids(t)-1$, et comme $\trees$ engendre $\Hr$, une r\'ecurrence montre qu'il existe au plus un morphisme d'alg\`ebres v\'erifiant (\ref{equphi}).\\

Existence : $\trees$ engendre librement $\Hr$, donc il existe un unique morphisme d'alg\`ebres de $\Hr$ dans $A$ v\'erifiant (\ref{equphi}).
Il s'agit de montrer que c'est aussi un morphisme de cog\`ebres, c'est-\`a-dire $\Delta_A\left(\varphi(t)\right)= \varphi \otimes \varphi \left(\Delta(t)\right) \esp \forall t\in \trees.$

Proc\'edons par r\'ecurrence sur $poids(t)$ : si $poids(t)=1$, alors il existe $ d \in \cal D$, tel que $t=\bullet_d$.
\begin{eqnarray*}
\mbox{Alors }\Delta_A\left(\varphi(\bullet_d)\right)&=&\Delta_A\left(\L_d \circ \varphi(1)\right)\\
                                          &=&\Delta_A\left( L_d(1)\right)\\
                                          &=&L_d(1) \otimes 1 + 1 \otimes L_d(1)\\
                                          &=&\varphi \otimes \varphi \left(\Delta(\bullet_d)\right).  
\end{eqnarray*}
Supposons l'hypoth\`ese vraie pour tout $t'$ de poids $<n$ et soit $t$ de poids $n$. Soit $d$ la d\'ecoration de la racine de $t$ et soit $F=B^-(t)$ ; alors $t=B_d^+(F)$. 
On a :
\begin{eqnarray*}
\Delta_A\left(\varphi(t)\right)&=&\Delta_A\left(\L_d\circ \varphi(F)\right)\\
                                          &=&L_d(\varphi(F)) \otimes 1 + (Id\otimes L_d) \circ\Delta_A(\varphi(F)) \\
                                          &=&\varphi(t) \otimes 1 + (Id\otimes L_d) \circ ( \varphi \otimes \varphi)\circ \Delta(F)\\
                                          &=&(\varphi \otimes \varphi) \left(t\otimes 1 + (Id\otimes B^+_d)\circ \Delta(F)\right)\\ 
                                           &=&\varphi \otimes \varphi \left(\Delta(t)\right).  
 \end{eqnarray*}
(On a utilis\'e le fait que $L_d$ soit un 1-cocycle pour la deuxi\`eme \'egalit\'e, l'hypoth\`ese de r\'ecurrence pour la troisi\`eme \'egalit\'e,
l'\'equation (\ref{equphi}) pour la troisi\`eme et la quatri\`eme \'egalit\'e, et le fait que $B_d^+$ soit un 1-cocycle pour la
derni\`ere \'egalit\'e.)\\

Montrons que $\varepsilon_A \circ \varphi(t)=\varepsilon(t)=0$, $\forall t \in \trees$. Il suffit de montrer que $\varepsilon_A \circ L_d=0$, $\forall d \in \cal D$.
Soit $x \in A$. On pose $\Delta_A(x)=\sum x'\otimes x''$. On a :
\begin{eqnarray*}
(\varepsilon_A\otimes Id)\circ \Delta_A(L_d(x))&=&(\varepsilon_A\otimes Id)\left(L_d(x)\otimes 1 +Id\otimes L_d(\sum x'\otimes x'')\right)\\
&=&\varepsilon_A(L_d(x))1+Id \otimes L_d(\sum  \varepsilon_A(x')\otimes x'')\\
&=&\varepsilon_A(L_d(x))1+ L_d(\sum \varepsilon_A(x') x'')\\
&=&\varepsilon_A(L_d(x))1+ L_d(x)\\
&=& L_d(x).
\end{eqnarray*}
Donc $\varepsilon_A \circ L_d=0$.\\

Supposons que $A$ a une antipode $S_A$. Montrons que $\varphi \circ S(t)=S_A \circ\varphi(t)$ pour tout $t  \in \trees$ par r\'ecurrence sur $n=poids(t)$.
Si $poids(t)=1$, alors $t$ est primitif, et donc $\varphi(t)$ est primitif,  d'o\`u :
$$\varphi \circ S(t)=-\varphi(t)=-S_A\circ \varphi(t).$$
Supposons la propri\'et\'e vraie pour tout arbre de poids inf\'erieur \`a $n$. Posons $t=B_d^+(F)$, $\Delta(F)=1 \otimes F+F\otimes 1 +\sum F' \otimes F''$. On a :
\begin{eqnarray*}
m \circ (S\otimes Id)\circ \Delta(t)&=&m\circ (S\otimes Id)( t \otimes 1 +1 \otimes t + F \otimes \bullet_d +\sum F' \otimes B_d^+(F''))\\
&=&S(t)+t+S(F) \bullet_d +\sum S(F')B_d^+(F'')\\
&=&\varepsilon(t)1=0.
\end{eqnarray*}
Et donc :
 \begin{eqnarray*}
\varphi \circ S(t)&=&-\varphi(t) -\varphi \circ S (F) \varphi(\bullet_d) -\sum \varphi\circ S(F') \varphi(B^+_d(F''))\\
&=& -\varphi(t) -\varphi \circ S (F) L_d(1) -\sum \varphi \circ S(F')  L_d (\varphi(F''))\\
&=& -\varphi(t) -S_A \circ \varphi(F) L_d(1) -\sum S_A ( \varphi(F)')  L_d(\varphi(F)'').
\end{eqnarray*}
(On a utilis\'e le fait que $\varphi$ est un morphisme de big\`ebres v\'erifiant $L_d \circ \varphi=\varphi \circ B_d^+$, et l'hypoth\`ese de r\'ecurrence.)

De plus :
\begin{eqnarray*}
m \circ (S_A \otimes Id) \circ \Delta_A(\varphi(t))&=&
S_A(\varphi(t)) + \varphi(t) +S_A\circ \varphi(F) L_d(1)\\
&& + \sum S_A (\varphi(F)' )L_d(\varphi(F)'')\\
&=&\varepsilon_A\circ \varphi(t)1\\
&=&\varepsilon(t)1\\
&=&0.
\end{eqnarray*}
 (Car $L_d\in \cocy(A^*)$.)\\
Donc on a bien $\varphi \circ S(t)=S_A \circ \varphi (t)$.
\\

2. Soit $\Psi :\Hr\longmapsto \cal H$ l'unique morphisme d'alg\`ebres de Hopf tel que $\Psi \circ B_d^+=b^+_d \circ \Psi$, et soit $\Psi' :{\cal H}\longmapsto \Hr$ l'unique morphisme d'alg\`ebres de Hopf v\'erifiant 
$\Psi'\circ b^+_d=B_d^+ \circ \Psi'$. On a alors :
\begin{eqnarray*}
\Psi' \circ \Psi  \circ B_d^+&=&\Psi' \circ b^+_d \circ \Psi\\
&=& B_d^+ \circ \Psi'\circ \Psi.
\end{eqnarray*}
Donc $\Psi'\circ \Psi$ est l'unique morphisme d'alg\`ebres de Hopf de $\Hr$ dans $\Hr$ commutant avec $B_d^+$ pour tout $d$ : c'est donc l'identit\'e de $\Hr$.
De m\^ eme, $\Psi \circ \Psi'$ est l'identit\'e de $\cal H$, et donc $\Psi $ et $\Psi'$ sont des isomorphismes r\'eciproques l'un de l'autre. $\Box$

\section{Dual gradu\'e de $\Hr$}
\subsection{Construction de la forme bilin\'eaire $(\hspace{1mm},)$}
\label{part55}
\spa Dans toute cette section, on suppose $\cal D$ fini de cardinal $D$.
Comme nous l'avons remarqu\'e \`a la fin de la partie \ref{rem} (remarque 1), on peut supposer que ${\cal D}=\{1, \ldots ,D\}$.

\begin{lemme}
\label{lemme18}
Soit $p=\sum_{F \in \forets} a_F F$ un primitif non nul de $\Hr$. Alors il existe $F \in \forets$, de la forme $t_1\ldots t_{n-1} \bullet_d$, $d \in \cal D$, telle que
$a_F$ soit non nul.
\end{lemme}
{\it Preuve :} soit $F=t_1 \ldots t_n \in \forets$ telle que :
\begin{description}
\item[\textnormal{a)}] $a_F \neq 0$,
\item[\textnormal{b)}] si $G=t_1'\ldots t'_m \in \forets$ est telle que que $a_G \neq 0$, alors $m \leq n$, et si $n=m$, $poids(t'_n)\geq poids(t_n)$. 
\end{description}

On note $d$ la d\'ecoration de la racine de $t_n$.
Soit $G=t'_1 \ldots t'_m \in \forets$. On suppose que $\ad(G)$ contient  une coupe $c=(c_1,\ldots, c_m)$ telle que $R^c(G)= t_1 \ldots t_{n-1} \bullet_d$ et $P^c(G)=B^-(t_n)$.
Remarquons que n\'ecessairement, $m\geq n$.
Trois cas se pr\'esentent : \begin{enumerate}
\item soit $m>n$ : dans ce cas $a_G=0$ ; 
\item soit $m=n$ et l'une des $c_i$, $i<n$, n'est pas vide : alors $poids(t'_n)<poids(t_n)$, et donc $a_G=0$ ;
\item soit $m=n$ et toutes les $c_i$, $i<n$, sont vides : alors $G=F$, et alors $c$ est unique.
\end{enumerate}
On en d\'eduit que dans la base $(F_1 \otimes F_2)_{F_i \in \forets}$ de $\Hr \otimes \Hr$, le coefficient dans l'\'ecriture de $\Delta(p)$ de $B^-(t_n) \otimes t_1 \ldots t_{n-1} \bullet_d$ est $a_F$.
Comme $p$ est primitif, n\'ecessairement $B^-(t_n)=1$, d'o\`u $F=t_1 \ldots t_{n-1} \bullet_d$. $\Box$ \\

Soit $d \in \cal D$. On d\'efinit $\gamma_d : \Hr \longmapsto \Hr$ par :
\begin{eqnarray}
\nonumber
\gamma_d(1)&=&0,\\
\label{eq16}
\gamma_d(t_1 \ldots t_n)&=& \left\{ \begin{array}{ccl}
                       0 & \mbox{si} & t_n \neq \bullet_d,\\
                       t_1 \ldots t_{n-1} &  \mbox{si} & t_n =\bullet_d.
\end{array} \right.
\end{eqnarray}
On pose $\overline{\gamma_d}={\gamma_d}_{\mid Prim(\Hr)}.$

\begin{prop}
\label{pro19}
\begin{enumerate}
\item $\gamma_d$ est sujective, homog\`ene de degr\'e $-1$.
\item $ \forall x,y \in \Hr, \esp \gamma_d(xy)=\gamma_d(x) \varepsilon(y) + x \gamma_d(y) \esp ;$
\item $ \forall p \in Prim(\Hr),\esp p\neq 0, \esp \exists d\in {\cal D}, \gamma_d(p) \neq 0.$
\end{enumerate}
\end{prop}
{\it Preuve :} 

1.  D\'ecoule imm\'ediatement de (\ref{eq16}). 

2. On remarque que $\gamma_d(t_1\ldots t_n)=t_1\ldots t_{n-1} \gamma_d(t_n)$. Le r\'esultat en d\'ecoule aussit\^ot.

3. D\'ecoule imm\'ediatement du lemme \ref{lemme18}. $\Box$\\

Soit $M$ l'id\'eal d'augmentation de $(\Hr)^{*g}$.
D'apr\`es la proposition \ref{prop6}, l'orthogonal de $(1) \oplus M^2$ est $Prim(\Hr)$. La dualit\'e entre $\Hr$ et $(\Hr)^{*g}$ induit donc une dualit\'e
entre $Prim(\Hr)$  et $(\Hr)^{*g}/\left((1)\oplus M^2\right)$.\\
$\overline{\gamma_d}$ est homog\`ene de degr\'e $-1$ ; par suite, $\overline{\gamma_d}^{\esp *g} :(\Hr)^{*g} \longmapsto  (\Hr)^{*g}/\left((1)\oplus M^2\right)$ est homog\`ene de degr\'e $1$.\\
Soit $\overline{\gamma} : \left \{ \begin{array}{ccl}
                 Prim(\Hr) &\longmapsto & \left(\Hr\right)^D\\
                    p & \longmapsto & (\overline{\gamma}_1(p), \ldots, \overline{\gamma}_d(p)).
                    \end{array} \right.$\\
Par la proposition \ref{pro19}-3, $\overline{\gamma}$ est injective, donc 
$\overline{\gamma}^{\esp *g} : \left( (\Hr)^{*g} \right)^D \longmapsto (\Hr)^{*g}/\left((1) \oplus M^2\right)$ est surjective.\\
Or $\forall (l_1, \ldots, l_D) \in \left( (\Hr)^{*g} \right)^D,$ $\forall p \in Prim(\Hr)$ :
\begin{eqnarray*}
\overline{\gamma}^{\esp *g}(l_1, \ldots, l_D)(p)&=& \sum_{d=1}^{D} l_d\left((\overline{\gamma_d}(p)\right)\\ 
&=&\sum_{d=1}^{D} \overline{\gamma_d}^{\esp *g}(l_d)(p).
\end{eqnarray*}
Donc $$\overline{\gamma}^{\esp *g}(l_1, \ldots, l_D)= \sum_{d=1}^{D} \overline{\gamma_d}^{\esp *g}(l_d),$$ d'o\`u :
\begin{eqnarray}
\label{eq12}
\nonumber Im(\overline{\gamma}^{\esp *g}) &=&\sum_{d=1}^{D} Im(\overline{\gamma_d}^{\esp *g})\\
&=& \frac{(\Hr)^{*g}}{(1) \oplus M^2}.
\end{eqnarray}

Soit $i$ l'injection canonique de $Prim(\Hr)$ dans $\Hr$ et soit $\pi$ la surjection canonique de $(\Hr)^{*g}$ sur  $(\Hr)^{*g}/\left((1)\oplus M^2\right)$. 
On a $i^{*g}=\pi$ ; comme $\overline{\gamma_d}= \gamma_d \circ i$, on a $\overline{\gamma_d}^{\esp *g} = \pi \circ \gamma_d^{*g}$ ; (\ref{eq12}) implique donc :
\begin{eqnarray}
\label{eq13}
\sum_{d=1}^{D} Im(\gamma_d^{*g}) + \left( (1)\oplus M^2 \right)&=& (\Hr)^{*g}.
\end{eqnarray}

Soient $f \in (\Hr)^{*g}$, $x,y \in  \Hr$.
\begin{eqnarray*}
\Delta\left(\gamma_d^{*g}(f)\right)(x \otimes y) &=& \gamma_d^{*g}(f)(xy)\\
&=& f \circ \gamma_d(xy)\\
&=& f (x \gamma_d(y))+f(\gamma_d(x)) \varepsilon(y)\\
&=& \left( (Id \otimes \gamma_d^{*g})\circ \Delta(f) + \gamma_d^{*g} (f) \otimes 1\right)(x \otimes y).
\end{eqnarray*}
(On a utilis\'e la proposition \ref{pro19}-2 pour la troisi\`eme \'egalit\'e).\\

Donc $\gamma_d^{*g}$ est un 1-cocycle de $(\Hr)^{*g}$. D'apr\`es le th\'eor\`eme \ref{theo10}, il existe un morphisme d'alg\`ebres de Hopf
$ \Psi : \Hr \longmapsto (\Hr)^{*g}$
v\'erifiant $\Psi \circ B_d^+= \gamma_d^{*g} \circ \Psi$.\\

Montrons que $\Psi$ est homog\`ene de poids $0$ : soit $F \in \forets$ de poids $n$, il faut montrer que $\Psi(F) \in {\cal H}_n^*$. Proc\'edons par r\'ecurrence sur $n$ : 
si $n=0$, alors $F=1$, $\Psi(1)=1$. Supposons la propri\'et\'e   vraie pour toute for\^ et de poids strictement inf\'erieur \`a $n$ ;
 comme $\Psi$ est un morphisme d'alg\`ebres on peut supposer $F \in \trees$. Alors $F=B_d^+(G)$, avec 
$poids(G)=n-1$ ; donc $\Psi(F)=\Psi \circ B_d^+(G)=\gamma_d^{*g} \circ \Psi(G)$ ; d'apr\`es l'hypoth\`ese de r\'ecurrence, $\Psi(G) \in {\cal H}_{n-1}^*$ ; comme $\gamma_d^{*g}$
est homog\`ene de poids $1$ ; $\Psi(F) \in {\cal H}_n^*$.\\

Montrons que $\Psi$ est surjectif : soit $f \in {\cal H}_n^*$, montrons que $f \in Im(\Psi)$. Proc\'edons par r\'ecurrence sur $n$. Si $n=0$, alors $f=\lambda 1$,
et donc $f=\Psi(\lambda 1)$. Supposons la propri\'et\'e vraie pour tout $i<n$. Si $f \in M^2$,  on peut alors se ramener \`a $f=f_1f_2$, $\varepsilon(f_i)=0$,  $poids(f_i)<n$. D'apr\`es l'hypoth\`ese de r\'ecurrence, il existe $x_i \in \Hr$, $\Psi(x_i)=f_i$ ; alors $\Psi(x_1x_2)=\Psi(x_1)\Psi(x_2)=f_1f_2$.  
Sinon, d'apr\`es (\ref{eq13}), on peut se ramener \`a $f = \gamma_d^{*g}(h)$, $h \in (\Hr)^{*g}$. Comme $\gamma_d^{*g}$ est homog\`ene de poids $1$, on peut supposer $poids(h)=n-1$. D'apr\`es l'hypoth\`ese de r\'ecurrence, il existe $x \in \Hr$, $h=\Psi(x)$.
Alors $\Psi(B_d^+(x))=\gamma_d^{*g} \circ \Psi(x)=\gamma_d^{*g}(h)=f$.

Comme $\Psi$ est surjectif, homog\`ene de poids z\'ero, et que les composantes homog\`enes de $\Hr$ et $(\Hr)^{*g}$ ont la m\^ eme dimension finie,
$\Psi$ est \'egalement injectif.

\begin{theo}
\label{theo27}
$\Psi$ est un isomorphisme d'alg\`ebres de Hopf gradu\'ees entre $\Hr$ et $(\Hr)^{*g}$.
\end{theo}

\begin{theo}
\label{theo28}
Il existe une unique forme bilin\'eaire $(\hspace{1mm},)$ sur $\Hr$ v\'erifiant :
\begin{enumerate}
\item $\forall x \in \Hr$, $(1,x)=\varepsilon(x)$ ;
\item $\forall x_1,x_2,y \in \Hr$, $(x_1x_2,y)=(x_1\otimes x_2, \Delta(y))$ ;
\item $\forall x,y \in \Hr$, $\forall d \in \cal D$, $(B_d^+(x),y)=(x,\gamma_d(y)).$\\
\\
De plus, $(\hspace{1mm},)$ v\'erifie :
\item Si $x,y \in \Hr$, homog\`enes avec des poids diff\'erents, $(x,y)=0$ ;
\item $(\hspace{1mm},)$ est sym\'etrique et non d\'eg\'en\'er\'ee ;
\item $\forall x,y \in \Hr$, $(S(x),y)=(x,S(y))$ ;
\item soit $(e_F)_{F \in \forets}$ la base de $\Hr$ d\'efinie par $(e_F,G)=\delta_{F,G}$. Alors :
\begin{description}
\item[\textnormal{a)}] $e_F$ est homog\`ene de poids \'egal au poids de $F$ ;
\item[\textnormal{b)}] $(e_t)_{t \in \trees}$ est une base de $prim(\Hr)$ ;
\item[\textnormal{c)}] $\forall t_1 \ldots t_n \in \forets$,
$$\hspace{-9cm}\Delta(e_{t_1 \ldots t_n})=\sum_{i=0}^{n} e_{t_1\ldots t_i}\otimes e_{t_{i+1} \ldots t_n}.$$
\end{description}
\end{enumerate}
\end{theo}
{\it Preuve :}

Unicit\'e : soit $F,G \in \forets$. Montrons par r\'ecurrence sur $poids(F)$ que $(F,G)$ est enti\`erement d\'etermin\'e. Si $F=1$, alors 2 permet de conclure.
Supposons que $(F',G)$ est d\'etermin\'e si $poids(F')<poids(F)$. Si $F \in \trees$, posons $F=B_d^+(F')$ ; alors $(F,G)=(F',\gamma_d(G))$. Sinon, \'ecrivons $F=F_1F_2$,
$poids(F_i)<poids(F)$ pour $i=1,2$. Alors $(F,G)=(F_1 \otimes F_2, \Delta(G))$.\\

Existence :
on pose $(x,y)=\Psi(x)(y)$. Montrons que $(\hspace{1mm},)$ v\'erifie 1-7.

1. $\Psi(1) =1_{(\Hr)^{*g}}=\varepsilon$.

2. $\Psi(x_1x_2)(y)=\left(\Psi(x_1)\Psi(x_2)\right)(y)=\left(\Psi(x_1) \otimes \Psi(x_2)\right)(\Delta(y))$ par d\'efinition du produit de $(\Hr)^{*g}$.\\
On a de plus :
\begin{eqnarray*}
\Psi(x)(y_1y_2)&=&\Delta(\Psi(x))(y_1 \otimes y_2) \mbox{ par d\'efinition du coproduit de $(\Hr)^{*g}$},\\
&=&\Psi \otimes \Psi(\Delta(x))(y_1 \otimes y_2) \mbox{ car $\Psi$ est un morphisme de cog\`ebres},
\end{eqnarray*}
et donc :
\begin{equation}
\label{eq14}
(x,y_1y_2)=(\Delta(x),y_1 \otimes y_2), \esp \forall x,y_1,y_2 \in \Hr.
\end{equation}

3. $\Psi(B_d^+(x))(y)=\gamma_d^{*g} \circ\Psi(x) (y)=\Psi(x) \left(\gamma_d(y)\right)$.

4. D\'ecoule du fait que $\Psi$ est homog\`ene.

5. Soient $F,G\in \forets$ ; il s'agit de montrer que $(F,G)=(G,F)$. Si $F$ et  $G$ ont des poids diff\'erents, alors $(F,G)=(G,F)=0$ d'apr\`es 4. Supposons donc que $F$ et $G$ ont m\^ eme poids $n$ et proc\'edons par r\'ecurrence sur $n$.
Si $n=0$, alors $F=G=1$, le r\'esultat est trivial. Si $n=1$, alors $F=\bullet_d$, $G= \bullet_{d'}$, 
et $(F,G)=(G,F)=\delta_{d,d'}$.
Supposons $n\geq 2$ et le r\'esultat acquis pour tout $i<n$. 3 cas se pr\'esentent :
\begin{description}
\item[\textnormal{a)}] $F,G \in \trees$ : alors si $d\in {\cal D}$ est la d\'ecoration de la racine de $F$,
$$(F,G)=(B^-(F), \gamma_d(G))=(B^-(F),0)=0.$$
De m\^ eme, $(G,F)=0$.
\item[\textnormal{b)}] $F \in \forets -\trees$ : on peut alors \'ecrire $F=F_1F_2$, $F_i \in \forets$, $poids(F_i)<n$. 
\begin{eqnarray*}
(F_1F_2,G)&=&(F_1 \otimes F_2,\Delta(G))\\
&=&(\Delta(G),F_1\otimes F_2)\\
&=&(G,F_1F_2)
\end{eqnarray*}
(On a utilis\'e 2 pour la premi\`ere \'egalit\'e, l'hypoth\`ese de r\'ecurrence pour la deuxi\`eme, et (\ref{eq14}) pour la troisi\`eme.)
\item[\textnormal{c)}] $G \in \forets-\trees$ : m\^eme calcul.
\end{description}
De plus, $(\Hr)^{\perp}=Ker(\Psi)=(0)$, et donc $(\hspace{1mm},)$ est non d\'eg\'en\'er\'ee.

6. $\Psi(S(x))(y)=S\left( \Psi(x)\right)(y)=\Psi(x)(S(y))$ par d\'efinition de l'antipode de $(\Hr)^{*g}$.

7. Soit $(Z_F)_{F \in \forets}$ la base de $(\Hr)^{*g}$ d\'efinie par $Z_F(G)=\delta_{F,G}$, $\forall G \in \forets$ (lemme \ref{lemme1}). Il s'agit d'une base form\'ee d'\'el\'ements homog\`enes.
On a imm\'ediatement $e_F =\Psi^{-1}(Z_F)$. Comme $\Psi$ est homog\`ene de poids z\'ero, $\Psi^{-1}$ l'est aussi, et donc $e_F$ est homog\`ene de poids \'egal au poids de $F$.

Dans $(\Hr)^{*g}$, on a :
\begin{eqnarray*}
vect(Z_t/t\in \trees)&=&[vect(F/F \in \forets-\trees)]^{\perp}\\
&=&((1)\oplus M^2)^{\perp}\\
&=&Prim((\Hr)^{*g}).
\end{eqnarray*}
Comme $\Psi^{-1}(vect(Z_t/t \in \trees))=vect(e_t/t \in \trees)$ et que $\Psi^{-1}$ est un isomorphisme d'alg\`ebres de Hopf, ceci d\'emontre le point b).
\begin{eqnarray*}
(\Delta(e_{t_1 \ldots t_n}),F \otimes G)&=&(e_{t_1 \ldots t_n},FG)\\
&=&\delta_{t_1 \ldots t_n,FG}\\
&=&(\sum_{i=0}^{n} e_{t_1\ldots t_i}\otimes e_{t_{i+1} \ldots t_n},F \otimes G).
\end{eqnarray*}
Comme $(\hspace{1mm},)$ est non d\'eg\'en\'er\'ee, on obtient le point c). $\Box$

\begin{prop}
Soit $F \in \forets$.
\begin{eqnarray*}
B_d^+(e_F)&=&e_{F \bullet_d} \esp ;\\
\gamma_d(e_F)&=&e_{B^-(F)} \mbox{ si $F\in \trees$ et si la racine de $F$ est d\'ecor\'ee par $d$,}\\
&=&0\mbox{ sinon.}
\end{eqnarray*}
\end{prop}
{\it Preuve :} soit $G \in \forets$.
\begin{eqnarray*}
(B_d^+(e_F),G)&=&(e_F, \gamma_d(G))\\
&=& \delta_{F,H}  \mbox{ si $G$ est de la forme $H\bullet_d$,}\\
&=&0\mbox{ sinon.}\\
(e_{F\bullet_d},G)&=&\delta_{F\bullet_d,G}\\
&=& \delta_{F,H}  \mbox{ si $G$ est de la forme $H\bullet_d$,}\\
&=&0\mbox{ sinon.}
\end{eqnarray*}
D'o\`u le premier r\'esultat.
\begin{eqnarray*}
(\gamma_d(e_F),G)&=&(e_F,B_d^+(G))\\
&=&\delta_{F,B_d^+(G)}.
\end{eqnarray*}
Par suite, si $F$ n'est pas un arbre dont la racine est d\'ecor\'ee par $d$, 
$(\gamma_d(e_F),G)=0$ $\forall G \in \forets$, et donc $\gamma_d(e_F)=0$.
Sinon, $(\gamma_d(e_F),G)=\delta_{B^-(F),G}$ et donc $\gamma_d(e_F)=e_{B^-(F)}$. $\Box$

\subsection{L'alg\`ebre de Lie $Prim(\Hr)$}
\spa Soient $t,t_1,t_2 \in \trees$. On note $n(t_1,t_2 ;t)$ le nombre de coupes \'el\'ementaires $c$ de $t$ telles que $P^c(t)=t_1$ et $R^c(t)=t_2$.
Pour $t_1,t_2$ fix\'es, seul un nombre fini d'\'el\'ements $t$ de $\trees$ v\'erifient $n(t_1,t_2 ;t)\neq 0$.

\begin{prop}
\label{prop29}
Soient $t_1,t_2 \in \trees$. On a :
$$ [e_{t_1},e_{t_2}]=\sum_{t \in \trees} (n(t_1,t_2 ;t)-n(t_2,t_1 ;t))e_t.$$
\end{prop}
{\it Preuve :} soit $F\in \forets$, $F\neq 1$.
$$(e_{t_1}e_{t_2},F)=\sum_{c \in \ad_*(F)} \delta_{t_1,P^c(F)}\delta_{t_2,R^c(F)}.$$

Supposons $F=t'_1 \ldots t'_m$, $m\geq 3$. Alors aucune coupe admissible de $F$ ne v\'erifie $P^c(F) \in \trees,$ $R^c(F) \in \trees$.
Donc $(e_{t_1}e_{t_2},F)=0$.\\

Supposons $F=t'_1t'_2$. Les coupes admissibles de $F$ v\'erifiant $P^c(F) \in \trees,$ $R^c(F) \in \trees$ sont $(c_v(t'_1),c_t(t'_2))$ et $(c_t(t'_1),c_v(t'_2))$, o\`u 
$c_v(t'_i)$ d\'esigne la coupe vide de $t'_i$, et $c_t(t'_i)$ d\'esigne la coupe totale de $t'_i$. Donc :
\begin{eqnarray*}
 (e_{t_1}e_{t_2},F)&=&(e_{t_1},t'_1)(e_{t_2},t'_2)+(e_{t_1},t'_2)(e_{t_2},t'_1)\\
 &=&\delta_{t_1t_2,t'_1t'_2}+\delta_{t_1t_2,t'_2t'_1}\\
&=& \delta_{t_1t_2,F}+\delta_{t_2t_1,F}.
\end{eqnarray*}

Supposons $F\in \trees$. Les seules coupes admissibles de $F$ telles que  $P^c(F) \in \trees,$ $R^c(F) \in \trees$
sont les coupes \'el\'ementaires, donc :
$$(e_{t_1}e_{t_2},F)=n(t_1,t_2 ;F).$$

Et donc :
$$e_{t_1}e_{t_2}=e_{t_1t_2}+e_{t_2t_1} +\sum_{t \in \trees} n(t_1,t_2 ;t)e_t.$$
Le r\'esultat annonc\'e en d\'ecoule imm\'ediatement. $\Box$\\

{\it Remarque :} on donnera une autre expression de $[e_{t_1}, e_{t_2}]$ dans le corollaire \ref{cor52}.

\section{Relation d'ordre sur $\forets$}

\spa On consid\`ere le sous-anneau de $\Hr$ form\'e par l'alg\`ebre libre engendr\'ee sur $\mathbb{Z}$ par $\trees$. On le note $\Ar$.
Une $\mathbb{Z}$-base de ce $\mathbb{Z}$-module libre est form\'ee des \'el\'ements de $\forets$. 
D'apr\`es (\ref{defDelta}), $\Delta$ envoie $\Ar$ sur $\Ar \otimes \Ar$. Il s'agit donc d'une $\mathbb{Z}$- alg\`ebre de Hopf. Le but de cette section est de montrer que
$(e_F)$ est une autre $\mathbb{Z}$-base de $\Ar$.

\subsection{D\'efinition}
\spa On suppose ici que ${\cal D}$ est fini, totalement ordonn\'e. Soient $F=t_1\ldots t_n$, $G=t'_1 \ldots t'_m$ deux \'el\'ements distincts de $\forets$. 
On pose $t_n=B_d^+(F')$, $t'_m=B_{d'}^+(G').$\\
\begin{tabular}{ccccl}
$F > G$& si & $poids(F)>poids(G) \esp ;$\\
  & ou si & $poids(F)=poids(G)$,& $n=1$, & $m\geq 2 \esp ;$\\
  & ou si & $poids(F)=poids(G)$,& $n=m=1$,& $d>d'$ ;\\
  & ou si & $poids(F)=poids(G)$,& $n=m=1$,& $d=d',\esp F' > G'$ ;\\
  & ou si & $poids(F)=poids(G)$,& $n \mbox{ ou }m >1$,&$\exists i,\esp t_n=t'_m, \ldots, t_{n-i+1}=t'_{m-i+1}$,\\
&&&& $t_{n- i}> t_{m-i}'$.
\end{tabular}
On d\'efinit ainsi par r\'ecurrence sur le poids une relation d'ordre totale sur $\forets$, v\'erifiant :
\begin{eqnarray*}
\forall F,G \in \forets,\esp F \geq G &\Rightarrow & B_d^+(F) \geq B_d^+(G) \esp ;\\
 \forall F,G,H \in \forets, \esp F \geq G &\Rightarrow & FH \geq GH,\\
 &&HF \geq HG.
\end{eqnarray*}

\begin{figure}[h]
\framebox(450,180){
\begin{picture}(0,0)(150,-50)
$1<$
\begin{picture}(10,5)(-3,0)
\put(0,0){\circle*{5}}
\end{picture}
$<$ 
\begin{picture}(15,15)(-3,0)
\put(0,0){\circle*{5}}
\put(10,0){\circle*{5}}
\end{picture}
$<$
\begin{picture}(15,10)(-3,0)
\put(0,0){\circle*{5}}
\put(0,10){\circle*{5}}
\put(0,0){\line(0,1){10}}
\end{picture}
$<$
\begin{picture}(25,10)(0,0)
\put(0,0){\circle*{5}}
\put(10,0){\circle*{5}}
\put(20,0){\circle*{5}}
\end{picture}
$<$
\begin{picture}(25,10)(-3,0)
\put(0,0){\circle*{5}}
\put(10,0){\circle*{5}}
\put(0,10){\circle*{5}}
\put(0,0){\line(0,1){10}}
\end{picture}
$<$
\begin{picture}(25,10)(-3,0)
\put(0,0){\circle*{5}}
\put(10,0){\circle*{5}}
\put(10,10){\circle*{5}}
\put(10,0){\line(0,1){10}}
\end{picture}
$<$
\begin{picture}(25,15)(-3,0)
\put(10,0){\circle*{5}}
\put(10,0){\line(-1,1){10}}
\put(10,0){\line(1,1){10}}
\put(0,10){\circle*{5}}
\put(20,10){\circle*{5}}
\end{picture}
$<$
\begin{picture}(15,30)(-3,0)
\put(0,0){\circle*{5}}
\put(0,10){\circle*{5}}
\put(0,0){\line(0,1){10}}
\put(0,10){\line(0,1){10}}
\put(0,20){\circle*{5}}
\end{picture}
$<$
\begin{picture}(0,0)(330,50)
\begin{picture}(35,10)(-3,0)
\put(0,0){\circle*{5}}
\put(10,0){\circle*{5}}
\put(20,0){\circle*{5}}
\put(30,0){\circle*{5}}
\end{picture}
$<$
\begin{picture}(35,10)(-3,0)
\put(0,0){\circle*{5}}
\put(10,0){\circle*{5}}
\put(0,10){\circle*{5}}
\put(00,0){\line(0,1){10}}
\put(20,0){\circle*{5}}
\end{picture}
$<$
\begin{picture}(35,10)(-3,0)
\put(0,0){\circle*{5}}
\put(10,0){\circle*{5}}
\put(10,10){\circle*{5}}
\put(10,0){\line(0,1){10}}
\put(20,0){\circle*{5}}
\end{picture}
$<$
\begin{picture}(35,15)(-3,0)
\put(10,0){\circle*{5}}
\put(10,0){\line(-1,1){10}}
\put(10,0){\line(1,1){10}}
\put(0,10){\circle*{5}}
\put(20,10){\circle*{5}}
\put(30,0){\circle*{5}}
\end{picture}
$<$
\begin{picture}(25,30)(-3,0)
\put(0,0){\circle*{5}}
\put(0,10){\circle*{5}}
\put(0,0){\line(0,1){10}}
\put(0,10){\line(0,1){10}}
\put(0,20){\circle*{5}}
\put(10,0){\circle*{5}}
\end{picture}
$<$
\begin{picture}(35,10)(-3,0)
\put(0,0){\circle*{5}}
\put(10,0){\circle*{5}}
\put(20,10){\circle*{5}}
\put(20,0){\line(0,1){10}}
\put(20,0){\circle*{5}}
\end{picture}
$<$
\begin{picture}(25,10)(-3,0)
\put(0,0){\circle*{5}}
\put(0,10){\circle*{5}}
\put(0,0){\line(0,1){10}}
\put(10,0){\circle*{5}}
\put(10,10){\circle*{5}}
\put(10,0){\line(0,1){10}}
\end{picture}
$<$
\end{picture}
\end{picture}
\begin{picture}(0,0)(150,60)
\begin{picture}(35,15)(-3,0)
\put(20,0){\circle*{5}}
\put(20,0){\line(-1,1){10}}
\put(20,0){\line(1,1){10}}
\put(10,10){\circle*{5}}
\put(30,10){\circle*{5}}
\put(0,0){\circle*{5}}
\end{picture}
$<$
\begin{picture}(25,30)(-3,0)
\put(10,0){\circle*{5}}
\put(10,10){\circle*{5}}
\put(10,0){\line(0,1){10}}
\put(10,10){\line(0,1){10}}
\put(10,20){\circle*{5}}
\put(0,0){\circle*{5}}
\end{picture}
$<$
\begin{picture}(25,10)(0,0)
\put(0,10){\circle*{5}}
\put(10,10){\circle*{5}}
\put(20,10){\circle*{5}}
\put(10,0){\circle*{5}}
\put(10,0){\line(-1,1){10}}
\put(10,0){\line(1,1){10}}
\put(10,0){\line(0,1){10}}
\end{picture}
$<$
\begin{picture}(25,10)(-3,0)
\put(0,10){\circle*{5}}
\put(20,10){\circle*{5}}
\put(0,20){\circle*{5}}
\put(0,10){\line(0,1){10}}
\put(10,0){\circle*{5}}
\put(10,0){\line(-1,1){10}}
\put(10,0){\line(1,1){10}}
\end{picture}
$<$
\begin{picture}(25,10)(-3,0)
\put(0,10){\circle*{5}}
\put(20,10){\circle*{5}}
\put(20,20){\circle*{5}}
\put(20,10){\line(0,1){10}}
\put(10,0){\circle*{5}}
\put(10,0){\line(-1,1){10}}
\put(10,0){\line(1,1){10}}
\end{picture}
$<$
\begin{picture}(25,15)(-3,0)
\put(10,10){\circle*{5}}
\put(10,10){\line(-1,1){10}}
\put(10,10){\line(1,1){10}}
\put(0,20){\circle*{5}}
\put(20,20){\circle*{5}}
\put(10,0){\circle*{5}}
\put(10,0){\line(0,1){10}}
\end{picture}
$<$
\begin{picture}(15,30)(-3,0)
\put(0,10){\circle*{5}}
\put(0,20){\circle*{5}}
\put(0,10){\line(0,1){10}}
\put(0,20){\line(0,1){10}}
\put(0,30){\circle*{5}}
\put(0,0){\circle*{5}}
\put(0,0){\line(0,1){10}}
\end{picture}
$< \ldots$
\end{picture}
}
\framebox(450,100){
\begin{picture}(0,0)(150,-20)
$1<$
\begin{picture}(15,30)(-5,0)
\put(0,0){\circle*{5}}
\put(5,-5){1}
\end{picture}
$<$
\begin{picture}(15,30)(-5,0)
\put(0,0){\circle*{5}}
\put(5,-5){2}
\end{picture}
$<$
\begin{picture}(30,30)(-5,0)
\put(0,0){\circle*{5}}
\put(5,-5){1}
\put(15,0){\circle*{5}}
\put(20,-5){1}
\end{picture}
$<$
\begin{picture}(30,30)(-5,0)
\put(0,0){\circle*{5}}
\put(5,-5){2}
\put(15,0){\circle*{5}}
\put(20,-5){1}
\end{picture}
$<$
\begin{picture}(30,30)(-5,0)
\put(0,0){\circle*{5}}
\put(5,-5){1}
\put(15,0){\circle*{5}}
\put(20,-5){2}
\end{picture}
$<$
\begin{picture}(30,30)(-5,0)
\put(0,0){\circle*{5}}
\put(5,-5){2}
\put(15,0){\circle*{5}}
\put(20,-5){2}
\end{picture}
$<$
\end{picture}
\begin{picture}(0,0)(80,30)
\begin{picture}(15,30)(-5,0)
\put(0,0){\circle*{5}}
\put(0,0){\line(0,1){15}}
\put(0,15){\circle*{5}}
\put(5,-5){1}
\put(5,10){1}
\end{picture}
$<$
\begin{picture}(15,30)(-5,0)
\put(0,0){\circle*{5}}
\put(0,0){\line(0,1){15}}
\put(0,15){\circle*{5}}
\put(5,-5){1}
\put(5,10){2}
\end{picture}
$<$
\begin{picture}(15,30)(-5,0)
\put(0,0){\circle*{5}}
\put(0,0){\line(0,1){15}}
\put(0,15){\circle*{5}}
\put(5,-5){2}
\put(5,10){1}
\end{picture}
$<$
\begin{picture}(15,30)(-5,0)
\put(0,0){\circle*{5}}
\put(0,0){\line(0,1){15}}
\put(0,15){\circle*{5}}
\put(5,-5){2}
\put(5,10){2}
\end{picture}
$<\ldots$
\end{picture}
}
\caption{{\it la relation d'ordre sur les \'el\'ements de ${\cal F}_{P,R}$ et $\forets$ avec ${\cal D}=\{1,2\}$.}}
\end{figure}

\subsection{Application \`a la forme bilin\'eaire $(\hspace{1mm},)$}

\begin{lemme}
\label{lem29}
$\forall F,G \in \forets,$ $(F,G) \in \mathbb{N}$.
\end{lemme}
{\it Preuve :} r\'ecurrence sur $poids(F)$. Si $F=1$, alors $(F,G)=\varepsilon(G)=0$ ou $1$.
Supposons la propri\'et\'e vraie pour toute for\^et de poids $<n$. Soit $F \in \forets$, de poids $n$.
Si $F=B^+_d(F')$, alors $(F,G)=(F',\gamma_d(G))$ ; comme $\gamma_d(G)$ est nul ou dans $\forets$, $(F,G) \in \mathbb{N}$. 
Si $F=F_1F_2$, $(F,G)=\sum_{c \in \ad(G)}(F_1,P^c(G))(F_2,R^c(G)) \in \mathbb{N}$. $\Box$ \\

Pour $F \in \forets$, on pose ${\cal M}_F=\{ G \in \forets/ (F,G)\neq 0\}$. Comme $(\hspace{1mm},)$ est non d\'eg\'en\'er\'ee,  
${\cal M}_F$ est non vide. Par la propri\'et\'e 4 du th\'eor\`eme \ref{theo28}, ${\cal M}_F$ est fini. On pose
$m(F)=\max \esp ({\cal M}_F)$ ; $m(F)$ est un \'el\'ement de $\forets$ de m\^eme poids que $F$.

\begin{prop}
\label{prop40}
\begin{enumerate}
\item $m(\bullet_d)=\bullet_d$, et $(\bullet_d,m(\bullet_d))=1$, $\forall d \in \cal D$.
\item Si $F=B_d^+(F')$, alors $m(F)=m(F')\bullet_d$, et $(F,m(F))=(F',m(F'))$.
\item Si $F=F' \bullet_d$, alors $m(F)=B_d^+(m(F'))$, et $(F,m(F))=(F',m(F'))$.
\item Si $F=F_1t$, $t = B_d^+(t_1)$, avec $t_1 \in \forets-\{1\}$, alors $m(F)=m(t_1)B_d^+(m(F_1))$, et $(F,m(F))=(F_1,m(F_1))(t_1,m(t_1))$.
\end{enumerate}
\end{prop}
{\it Preuve :}

1. En effet, $(\bullet_d,\bullet_{d'})=(1, \gamma_d(\bullet_{d'}))=\delta_{d,d'}(1,1)=\delta_{d,d'}$.\\

2.  $(F,m(F') \bullet_d)=(F', \gamma_d(m(F') \bullet_d))=(F',m(F'))\neq 0$, donc $m(F) \geq m(F') \bullet_d$.
Soit $G \in \forets$, $(F,G) \neq 0$. $(F,G)=(F', \gamma_d(G))$ donc $\gamma_d(G)\neq 0$ : posons $G=G' \bullet_d$.
Supposons $G' > m(F')$, alors $(F,G)=(F',G')=0$ par d\'efinition de $m(F')$. Donc $G'\leq m(F')$, 
d'o\`u $G'\bullet_d \leq m(F')\bullet_d$, et donc $m(F) \leq m(F') \bullet_d$.\\

3. $(F'\bullet_d, B^+_d(m(F')))=(\gamma_d(F' \bullet_d),m(F'))=(F',m(F'))\neq 0$. Donc $m(F)\geq B_d^+(m(F'))$ : \'etant donn\'ee la d\'efinition de $\geq$,
$m(F) \in \trees$ : posons $m(F)=B^+_{d'}(G)$. Comme $\gamma_{d'}(F\bullet_d)=0$ si $d'\neq d$, n\'ecessairement $d'=d$.
Supposons $G>m(F')$, alors $(F,B_d^+(G))=(F',G)=0$ par d\'efinition de $m(F')$, et donc $G \leq m(F')$. Comme $B^+_d$ est croissante, $m(F) \leq B^+_d(m(F'))$.\\

4. Soit $G \in \forets$, $d' \in \cal D$. $(F_1t,B_{d'}^+(G))=(\gamma_{d'}(F_1t),G)$. Comme $t \neq \bullet_{d'}$, ceci est nul. Par suite, $m(F) \notin \trees$. Posons
$m(F)=GS,$ $ G \neq 1$, $S=B_{d'}^+(S_1)$.\\
\begin{eqnarray*}
\mbox{Posons } \Delta(F_1)&=&F_1 \otimes 1 +1 \otimes F_1 + \sum F_1' \otimes F_1'', \esp F_1',F_1'' \mbox{ for\^ets non vides,}\\
\Delta(t)&=&t \otimes 1 + 1 \otimes t + \sum t' \otimes t'', \esp t' \mbox{ for\^et non vide et $t''\in \trees$}.
\\
\mbox{Alors } \Delta(F)&=&F\otimes 1 + F_1 \otimes t + \sum F_1 t' \otimes t'' + t \otimes F_1 + 1 \otimes F + \sum t' \otimes F_1t''\\
&&+ \sum F_1't \otimes F_1'' + \sum F_1' \otimes F_1''t + \sum\sum F_1't' \otimes F_1''t''.
\end{eqnarray*}
Comme $t \neq \bullet_{d'}$, on a $(t,S)=(F_1''t,S)=0.$ De plus, Si $t''\neq \bullet_{d'}$, alors $(t'',S)=(F_1t'',S)=(F_1''t'',S)=0$. 
Si $t''=\bullet_{d'}$, alors $t''=\bullet_d$ et $t'=B^-(t)=t_1$.
On a donc : 
\begin{eqnarray*}
(F,GS)&=&(\Delta(F),G\otimes S)\\
&=&(F_1t_1,G)(\bullet_d,S) + (t,G)(F_1,S) + ( t_1,G)(F_1 \bullet_d,S)+ \sum (F_1't,G)(F_1'',S)  \\
&&+\sum (F_1't_1,G)(F_1''\bullet_d,S).
\end{eqnarray*}
Pour $G=m(t_1)$ et $S=B_d^+(m(F_1))$ : alors $poids(G)=poids(t_1)=poids(t)-1$. A l'aide de la propri\'et\'e 4 du th\'eor\`eme \ref{theo28}, on a :
$$ (F_1t_1,G)=(t,G)=(F_1't,G)=(F_1't_1,G)=0.$$
Alors $(F,GS)=(t_1,m(t_1))(F_1 \bullet_d,B_d^+(m(F_1))=(t_1,m(t_1))(F_1,m(F_1)) \neq 0$.
 Donc $GS\geq m(T_1) B_d^+(m(F_1))$.\\
Etant donn\'ee la d\'efinition de $\geq$, si on avait $poids(S)<poids(B_d^+(m(F_1))$, on aurait $GS<m(t_1)B_d^+(m(F_1))$. 
Donc $poids(S)\geq poids(B_d^+(m(F_1))=1+poids(F_1)$ : par suite, $(\bullet_d,S)=(F_1,S)=(F_1'',S)=(F_1''\bullet_d,S)=0$.
Donc $(F,GS)=(t_1,G)(F_1\bullet_d,S)=(t_1,G)(\gamma_{d'}(F_1\bullet_d),S_1)\neq 0$. Par suite, $d=d'$. De plus, $S_1\leq m(F_1)$ et $G \leq m(T_1)$,
d'o\`u $GS \leq m(t_1) B_d^+(m(F_1))$. $\Box$
  
\begin{prop}
\begin{enumerate}
\item $(F,m(F))=1$, $\forall F \in \forets$.
\item $m(m(F))=F$, $\forall F \in \forets$.
\end{enumerate}
\end{prop}
{\it Preuve :}\\
1. R\'ecurrence facile sur le poids de $F$.\\
2. R\'ecurrence sur le poids de $F$. Si $poids(F) =0$ ou $1$, c'est imm\'ediat. Supposons $m(m(F'))=F'$, $\forall F' \in \forets$, $poids(F')<n$.
Soit $F$ de poids $n$. 

Si $F=B_d^+(F')$, alors $m(F)=m(F') \bullet_d$, et donc $m(m(F))=B_d^+(m(m(F')))=B_d^+(F')=F$.

Si $F=F' \bullet_d$, alors $m(F)=B_d^+(m(F'))$, et donc $m(m(F))=m(m(F'))\bullet_d=F' \bullet_d=F$.

Si $F=F_1t$, $t = B_d^+(t_1)$, avec $t_1 \in \forets-\{1\}$, alors $m(F)=m(t_1)B_d^+(m(F_1))$, et 
$m(m(F))=m(m(F_1))B_d^+(m(m(t_1)))=F_1 B_d^+(t_1)=F_1 t=F$. $\Box$\\

Par suite, $m :\{F \in \forets/ poids(F)=n\} \longmapsto \{F \in \forets/ poids(F)=n\}$ est une bijection.
Soit ${\cal B}_n=(F_i)_{i \leq r_n}$ la base de ${\cal H}_n$ form\'ee des for\^ets de poids $n$, index\'ees de sorte que 
$m(F_1) < \ldots < m(F_{r_n})$. Soit $A_n$ la matrice de la forme bilin\'eaire $(\hspace{1mm},)$ restreinte sur ${\cal H}_n\times {\cal H}_n$ dans cette base.
Les coefficients de $A_n$ sont entiers (lemme \ref{lem29}).

 Soit $M_n=(\delta_{F_i,m(F_j)})_{1\leq i,j\leq r_n}$ la  matrice de permutation associ\'ee \`a $m_{\mid \{F \in \forets/ poids(F)=n\} }$.
Soit $B_n=A_nM_n$.
\begin{eqnarray*}
(B_n)_{i,j}&=&\sum_{k} (A_n)_{i,k}\delta_{F_k,m(F_j)}\\
&=&\sum_{k} (F_i,F_k) \delta_{F_k,m(F_j)}\\
&=&(F_i,m(F_j)).
\end{eqnarray*}
Comme $(F_i,m(F_j))=0$ si $m(F_j)>m(F_i)$, c'est-\`a-dire si $j>i$, et que $(F_i,m(F_i))=1$ :
$$B_n = \left[ \begin{array}{cccc}
   1&0&\ldots &0\\
   \vdots&1&\ddots & \vdots \\
    \vdots &&\ddots &0 \\
    \ldots &\ldots & \ldots &1
\end{array} \right]. $$
Par suite, $det(B_n)=1$. Comme $M_n$ est une matrice de permutation, $det(M_n)=\pm 1$, et donc $det(A_n)=\pm1$.
D'o\`u $A_n \in GL_{r_n}(\mathbb{Z})$.\\
Soit $P_n=Pass((F_i)_{i\leq r_n},(e_{F_i})_{i\leq r_n})$, c'est-\`a-dire $e_{F_j}=\sum_i (P_n)_{i,j} F_i$.
\begin{eqnarray*}
(A_n P_n)_{i,j}&=&\sum_k (F_i,F_k) (P_n)_{k,j}\\
&=& (F_i, \sum_k (P_n)_{k,j}F_k)\\
&=& (F_i,e_{F_j})\\
&=& \delta_{i,j}, \mbox{ par d\'efinition de la base $(e_F)_{F \in \forets}$.}
\end{eqnarray*}
Donc $P_n=A_n^{-1}$, et donc $P_n \in GL_{r_n}(\mathbb{Z})$. D'o\`u :
\begin{theo}
\label{theo33}
$(e_F)_{F \in \forets}$ est une $\mathbb{Z}$-base de $\Ar$.
\end{theo} 
\subsection{Cas o\`u l'ensemble des d\'ecorations $\cal D$ est infini}
\spa Dans ce cas, $\Hr$ ne v\'erifie pas la condition $(C_2)$. On ne peut alors plus munir $(\Hr)^{*g}$ d'un coproduit. Cependant, $\Hr$ est l'union des ${\cal H}_{P,R}^{{\cal D}'}$,
o\`u ${\cal D}'$ parcourt l'ensemble des parties finies de $\cal D$.
On note $(\hspace{1mm},)_{\cal D'}$ la forme bilin\'eaire d\'ecrite par le th\'eor\`eme \ref{theo28} sur ${\cal H}_{P,R}^{{\cal D}'}$.
Par l'unicit\'e dans ce th\'eor\`eme, $(\hspace{1mm},)_{\cal D'}$ et $(\hspace{1mm},)_{\cal D''}$ co\"incident sur ${\cal H}_{P,R}^{{\cal D}'\cap {\cal D}''}$ si ${\cal D'} \cap {\cal D''}\neq \emptyset$.
Pour $x,y \in \Hr$, on peut donc d\'efinir :
$$ (x,y)=(x,y)_{\cal D'} \mbox{ o\`u $\cal D'$ est choisi tel que $x, y \in {\cal H}_{P,R}^{{\cal D}'}$}.$$
Cette forme bilin\'eaire v\'erifie les 7 propri\'et\'es du th\'eor\`eme \ref{theo28} ; on a donc une injection de $\Hr$ dans son dual de Hopf.
De la m\^eme mani\`ere, le th\'eor\`eme \ref{theo33} reste vrai.

\section{Relations d'ordre sur les sommets d'une for\^et}

\spa Dans cette section, $\cal D$ est un ensemble non vide quelconque, non n\'ecessairement fini.

\subsection{D\'efinitions}
\spa Soit $F\in \forets$, $F \neq 1$. On note $som(F)$ l'ensemble des sommets de $F$.\\
$F$ est repr\'esent\'ee par un graphe orient\'e. 
Soit $x,y \in som(F)$.
 On dira que $x\geq_{haut} y$ si il existe un trajet dans ce graphe d'origine $y$ et d'arriv\'ee $x$.
$\geq_{haut}$ est une relation d'ordre (non n\'ecessairement totale) sur $som(F)$.\\

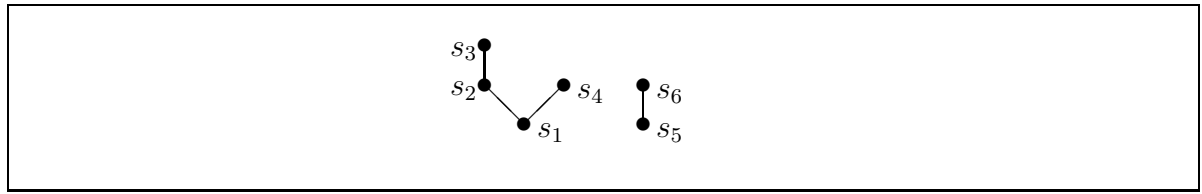
\begin{figure}[h]
\framebox(450,70){
\begin{picture}(60,100)(0,-40)
\put(0,0){\circle*{5}}
\put(0,0){\line(1,1){15}}
\put(0,0){\line(-1,1){15}}
\put(-15,15){\circle*{5}}
\put(15,15){\circle*{5}}
\put(-15,15){\line(0,1){15}}
\put(-15,30){\circle*{5}}
\put(45,0){\circle*{5}}
\put(45,15){\circle*{5}}
\put(45,0){\line(0,1){15}}
\put(50,-5){$s_5$}
\put(50,10){$s_6$}
\put(5,-5){$s_1$}
\put(20,10){$s_4$}
\put(-28,11){$s_2$}
\put(-28,26){$s_3$}
\end{picture}
}
\caption{{\it un exemple de for\^et plane enracin\'ee.}}
\end{figure}

{\it Exemple :} on note $x \nless \ngtr_{haut} y$ quand $x$ et $y$ ne sont pas comparables pour $\geq_{haut}$.
Pour la for\^et de la figure 6, on a :
$$\begin{array}{cccccc}
s_6&\geq_{haut}&s_5 \esp ;&s_6&\nless \ngtr_{haut} &s_4 \esp ;\\
s_6&\nless \ngtr_{haut} &s_3 \esp ;&s_6&\nless \ngtr_{haut} &s_2 \esp ;\\
s_6&\nless \ngtr_{haut} &s_1 \esp ;&s_5&\nless \ngtr_{haut} &s_4 \esp ;\\
s_5&\nless \ngtr_{haut} &s_3 \esp ;&s_5&\nless \ngtr_{haut} &s_2 \esp ;\\
s_5&\nless \ngtr_{haut} &s_1 \esp ;&s_4&\nless \ngtr_{haut} &s_3 \esp ;\\
s_4&\nless \ngtr_{haut} &s_2 \esp ;&s_4&\geq_{haut} &s_1 \esp ;\\
s_3&\geq_{haut} &s_2 \esp ;&s_3&\geq_{haut} &s_1 \esp ;\\
s_2&\geq_{haut}&s_1.
\end{array}$$

On d\'efinit une deuxi\`eme relation d'ordre $\geq_{gauche}$ sur $som(F)$ par r\'ecurrence sur $poids(F)$.\\
Si $F$=$\bullet_d$, $som(F)$ est r\'eduit \`a un seul \'el\'ement.\\
Si $poids(F)\geq 2$, soient $x,y$ deux sommets diff\'erents de $F$. On pose $F=t_1 \ldots t_n$ ;
on suppose que $x$ est un sommet de $t_i$ et $y$ un sommet de $t_j$ : 

si $i<j$, alors $x \geq_{gauche} y$ ; si $i>j$, alors $y \geq_{gauche} x$ ;

si $i=j$, $x$ ou $y$ est la racine de $t_i$ : alors $x$ et $y$ ne sont pas comparables pour $\geq_{gauche}$ ;

si $i=j$, $x$ et $y$ ne sont pas \'egaux \`a la racine de $t_i$ : on les compare dans $(B^-(t_i),\geq_{gauche})$.\\
$\geq_{gauche}$ est une relation d'ordre (non n\'ecessairement totale) sur $som(F)$.\\

{\it Exemple :} on note $x \nless \ngtr_{gauche} y$ quand $x$ et $y$ ne sont pas comparables pour $\geq_{gauche}$.
Pour la for\^et de la figure 6, on a :
$$\begin{array}{cccccc}
s_5&\nless \ngtr_{gauche}&s_6 \esp ;&s_4&\geq_{gauche} &s_6 \esp ;\\
s_3&\geq_{gauche}&s_6 \esp ;&s_2&\geq_{gauche} &s_6 \esp ;\\
s_1&\geq_{gauche} &s_6 \esp ;&s_4&\geq_{gauche} &s_5 \esp ;\\
s_3&\geq_{gauche}&s_5 \esp ;&s_2&\geq_{gauche} &s_5 \esp ;\\
s_1&\geq_{gauche} &s_5 \esp ;&s_3&\geq_{gauche}&s_4 \esp ;\\
s_2&\geq_{gauche} &s_4 \esp ;&s_1&\nless \ngtr_{gauche} &s_4 \esp ;\\
s_2&\nless \ngtr_{gauche} &s_3 \esp ;&s_1&\nless \ngtr_{gauche} &s_3 \esp ;\\
s_1&\nless \ngtr_{gauche}&s_2.
\end{array}$$

\subsection{Expression combinatoire de $(F,G)$}

\begin{theo}
Soit $F,G \in \forets$. On pose
${\cal I}(F,G)$ l'ensemble des bijections $f$ de $som(F)$ vers $som(G)$ v\'erifiant : \begin{enumerate}
\item $ \forall x,y \in som(F), x\geq_{haut} y \Rightarrow \esp f(x) \geq_{gauche} f(y),$
\item  $\forall x,y \in som(F), f(x)\geq_{haut} f(y) \Rightarrow\esp x \geq_{gauche} y,$
\item $\forall x \in som(F),x$ et $f(x)$ ont la m\^eme d\'ecoration.
\end{enumerate}
Alors $(F,G)=card({\cal I}(F,G))$.
\end{theo}
{\it Preuve :} c'est vrai si $poids(F)\neq poids(G)$, car alors $(F,G)=0$, et il n'y a aucune bijection de $som(F)$ vers $som(G)$.
Supposons donc $poids(F)=poids(G)=n$, et proc\'edons par r\'ecurrence sur $n$. Si $n=1$, alors $( \bullet_d, \bullet_{d'})=\delta_{d,d'}=card({\cal I}(\bullet_d,\bullet_{d'}))$ par  la condition 3. 
Supposons la propri\'et\'e v\'erifi\'ee pour tout $k<n$, et soient $F,G\in \forets$ de poids $n$. Posons $F=t_1 \ldots t_m$.\\

Si $m=1$ : posons $F=t_1=B_d^+(F')$. Alors $(F,G)=(F',\gamma_d(G))$.

Si $G$ n'est pas de la forme $G'\bullet_d$ : alors $(F,G)=0$. Supposons ${\cal I}(F,G)$ non vide et soit $f \in {\cal I}(F,G)$.
Remarquons que $\forall x \in som(F)$, $x \geq_{haut}$ $racine$ $de$ $t_1$. Par la condition 1, $\exists x' \in som(G)$, $y' \geq_{gauche} x'$, $\forall y' \in som(G)$ ($x'$ est l'image par $f$ de la racine de $t_1$).
Donc $G=G'\bullet_{d'}$, et $f(racine$ $de$ $t_1)=sommet$ $de$ $\bullet_{d'}$. Par la condition 3, $d'=d$ : on aboutit \`a une contradiction. Dans ce cas, on a bien $(F,G)=card({\cal I}(F,G))=0$.

Si $G=G'\bullet_d$ : alors pour toute $f \in {\cal I}(F,G)$, $f(racine$ $de$ $t_1)=sommet$ $de$ $\bullet_{d}$. Donc on a une bijection :
$ {\cal I}(F,G) \longmapsto {\cal I}(F',G')$, envoyant $f$ sur sa restriction \`a $som(F')$. 
Comme $(F,G)=(F',G')$, le r\'esultat est acquis.\\

Si $m>1$ : posons $F'=t_1 \ldots t_{m-1}$. Alors :
\begin{eqnarray}
\nonumber
(F,G)&=&\sum_{c \in \ad(G)} (F',P^c(G))(t_m,R^c(G))\\
\label{eqn26} &=&\sum_{c \in \ad_*(G)} (F',P^c(G))(t_m,R^c(G)).
\end{eqnarray}
Soit $f\in {\cal I}(F,G)$ ; consid\'erons $f(som(t_m))$. 
Soit $r$ un sommet de $G$, tel qu'il existe $x\in som(t_m)$, $f(x)\geq_{haut} r$.
Supposons $r \notin f(som(t_m))$. Alors $r \in f(som(F'))$. Soit $y\in som(F')$, tel que $f(y)=r$.
$f(x) \geq_{haut} f(y)$, donc d'apr\`es la condition 2, $x \geq_{gauche} y$. Or $x \in som(t_m)$, $y \in som(F')$, donc 
$y \geq_{gauche} x$, et donc $x=y$ : contradiction, car $x\in som(t_m)$, $y\in som(F')$. Donc $r \in f(som(t_m))$.
Par suite, il existe une coupe admissible $c_f$ de $G$, telle que $R^{c_f}(G)=f(som(t_m))$.
Etant donn\'ee la d\'efinition d'une coupe admissible, $c_f$ est enti\`erement d\'etermin\'ee par $R^{c_f}(G)$, et donc $c_f$ est unique.\\
De plus, $f :som(t_m) \longmapsto som(R^{c_f}(G))$ $\in {\cal I}(t_m,R^{c_f}(G))$, et 
$f :som(F') \longmapsto som(P^{c_f}(G))$ $\in {\cal I}(F',P^{c_f}(G))$.
On a donc une application :
\begin{eqnarray*}
\beta :{\cal I}(F,G)&\longmapsto & \bigcup_{c \in \ad_*(G)} {\cal I}(F',P^c(G))\times {\cal I}(t_m,R^c(G))\\
f&\longmapsto&( f_{\mid som(F')} , f_{\mid som(t_m)}) \in  {\cal I}(F',P^{c_f}(G))\times {\cal I}(t_m,R^{c_f}(G)).
\end{eqnarray*}
$\beta$ est \'evidemment injective. Montrons qu'elle est surjective. Soit $c \in \ad_*(G)$, $(f_1,f_2) \in {\cal I}(F',P^c(G))\times {\cal I}(t_m,R^c(G))$.
Soit $f :som(F) \longmapsto som(G)$, d\'efinie par $f_{\mid som(F')}=f_1$ et $f_{\mid som(t_m)}=f_2$. $f$ est une bijection. De plus, comme $f_1$ et $f_2$ v\'erifient 3, $f$ v\'erifie 3.\\

Soient $x,y\in som(F)$. Supposons que $x \geq_{haut} y$. Les sommets de $t_m$ et les sommets de $F'$ ne sont pas comparables pour $\geq_{haut}$,
et donc soit $x,y \in som(F')$, soit $x,y \in som(t_m)$. Comme $f_1$ et $f_2$ v\'erifient 1, on a $f(x) \geq_{gauche} f(y)$.\\

Supposons $f(x) \geq_{haut} f(y)$. 3 cas se pr\'esentent : \begin{enumerate}
\item Si $x,y \in som(F')$ ou $x,y \in som(t_m)$ : alors comme $f_1$ et $f_2$ v\'erifient 2, $x \geq_{gauche} y$.
\item Si $x \in som(F')$ et $y \in som(t_m)$ : alors $x \geq_{gauche} y$.
\item Si $x \in som(t_m)$ et $y \in som(F')$ : alors  $f(x) \in R^c(G)$ et $f(y) \in P^c(G)$.
Par suite, soit $f(x)$ et $f(y)$ ne sont pas comparables pour $\geq_{haut}$, soit $f(y)\geq_{haut} f(x)$.
Comme $f(x) \geq_{haut} f(y)$, on a $f(x)=f(y)$ et donc $x=y$ : on aboutit \`a une contradiction et donc ce cas est impossible.
\end{enumerate}
Par suite, $f \in {\cal I}(F,G)$, et $\beta(f)=(f_1,f_2)$. $\beta $ \'etant une bijection,
\begin{eqnarray*}
card({\cal I}(F,G))&=&\sum_{c \in \ad_*(G)} card({\cal I}(F',P^c(G))) \times card({\cal I}(t_m,R^c(G)))\\
&=&\sum_{c \in \ad_*(G)} (F',P^c(G)) (t_m,R^c(G)) \mbox{ d'apr\`es l'hypoth\`ese de r\'ecurrence},\\
&=& (F,G) \mbox{ d'apr\`es (\ref{eqn26}).} \esp \Box
\end{eqnarray*}

\subsection{Relation d'ordre totale sur les sommets de $F$}
\begin{lemme}
Soit $F \in \forets$.
\label{lemme35}
\begin{enumerate}
\item Soient $a$ et $b$ deux sommets diff\'erents de $F$. Alors :
$$ a,b \mbox{ comparables pour }\geq_{haut} \esp \Leftrightarrow a,b \mbox{ non comparables pour }\geq_{gauche}.$$
\item soient $a,a',b,b' \in som(F)$, $b \neq b'$. Alors :
$$ a\geq_{haut} b, \esp a' \geq_{haut} b', \esp b \geq_{gauche} b' \esp \Rightarrow a \geq_{gauche} a'.$$
\end{enumerate}
\end{lemme}
{\it Preuve :}

1. $\Rightarrow$ : supposons $a \geq_{haut} b$ ; quitte \`a effectuer une coupe \'el\'ementaire on peut supposer que $F\in \trees$
et que $b$ est la racine de $F$. Comme $a \neq b$, par d\'efinition de $\geq_{gauche}$, $a$ et $b$ ne sont pas comparables pour $\geq_{gauche}$.\\
$\Leftarrow$ : supposons $a$ et $b$ non comparables pour $\geq_{gauche}$. Quitte \`a effectuer une coupe \'el\'ementaire, on peut supposer que $F \in \trees$, et $a$ ou $b$
est la racine de $F$. Alors par d\'efinition de $\geq_{haut}$, $a\geq_{haut} b $ ou  $b\geq_{haut} a$.

2. Soit $c$ (respectivement $c'$) la coupe portant sur l'ar\^ete arrivant \`a $b$ (respectivement $b'$) si $b$  (respectivement $b'$) n'est  pas une racine, ou la coupe totale de l'arbre de racine $b$ (respectivement $b'$) sinon. Soit $c''=c \cup c'$. Comme $b>_{gauche} b'$, $c''$ est admissible,
et $P^{c''}(F)=tt'$, $b$ \'etant la racine de $t$, $b'$ la racine de $t'$. Comme $a \geq_{haut} b$, $a \in som(t)$ ; de m\^eme, $a' \in som(t')$.
Donc $a \geq_{gauche} a'$.
\begin{prop}
soit $F \in \forets$, $x,y \in som(F)$. On notera $x \geq_{tot} y$ si $x \geq_{haut} y$ ou $y \geq_{gauche} x$.
Alors $\geq_{tot}$ d\'efinit une relation d'ordre totale sur $som(F)$.
\end{prop}

{\it Exemple :} pour la for\^et de la figure 6, on a $s_6 \geq_{tot}s_5 \geq_{tot}s_4 \geq_{tot}s_3 \geq_{tot}s_2 \geq_{tot}s_1$.\\

{\it Preuve :} r\'eflexivit\'e : \'evident.

Transitivit\'e : supposons $x \geq_{tot}y$ et $y \geq_{tot} z$. On se ram\`ene au cas $x\neq y$, $y\neq z$.
\begin{enumerate}
\item Si $x \geq_{haut} y$ et $y \geq_{haut} z$, alors $x \geq_{haut} z$, et donc $x\geq_{tot} z$.
\item Si $y \geq_{gauche} x$ et $z \geq_{gauche} y$, alors $z \geq_{gauche} x$, et donc $x\geq_{tot} z$.
\item Si $x \geq_{haut} y$ et $z \geq_{gauche} y$ : d'apr\`es le lemme \ref{lemme35}-2 avec $a=z$, $a'=x$, $b=z$, $b'=y$, on a
$z \geq_{gauche} x$, et donc $x \geq_{tot} z$.
\item Si $y \geq_{gauche} x$ et $y \geq_{haut} z$ : si $x \geq_{haut} z$, alors $x \geq_{tot} z$.
Supposons que l'on n'ait pas $x \geq_{haut} z$.
Soit $c$ la coupe portant sur les (\'eventuelles) ar\^etes arrivant \`a $x$ et $z$. Supposons $c$ non admissible ; alors 
$z \geq_{haut} x$, et alors $y \geq_{haut} x$ par transitivit\'e de $\geq_{haut}$. Par le lemme \ref{lemme35}-1, 
n\'ecessairement $x=y$, cas que nous avons exclus. Donc $c$ est admissible. Alors $P^c(F)=tt'$, avec $x,z$ racines de $t,t'$.
Comme $y \geq_{gauche} x$, n\'ecessairement $y \in som(t)$ ; comme $y \geq_{haut} z$, $z$ est la racine de $t$, et donc $x$ est la racine de $t'$.
Par suite $z \geq_{gauche} x$, et donc $x \geq_{tot} z$.
\end{enumerate}

Antisym\'etrie : supposons $x \geq_{tot} y $ et $y \geq_{tot} x$. \begin{enumerate}
\item Si $x \geq_{haut} y$ et $y \geq_{haut} x$, alors $x=y$. 
\item Si $y \geq_{gauche} x$ et $x \geq_{gauche} y$, alors $x=y$
\item Si $x \geq_{haut} y$ et $x \geq_{gauche} y$ : par le lemme \ref{lemme35}-1, $x=y$.
\item Si $y \geq_{gauche} x$ et $y \geq_{haut} x$ : m\^eme raisonnement.
\end{enumerate}

Enfin, $ \geq_{tot}$ est totale : si $x,y \in som(F)$, alors d'apr\`es le lemme \ref{lemme35}-1,  ils sont comparables pour
$\geq_{haut}$ ou $\geq_{gauche}$, donc ils le sont pour $\geq_{tot}$. $\Box$
 
\subsection{Application au calcul de l'antipode}

\spa NB : dans cette section, la coupe totale d\'efinie dans la partie \ref{partie2} n'est pas consid\'er\'ee comme une coupe.

Soit $F$ for\^et, $c$ une coupe de $F$. On note $t_1,\ldots, t_m$ les diff\'erentes composantes connexes de $F$ apr\`es l'action de $c$ ; 
on note $x_i$ la racine de $t_i$. On suppose qu'avant l'action de $c$, on avait $x_1 \leq_{tot} \ldots \leq_{tot} x_m$.
On pose alors $W^c(F)=t_m \ldots t_1$.\\

{\it Exemple :} $F=t_1 \ldots t_m$. On fait agir la coupe vide $c_v$. Les composantes connexes sont les $t_i$, et 
$x_1 \geq_{gauche} \ldots \geq_{gauche} x_m$, d'o\`u $x_1 \leq_{tot} \ldots \leq_{tot} x_m$. Donc $W^{c_v}(t_1 \ldots t_m)=t_m \ldots t_1$.\\

On note $n_c$ le nombre de coupes \'el\'ementaires constituant $c$.
\begin{theo}
\label{theo37}
Soit $F \in \forets$, $F=t_1 \ldots t_m$.
\begin{eqnarray}
\label{eqn27}
S(F)&=&(-1)^m \sum_{ \mbox{$c$ coupe de $F$}} (-1)^{n_c}\esp  W^c(F).
\end{eqnarray}
\end{theo}
{\it Preuve :} on note $lg(t_1 \ldots t_m)=m$, $\forall t_1 \ldots t_m \in \forets$.
Soit $S'$ l'op\'erateur de $\Hr$ d\'efini par le second membre de (\ref{eqn27}).
Soit $c$ une coupe de $F=t_1 \ldots t_m$ ; on note $c_i$ la restriction de $c$ \`a $t_i$.
Alors par d\'efinition de $\geq_{tot}$, on a $W^c(F)=W^c(t_m) \ldots W^c(t_1)$. De plus, $n_c = n_{c_1} + \ldots +n_{c_m}$,
et donc :
\begin{eqnarray*}
S'(F)&=&(-1)^{lg(F)} \sum_{(c_1,\ldots c_m)} \esp \prod_{i=m}^1 (-1)^{n_{c_i}}W^{c_i}(t_i)\\
&=&S'(t_m) \ldots S'(t_1).
\end{eqnarray*}
Donc $S'$, tout comme $S$, est un antimorphisme d'alg\`ebres. Il suffit donc de montrer (\ref{eqn27}) pour $t \in \trees$.
Si $poids(t)=1$, alors $t$ est primitif, et $S'(t)=S(t)=-t$. Supposons (\ref{eqn27}) vraie pour toute for\^et de poids inf\'erieur ou \'egal \`a $n$,
et soit $t\in \trees$, de poids $n+1$.
\begin{eqnarray}
\nonumber m\circ (S \otimes Id)\circ \Delta(t)&=&S(t) + t + \sum_{c \in \ad_*(t)} S(P^c(t)) R^c(t)\\
\nonumber &=&\varepsilon(t)1\\
\label{eqn28} &=&0.
\end{eqnarray}
Soit $c$ une coupe non vide de $t$ ; on note $W^c(t)=t_1 \ldots t_m$. La composante connexe de la racine de $t$ apr\`es l'action de $c$ est $t_m$,
car tout sommet de $t$ est sup\'erieur (pour $\geq_{tot}$) \`a la racine de $t$.
Or il existe une unique coupe admissible $c'\in \ad_*(t)$ telle que $t_m=R^{c'}(t)$ ; $c$ s'\'ecrit alors $c=c' \cup c''$, avec $c''=c_{\mid P^{c'}(t)}$.
On a $n_c=n_{c'}+n_{c''}$, et $n_{c'}=lg(P^{c'}(t))$.
Alors :
\begin{eqnarray*}
S'(t)&=&-t -\sum_{\mbox{$c' \in \ad_*(t)$}} \left(\sum_{\mbox{$c''$ coupe de $R^{c'}(t)$}} (-1)^{lg(P^{c'}(t))}\esp (-1)^{n_{c''}}W^{c''}(P^{c'}(t))\esp  R^{c'}(t)\right)\\
&=&-t-\sum_{\mbox{$c' \in \ad_*(t)$}} S'(P^{c'}(t)) R^{c'}(t)\\
&=&-t -\sum_{\mbox{$c' \in \ad_*(t)$}} S(P^{c'}(t)) R^{c'}(t)\\
&=& S(t).
\end{eqnarray*}
($-t$ provient de la coupe vide. On a utilis\'e l'hypoth\`ese de r\'ecurrence pour la troisi\`eme \'egalit\'e et (\ref{eqn28}) pour la derni\`ere). $\Box$

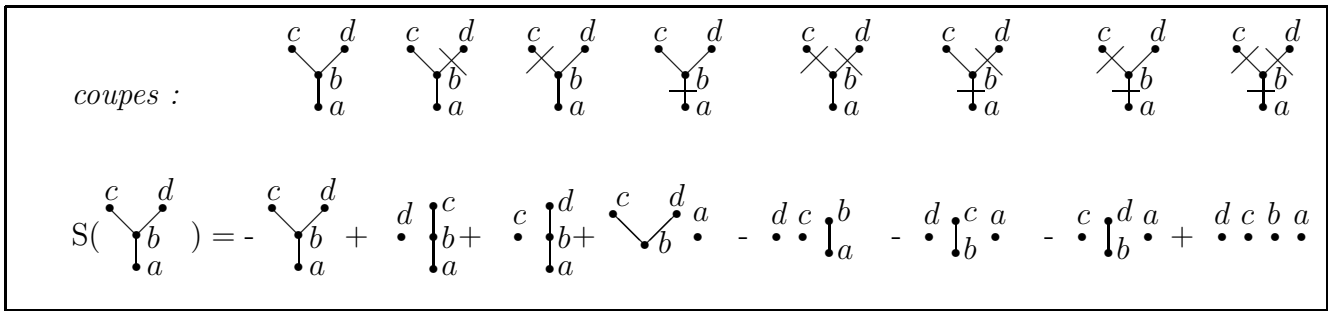
\begin{figure}[h]
{\it   coupes :}
\begin{picture}(27,30)(-52,0)
\put(-3,-1){\circle*{3}}
\put(-3,-1){\line(0,1){12}}
\put(-3,11){\circle*{3}}
\put(-3,11){\line(-1,1){10}}
\put(-3,11){\line(1,1){10}}
\put(-13,21){\circle*{3}}
\put(7,21){\circle*{3}}
\put(1,-4){$a$}
\put(1,6){$b$}
\put(5,24){$d$}
\put(-15,24){$c$}
\end{picture}
\begin{picture}(18,30)(-66,0)
\put(-3,-1){\circle*{3}}
\put(-3,-1){\line(0,1){12}}
\put(-3,11){\circle*{3}}
\put(-3,11){\line(-1,1){10}}
\put(-3,11){\line(1,1){10}}
\put(-13,21){\circle*{3}}
\put(7,21){\circle*{3}}
\put(-2,21){\line(1,-1){10}}
\put(1,-4){$a$}
\put(1,6){$b$}
\put(5,24){$d$}
\put(-15,24){$c$}
\end{picture}
\begin{picture}(27,30)(-90,0)
\put(-3,-1){\circle*{3}}
\put(-3,-1){\line(0,1){12}}
\put(-3,11){\circle*{3}}
\put(-3,11){\line(-1,1){10}}
\put(-3,11){\line(1,1){10}}
\put(-13,21){\circle*{3}}
\put(7,21){\circle*{3}}
\put(-15,12){\line(1,1){10}}
\put(1,-4){$a$}
\put(1,6){$b$}
\put(5,24){$d$}
\put(-15,24){$c$}
\end{picture}
\begin{picture}(27,30)(-107,0)
\put(-3,-1){\circle*{3}}
\put(-3,-1){\line(0,1){12}}
\put(-3,11){\circle*{3}}
\put(-3,11){\line(-1,1){10}}
\put(-3,11){\line(1,1){10}}
\put(-13,21){\circle*{3}}
\put(7,21){\circle*{3}}
\put(-9,5){\line(1,0){10}}
\put(1,-4){$a$}
\put(1,6){$b$}
\put(5,24){$d$}
\put(-15,24){$c$}
\end{picture}
\begin{picture}(27,30)(-132,0)
\put(-3,-1){\circle*{3}}
\put(-3,-1){\line(0,1){12}}
\put(-3,11){\circle*{3}}
\put(-3,11){\line(-1,1){10}}
\put(-3,11){\line(1,1){10}}
\put(-13,21){\circle*{3}}
\put(7,21){\circle*{3}}
\put(-15,12){\line(1,1){10}}
\put(-2,21){\line(1,-1){10}}
\put(1,-4){$a$}
\put(1,6){$b$}
\put(5,24){$d$}
\put(-15,24){$c$}
\end{picture}
\begin{picture}(27,30)(-154,0)
\put(-3,-1){\circle*{3}}
\put(-3,-1){\line(0,1){12}}
\put(-3,11){\circle*{3}}
\put(-3,11){\line(-1,1){10}}
\put(-3,11){\line(1,1){10}}
\put(-13,21){\circle*{3}}
\put(7,21){\circle*{3}}
\put(-2,21){\line(1,-1){10}}
\put(-9,5){\line(1,0){10}}
\put(1,-4){$a$}
\put(1,6){$b$}
\put(5,24){$d$}
\put(-15,24){$c$}
\end{picture}
\begin{picture}(27,30)(-182,0)
\put(-3,-1){\circle*{3}}
\put(-3,-1){\line(0,1){12}}
\put(-3,11){\circle*{3}}
\put(-3,11){\line(-1,1){10}}
\put(-3,11){\line(1,1){10}}
\put(-13,21){\circle*{3}}
\put(7,21){\circle*{3}}
\put(-9,5){\line(1,0){10}}
\put(-15,12){\line(1,1){10}}
\put(1,-4){$a$}
\put(1,6){$b$}
\put(5,24){$d$}
\put(-15,24){$c$}
\end{picture}
\begin{picture}(27,30)(-202,0)
\put(-3,-1){\circle*{3}}
\put(-3,-1){\line(0,1){12}}
\put(-3,11){\circle*{3}}
\put(-3,11){\line(-1,1){10}}
\put(-3,11){\line(1,1){10}}
\put(-13,21){\circle*{3}}
\put(7,21){\circle*{3}}
\put(-15,12){\line(1,1){10}}
\put(-2,21){\line(1,-1){10}}
\put(-9,5){\line(1,0){10}}
\put(1,-4){$a$}
\put(1,6){$b$}
\put(5,24){$d$}
\put(-15,24){$c$}
\end{picture}

 S(
\begin{picture}(27,30)(-13,7)
\put(-3,-1){\circle*{3}}
\put(-3,-1){\line(0,1){12}}
\put(-3,11){\circle*{3}}
\put(-3,11){\line(-1,1){10}}
\put(-3,11){\line(1,1){10}}
\put(-13,21){\circle*{3}}
\put(7,21){\circle*{3}}
\put(1,-4){$a$}
\put(1,6){$b$}
\put(5,24){$d$}
\put(-15,24){$c$}
\end{picture}
) = -
\begin{picture}(27,30)(-16,7)
\put(-3,-1){\circle*{3}}
\put(-3,-1){\line(0,1){12}}
\put(-3,11){\circle*{3}}
\put(-3,11){\line(-1,1){10}}
\put(-3,11){\line(1,1){10}}
\put(-13,21){\circle*{3}}
\put(7,21){\circle*{3}}
\put(1,-4){$a$}
\put(1,6){$b$}
\put(5,24){$d$}
\put(-15,24){$c$}
\end{picture}
 + 
\begin{picture}(27,35)(-12,-3)
\put(-3,0){\circle*{3}}
\put(9,-12){\circle*{3}}
\put(9,-12){\line(0,1){12}}
\put(9,0){\circle*{3}}
\put(9,0){\line(0,1){12}}
\put(9,12){\circle*{3}}
\put(12,-14){$a$}
\put(12,-4){$b$}
\put(-5,5){$d$}
\put(12,10){$c$}
\end{picture}
 +
\begin{picture}(27,30)(-13,-3)
\put(-3,0){\circle*{3}}
\put(9,-12){\circle*{3}}
\put(9,-12){\line(0,1){12}}
\put(9,0){\circle*{3}}
\put(9,0){\line(0,1){12}}
\put(9,12){\circle*{3}}
\put(12,-14){$a$}
\put(12,-4){$b$}
\put(-5,5){$c$}
\put(12,10){$d$}
\end{picture}
 +
\begin{picture}(47,30)(-10,-3)
\put(25,0){\circle*{3}}
\put(5,-3){\circle*{3}}
\put(5,-3){\line(1,1){12}}
\put(5,-3){\line(-1,1){12}}
\put(-7,9){\circle*{3}}
\put(17,9){\circle*{3}}
\put(23,6){$a$}
\put(10,-6){$b$}
\put(14,12){$d$}
\put(-6,12){$c$}
\end{picture}
 -
\begin{picture}(47,30)(-10,-3)
\put(-3,0){\circle*{3}}
\put(7,0){\circle*{3}}
\put(17,-6){\circle*{3}}
\put(17,6){\circle*{3}}
\put(17,-6){\line(0,1){12}}
\put(20,-8){$a$}
\put(20,6){$b$}
\put(-5,5){$d$}
\put(5,5){$c$}
\end{picture}
-
\begin{picture}(47,50)(-10,-3)
\put(-3,0){\circle*{3}}
\put(22,0){\circle*{3}}
\put(7,-6){\circle*{3}}
\put(7,6){\circle*{3}}
\put(7,-6){\line(0,1){12}}
\put(10,-8){$b$}
\put(10,6){$c$}
\put(-5,5){$d$}
\put(20,5){$a$}
\end{picture}
 -
\begin{picture}(37,50)(-10,-3)
\put(-3,0){\circle*{3}}
\put(22,0){\circle*{3}}
\put(7,-6){\circle*{3}}
\put(7,6){\circle*{3}}
\put(7,-6){\line(0,1){12}}
\put(10,-8){$b$}
\put(10,6){$d$}
\put(-5,5){$c$}
\put(20,5){$a$}
\end{picture}
 + 
\begin{picture}(27,50)(-10,-3)
\put(-3,0){\circle*{3}}
\put(7,0){\circle*{3}}
\put(17,0){\circle*{3}}
\put(27,0){\circle*{3}}
\put(24,5){$a$}
\put(14,5){$b$}
\put(4,5){$c$}
\put(-6,5){$d$}
\end{picture}

\begin{picture}(-800,-200)
\put(-25,-10){\line(0,1){115}}
\put(-25,105){\line(1,0){500}}
\put(-25,-10){\line(1,0){500}}
\put(475,-10){\line(0,1){115}}
\end{picture}
\caption{{\it un calcul d'antipode dans $\Hr$, avec ${\cal D}=\{a,b,c,d\}$. Si on note $s_l$ le sommet d\'ecor\'e par $l$ dans l'arbre choisi,
on a $s_a \leq s_b \leq s_c \leq s_d$.}}
\end{figure}

\section{Alg\`ebre de Frabetti-Brouder}

\subsection{Construction}
\spa On utilise les notations de \cite{Frabetti}.

\begin{defi}
\begin{enumerate}
\item \textnormal{Un arbre binaire est un arbre enracin\'e plan dont chaque
sommet int\'erieur est trivalent. L'ensemble des arbres binaires sera not\'e
$\treesb$. Le degr\'e d'un arbre binaire est le nombre de ses sommets int\'erieurs.
En particulier, l'arbre binaire form\'e uniquement de sa racine est de degr\'e $0$.
Il sera not\'e $\mid$.}

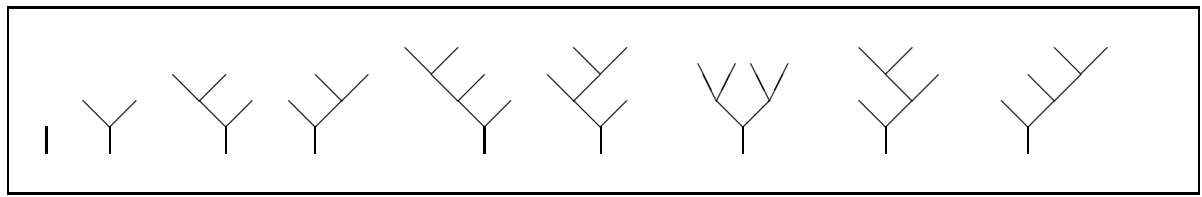
\begin{figure}[h]
\framebox(450,70){
\begin{picture}(20,40)(0,0)
\put(0,0){\line(0,0){10}}
\end{picture}
\begin{picture}(40,40)(0,0)
\put(0,0){\line(0,0){10}}
\put(0,10){\line(1,1){10}}
\put(0,10){\line(-1,1){10}}
\end{picture}
\begin{picture}(30,40)(0,0)
\put(0,0){\line(0,0){10}}
\put(0,10){\line(1,1){10}}
\put(0,10){\line(-1,1){10}}
\put(-10,20){\line(1,1){10}}
\put(-10,20){\line(-1,1){10}}
\end{picture}
\begin{picture}(60,40)(0,0)
\put(0,0){\line(0,0){10}}
\put(0,10){\line(1,1){10}}
\put(0,10){\line(-1,1){10}}
\put(10,20){\line(1,1){10}}
\put(10,20){\line(-1,1){10}}
\end{picture}
\begin{picture}(40,40)(0,0)
\put(0,0){\line(0,0){10}}
\put(0,10){\line(1,1){10}}
\put(0,10){\line(-1,1){10}}
\put(-10,20){\line(1,1){10}}
\put(-10,20){\line(-1,1){10}}
\put(-20,30){\line(1,1){10}}
\put(-20,30){\line(-1,1){10}}
\end{picture}
\begin{picture}(50,40)(0,0)
\put(0,0){\line(0,0){10}}
\put(0,10){\line(1,1){10}}
\put(0,10){\line(-1,1){10}}
\put(-10,20){\line(1,1){10}}
\put(-10,20){\line(-1,1){10}}
\put(0,30){\line(1,1){10}}
\put(0,30){\line(-1,1){10}}
\end{picture}
\begin{picture}(50,40)(0,0)
\put(0,0){\line(0,0){10}}
\put(0,10){\line(1,1){10}}
\put(0,10){\line(-1,1){10}}
\put(10,20){\line(1,2){7}}
\put(10,20){\line(-1,2){7}}
\put(-10,20){\line(1,2){7}}
\put(-10,20){\line(-1,2){7}}
\end{picture}
\begin{picture}(50,40)(0,0)
\put(0,0){\line(0,0){10}}
\put(0,10){\line(1,1){10}}
\put(0,10){\line(-1,1){10}}
\put(10,20){\line(1,1){10}}
\put(10,20){\line(-1,1){10}}
\put(0,30){\line(1,1){10}}
\put(0,30){\line(-1,1){10}}
\end{picture}
\begin{picture}(50,40)(0,0)
\put(0,0){\line(0,0){10}}
\put(0,10){\line(1,1){10}}
\put(0,10){\line(-1,1){10}}
\put(10,20){\line(1,1){10}}
\put(10,20){\line(-1,1){10}}
\put(20,30){\line(1,1){10}}
\put(20,30){\line(-1,1){10}}
\end{picture}
}
\caption{\it les arbres binaires de degr\'e 0,1,2,3.}
\end{figure}

\item \textnormal{Soit $\vee:\treesb \times \treesb \longmapsto
\treesb$ l'application qui greffe deux arbres binaires sur une racine commune.
Alors tout arbre binaire $t\neq \mid$ peut s'\'ecrire $t^l \vee t^r$, avec 
$t^l$ et $t^r$ deux arbres binaires de degr\'e strictement inf\'erieur au degr\'e de
$t$.}
\end{enumerate}
\end{defi}

Soit $H^{\gamma}$ l'espace gradu\'e engendr\'e sur $\mathbb{Q}$ par $\treesb$, muni du produit
d\'efini par les relations de r\'ecurrence :
\begin{eqnarray*}
st &=&(s^lt)\vee s^r \mbox{ o\`u } s=s^l\vee s^r \esp ;\\
\mid t&=&t.
\end{eqnarray*}
(C'est-\`a-dire qu'on greffe $t$ sur la feuille la plus \`a gauche de $s$).\\

Soit $V:\treesb \longmapsto \treesb$, tel que $V(t)=\mid\vee t$.
On munit $H^{\gamma}$ du coproduit d\'efini par les relations de r\'ecurrence :
\begin{eqnarray}
\label{eq18}
\Delta^{\gamma}(\mid)&=&\mid \otimes \mid,\\
\label{eq19}
\Delta^{\gamma}(V(t))&=&V(t)\otimes \mid +(Id \otimes V)\circ
\Delta^{\gamma}(t)-(Id\otimes V)[(V(t^r)\otimes \mid) 
\Delta^{\gamma}(t^l)],\\
\label{eq20}
\Delta^{\gamma}(t\vee s)&=&\Delta^{\gamma}(V(s))\Delta^{\gamma}(t).
\end{eqnarray}
On montre que $H^{\gamma}$ est une alg\`ebre de Hopf gradu\'ee (voir
\cite{Frabetti}).

\subsection{Alg\`ebre de Hopf $\Hfr$}
\begin{prop}
Soit $f:\treesb \longmapsto \sforets$, d\'efinie par r\'ecurrence sur le degr\'e par :
\begin{eqnarray*}
f(\mid)&=&1,\\
f(t^l \vee t^r)&=&B^+\circ f(t^r) f(t^l).
\end{eqnarray*}
Alors $f$ est une bijection, v\'erifiant $f \circ V= B^+\circ f$, et
$f(st)=f(s)f(t)$, $\forall s,t \in \treesb$.
De plus, $poids(f(t))=deg(t)$, $\forall t \in \treesb$.
\end{prop}
{\it Preuve :}  soit $g:\sforets \longmapsto \treesb$ d\'efinie par r\'ecurrence sur le poids :
\begin{eqnarray*}
g(1)&=&\mid,\\
g(tF)&=&g(F)\vee g(B^-(t)), \spa \forall t \in \strees, \esp F \in \sforets.
\end{eqnarray*}
On montre facilement par r\'ecurrence que $g\circ f=Id_{\treesb}$ et $f \circ g=Id_{\sforets}$. Donc $f$ est une bijection.\\

Soit $t \in \treesb$ :
\begin{eqnarray*}
f (V(t))&=&f(\mid \vee t)\\
&=&B^+(f(t)) f(\mid)\\
&=&B^+(f(t)).
\end{eqnarray*}

Soient $s,t \in \treesb$. Montrons que $f(st)=f(s)f(t)$. Supposons d'abord $s$ de la forme $V(s')$. Alors :
\begin{eqnarray*}
f(st)&=&f(t \vee s')\\
&=&B^+ (f(s'))f(t)\\
&=&f(V(s')) f(t)\\
&=&f(s)f(t).
\end{eqnarray*}
Supposons maintenant $s$ quelconque ; alors $s$ s'\'ecrit de mani\`ere unique $s=s_1 \ldots s_n$, avec $s_i =V(s_i')$.
Une r\'ecurrence sur $n$ permet de conclure.

Enfin, les deux propri\'et\'es pr\'ec\'edentes permettent de d\'eduire que $poids(f(t))=deg(t)$ par une r\'ecurrence simple sur le degr\'e. $\Box$
\\

Soit $\Hfr$ l'alg\`ebre librement engendr\'ee sur $\mathbb{Q}$ par les \'el\'ements de $\strees$. Une base de $\Hfr$ est $\sforets$.
Prolongeons $f:H^{\gamma}\longmapsto \Hfr$ par lin\'earit\'e : $f$ devient un isomorphisme d'alg\`ebres gradu\'ees. On munit 
$\Hfr$ d'un coproduit faisant de $f$ un isomorphisme d'alg\`ebres de Hopf. Ce coproduit est not\'e $\Delta_{Fr}$.

Soit $D:\Hfr \longmapsto \Hfr$ et $G:\Hfr \longmapsto \Hfr$  d\'efini par $D(t_1\ldots t_n)=B^-(t_1)$ et $G(t_1\ldots t_n)=t_2 \ldots t_n$, $\forall t_1 \ldots t_n \in \sforets$.
Pour tout $t \in \treesb$, $t \neq \mid$, on a :
\begin{eqnarray*}
D(f(t))&=&D(B^+(f(t^r))f(t^l))\\
&=&f(t^r),\\
G(f(t))&=&G(B^+(f(t^r))f(t^l))\\
&=&f(t^l).
\end{eqnarray*}
$\Delta_{fr}$ peut donc \^etre d\'efini par r\'ecurrence sur le poids en utilisant (\ref{eq18},\ref{eq19},\ref{eq20}) :
\begin{eqnarray}
\label{eq21}
\Delta_{Fr}(1)&=&1\otimes 1,\\
\label{eq22}
\Delta_{Fr} \circ B^+(F)&=&B^+(F) \otimes 1 + (Id \otimes B^+) \circ \Delta_{Fr}(F)\\
\nonumber&& -(Id\otimes B^+)\left[\left(B^+(D(F))\otimes 1\right)\Delta_{Fr}(G(F))\right],\\
\label{eq23}
\Delta_{Fr}(tF)&=&\Delta_{Fr}(t)\Delta_{Fr}(F)
\end{eqnarray}
o\`u $t \in \strees$ et $F\in \sforets$, $F \neq 1$. 

\begin{defi}
\textnormal{Soit $F \in \sforets$, et soit $a$ une ar\^ ete de $F$ ; $a$ est dite {\it ar\^ ete gauche} de $F$ si $a$ est l'ar\^ ete la plus \`a gauche parmi les ar\^etes ayant m\^ eme origine que $a$. Soit $c$ une coupe de $F$ ; $c$ est dite {\it admissible gauche} si $c$ est admissible et ne coupe que des ar\^etes qui ne sont pas gauches. L'ensemble des coupes admissibles gauches de $F$ est not\'e ${\cal A}d^G(F)$ ; l'ensemble des coupes admissibles gauches non vides et non totales de $F$ est not\'e
${\cal A}d_*^G(F)$.}
\end{defi}

\begin{prop}
Pour toute for\^et $F \in \sforets$ :
\begin{eqnarray}
\nonumber \Delta_{Fr}(F)&=&\sum_{c \in {\cal A}d^G(F)} P^c(F) \otimes R^c(F)\\
\label{eq30}&=&F\otimes 1+1 \otimes F+\sum_{c \in {\cal A}d_*^G(F)} P^c(F) \otimes R^c(F).
\end{eqnarray}
\end{prop}
{\it Preuve :} on note $\Delta'_{Fr}:\Hfr\longmapsto \Hfr\otimes \Hfr$ d\'efini par  le second membre de (\ref{eq30}). Il suffit de v\'erifier que $\Delta'_{Fr}$ v\'erifie les \'equations
de r\'ecurrence (\ref{eq21},\ref{eq22},\ref{eq23}). Il est imm\'ediat que (\ref{eq21}) et (\ref{eq23}) sont v\'erifi\'ees. 
Soit $F=t_1\ldots t_n \in \forets$.

 On a une bijection $\alpha:{\cal A}d(F) \longmapsto {\cal A}d(B^+(F))-\{\mbox{coupe totale de $B^+(F)$}\}$ telle que
$P^{\alpha(c)}(B^+(F))=P^c(F)$ et $R^{\alpha(c)}(B^+(F))=B^+(P^c(F))$.  Soit $c \in {\cal A}d^G(F)$. 
Deux cas se pr\'esentent :
\begin{enumerate}
\item Si $c_{\mid t_1}$ n'est pas la coupe totale de $t_1$, alors $\alpha(c) \in {\cal A}d^G(B^+(F))$, et $\alpha(c)$ n'est pas la coupe totale de $B^+(F)$. 
De plus, toute coupe dans ${\cal A}d^G(B^+(F))$ non totale  est atteinte.
\item  Sinon, $c$ ne coupe que des ar\^etes non gauches et l'ar\^ete $a$ menant de la racine de $B^+(F)$ vers la racine de $t_1$ (cette ar\^ete est gauche). Soit  
$c'=c_{\mid t_2 \ldots t_n}$. C'est une  coupe admissible gauche de $t_2 \ldots t_n$. On a alors  
$P^{\alpha(c)}(B^+(F))=t_1 P^{c'}(t_2\ldots t_n)$ et $R^{\alpha(c)}(B^+(F))=B^+( P^{c'}(t_2\ldots t_n))$. De plus, toute coupe de ${\cal A}d(B^+(F))$ de cette forme est atteinte.
 \end{enumerate}
On a donc :
\begin{eqnarray*}
 (Id \otimes B^+) \circ \Delta'_{Fr}(F)&=&(Id \otimes B^+) (\sum_{c \in {\cal A}d^G(F)} P^c(F) \otimes R^c(F))\\
&=&\sum_{c \in {\cal A}d^G(F)} P^{\alpha(c)}(B^+(F)) \otimes R^{\alpha(c)}(B^+(F)) \\   
&=&\sum_{c \in {\cal A}d^G(B^+(F))} P^c(B^+(F))\otimes R^c(B^+(F)) - B^+(F) \otimes 1\\ 
&&+ \sum_{c' \in {\cal A}d^G(t_2 \ldots t_n)} t_1 P^{c'}(t_2\ldots t_n)\otimes B^+(R^{c'}(t_2\ldots t_n)) \\
&=&\Delta'_{Fr}(B^+(F)) -B^+(F)\otimes 1 \\
&&+ (Id\otimes B^+)[(t_1 \otimes 1) \Delta'_{Fr}(t_2\ldots t_n)]\\
&=&\Delta'_{Fr}(B^+(F)) -B^+(F)\otimes 1 \\
&&+(Id\otimes B^+)\left[\left(B^+(D(F))\otimes 1\right)\Delta'_{Fr}(G(F))\right].
\end{eqnarray*}
Et donc (\ref{eq22})  est v\'erifi\'ee. $\Box$

\begin{prop}
$(\Hfr,m, \eta,\Delta_{Fr},\varepsilon,S_{Fr})$ est une alg\`ebre de Hopf gradu\'ee isomorphe \`a $H^{\gamma}$ ($m, \eta, \varepsilon$ et la gradu\-ation \'etant les m\^ emes que dans $\sHr$).
\end{prop}

\begin{figure}[h]
\framebox(450,180){
\begin{picture}(0,0)(-20,-60)
\begin{picture}(0,0)
\put(-20,15){\it Coupes admissibles}
\put(0,0){\it gauches :}
\end{picture}
\begin{picture}(40,50)(-90,0)
\put(2,10){\line(1,-1){10}}
\put(0,0){\circle*{2}}
\put(0,0){\line(0,1){10}}
\put(0,0){\line(1,1){10}}
\put(0,0){\line(-1,1){10}}
\put(-10,10){\circle*{2}}
\put(0,10){\circle*{2}}
\put(10,10){\circle*{2}}
\put(0,10){\line(0,1){10}}
\put(0,20){\line(-1,1){10}}
\put(0,20){\line(1,1){10}}
\put(0,20){\circle*{2}}
\put(-10,30){\circle*{2}}
\put(10,30){\circle*{2}}
\end{picture}
\begin{picture}(40,50)(-120,0)
\put(-5,7){\line(1,0){10}}
\put(0,0){\circle*{2}}
\put(0,0){\line(0,1){10}}
\put(0,0){\line(1,1){10}}
\put(0,0){\line(-1,1){10}}
\put(-10,10){\circle*{2}}
\put(0,10){\circle*{2}}
\put(10,10){\circle*{2}}
\put(0,10){\line(0,1){10}}
\put(0,20){\line(-1,1){10}}
\put(0,20){\line(1,1){10}}
\put(0,20){\circle*{2}}
\put(-10,30){\circle*{2}}
\put(10,30){\circle*{2}}
\end{picture}
\begin{picture}(40,50)(-140,0)
\put(2,30){\line(1,-1){10}}
\put(0,0){\circle*{2}}
\put(0,0){\line(0,1){10}}
\put(0,0){\line(1,1){10}}
\put(0,0){\line(-1,1){10}}
\put(-10,10){\circle*{2}}
\put(0,10){\circle*{2}}
\put(10,10){\circle*{2}}
\put(0,10){\line(0,1){10}}
\put(0,20){\line(-1,1){10}}
\put(0,20){\line(1,1){10}}
\put(0,20){\circle*{2}}
\put(-10,30){\circle*{2}}
\put(10,30){\circle*{2}}
\end{picture}
\begin{picture}(40,50)(-160,0)
\put(2,10){\line(1,-1){10}}
\put(-5,7){\line(1,0){10}}
\put(0,0){\circle*{2}}
\put(0,0){\line(0,1){10}}
\put(0,0){\line(1,1){10}}
\put(0,0){\line(-1,1){10}}
\put(-10,10){\circle*{2}}
\put(0,10){\circle*{2}}
\put(10,10){\circle*{2}}
\put(0,10){\line(0,1){10}}
\put(0,20){\line(-1,1){10}}
\put(0,20){\line(1,1){10}}
\put(0,20){\circle*{2}}
\put(-10,30){\circle*{2}}
\put(10,30){\circle*{2}}
\end{picture}
\begin{picture}(40,50)(-180,0)
\put(2,10){\line(1,-1){10}}
\put(2,30){\line(1,-1){10}}
\put(0,0){\circle*{2}}
\put(0,0){\line(0,1){10}}
\put(0,0){\line(1,1){10}}
\put(0,0){\line(-1,1){10}}
\put(-10,10){\circle*{2}}
\put(0,10){\circle*{2}}
\put(10,10){\circle*{2}}
\put(0,10){\line(0,1){10}}
\put(0,20){\line(-1,1){10}}
\put(0,20){\line(1,1){10}}
\put(0,20){\circle*{2}}
\put(-10,30){\circle*{2}}
\put(10,30){\circle*{2}}
\end{picture}
\end{picture}
\begin{picture}(420,50)
$\Delta_{Fr}${\huge(}
\begin{picture}(25,50)(-10,0)
\put(0,0){\circle*{2}}
\put(0,0){\line(0,1){10}}
\put(0,0){\line(1,1){10}}
\put(0,0){\line(-1,1){10}}
\put(-10,10){\circle*{2}}
\put(0,10){\circle*{2}}
\put(10,10){\circle*{2}}
\put(0,10){\line(0,1){10}}
\put(0,20){\line(-1,1){10}}
\put(0,20){\line(1,1){10}}
\put(0,20){\circle*{2}}
\put(-10,30){\circle*{2}}
\put(10,30){\circle*{2}}
\end{picture}
{\huge)}
$=$
\begin{picture}(7,5)(0,0)
\put(0,0){\circle*{2}}
\end{picture}
$\otimes$
\begin{picture}(20,50)(-10,0)
\put(0,0){\circle*{2}}
\put(0,0){\line(0,1){10}}
\put(0,0){\line(-1,1){10}}
\put(-10,10){\circle*{2}}
\put(0,10){\circle*{2}}
\put(0,10){\line(0,1){10}}
\put(0,20){\line(-1,1){10}}
\put(0,20){\line(1,1){10}}
\put(0,20){\circle*{2}}
\put(-10,30){\circle*{2}}
\put(10,30){\circle*{2}}
\end{picture}
$+$
\begin{picture}(20,50)(-10,0)
\put(0,0){\circle*{2}}
\put(0,0){\line(0,1){10}}
\put(0,10){\line(-1,1){10}}
\put(0,10){\line(1,1){10}}
\put(0,10){\circle*{2}}
\put(-10,20){\circle*{2}}
\put(10,20){\circle*{2}}
\end{picture}
$\otimes$
\begin{picture}(25,50)(-10,0)
\put(0,0){\circle*{2}}
\put(0,0){\line(1,1){10}}
\put(0,0){\line(-1,1){10}}
\put(-10,10){\circle*{2}}
\put(10,10){\circle*{2}}
\end{picture}
$+$
\begin{picture}(7,5)(0,0)
\put(0,0){\circle*{2}}
\end{picture}
$\otimes$
\begin{picture}(25,50)(-10,0)
\put(0,0){\circle*{2}}
\put(0,0){\line(0,1){10}}
\put(0,0){\line(-1,1){10}}
\put(0,0){\line(1,1){10}}
\put(10,10){\circle*{2}}
\put(-10,10){\circle*{2}}
\put(0,10){\circle*{2}}
\put(0,10){\line(0,1){10}}
\put(0,20){\line(0,1){10}}
\put(0,20){\circle*{2}}
\put(0,30){\circle*{2}}
\end{picture}
$+$
\begin{picture}(24,50)(-10,0)
\put(0,0){\circle*{2}}
\put(0,0){\line(0,1){10}}
\put(0,10){\line(-1,1){10}}
\put(0,10){\line(1,1){10}}
\put(0,10){\circle*{2}}
\put(-10,20){\circle*{2}}
\put(10,20){\circle*{2}}
\put(12,0){\circle*{2}}
\end{picture}
$\otimes$
\begin{picture}(10,50)(0,0)
\put(0,0){\circle*{2}}
\put(0,0){\line(0,1){10}}
\put(0,10){\circle*{2}}
\end{picture}
$+$
\begin{picture}(7,5)(0,0)
\put(0,0){\circle*{2}}
\put(7,0){\circle*{2}}
\end{picture}
$\otimes$
\begin{picture}(25,50)(-10,0)
\put(0,0){\circle*{2}}
\put(0,0){\line(0,1){10}}
\put(0,0){\line(-1,1){10}}
\put(-10,10){\circle*{2}}
\put(0,10){\circle*{2}}
\put(0,10){\line(0,1){10}}
\put(0,20){\line(0,1){10}}
\put(0,20){\circle*{2}}
\put(0,30){\circle*{2}}
\end{picture}
\end{picture}
\begin{picture}(0,0)(270,50)
$+$\begin{picture}(25,60)(-10,0)
\put(0,0){\circle*{2}}
\put(0,0){\line(0,1){10}}
\put(0,0){\line(1,1){10}}
\put(0,0){\line(-1,1){10}}
\put(-10,10){\circle*{2}}
\put(0,10){\circle*{2}}
\put(10,10){\circle*{2}}
\put(0,10){\line(0,1){10}}
\put(0,20){\line(-1,1){10}}
\put(0,20){\line(1,1){10}}
\put(0,20){\circle*{2}}
\put(-10,30){\circle*{2}}
\put(10,30){\circle*{2}}
\end{picture}
$\otimes1+1\otimes$
\begin{picture}(25,50)(-10,0)
\put(0,0){\circle*{2}}
\put(0,0){\line(0,1){10}}
\put(0,0){\line(1,1){10}}
\put(0,0){\line(-1,1){10}}
\put(-10,10){\circle*{2}}
\put(0,10){\circle*{2}}
\put(10,10){\circle*{2}}
\put(0,10){\line(0,1){10}}
\put(0,20){\line(-1,1){10}}
\put(0,20){\line(1,1){10}}
\put(0,20){\circle*{2}}
\put(-10,30){\circle*{2}}
\put(10,30){\circle*{2}}
\end{picture}.
\end{picture}
}
\caption{\it un calcul de coproduit dans $\Hfr$.}
\end{figure}

\subsection{Isomorphisme entre $\sHr$ et $\Hfr^{*g}$}

\spa Soit $F=t_1 \ldots t_n  \in \sforets$. Posons $t_n=B^+(s_1 \ldots s_m)$, $s_1 \ldots s_m \in \strees$. Soient $r(F)=m$ (nombre d'ar\^ etes ayant pour origine la racine de l'arbre le plus \`a droite de $F$), et $p(F)$ le nombre de $t_i$ \'egaux \`a $\bullet$. En particulier, si $F=1$, on a $r(F)=p(F)=0$, et si $t_n=\bullet$, alors $r(n)=0$.

Soit $x$ un \'el\'ement de $\Hfr$. Soit $x=\sum a_FF$ sa d\'ecomposition dans la base des for\^ets. On note ${\cal F}(x)=\{F \in \sforets, a_F \neq 0\}$. On note
$p(x)=\max\{p(F)/ F \in {\cal F}(x)\}$, et ${\cal P}(x)=\{ F \in {\cal F}(x), p(F)=p(x)\}$. Si $x \neq 0$, ${\cal P}(x)$ est non vide.

\begin{lemme}
Soit $x$ un \'el\'ement primitif non nul de $\Hfr$. Alors : soit $\bullet \in {\cal P}(x)$, soit 
il existe $F\in {\cal P}(x)$, tel que $r(F)=1$ (ces deux cas ne s'excluant pas mutuellement).
\end{lemme}
{\it Preuve :} soit $F=t_1 \ldots t_n \in {\cal P}(x)$. Si $F=\bullet$, c'est termin\'e. Supposons $F \neq \bullet$.

Supposons que $F=\bullet^i$, $i\geq 2$. Alors $F$ est la seule for\^et ayant une coupe admissible gauche $c$ telle que
$P^c(F)=\bullet^{i-1}$, $R^{c}(F)=\bullet$. Donc $\bullet^{i-1}\otimes \bullet$ doit apparaitre dans l'\'ecriture de $\Delta_{Fr}(x)$ dans la base des for\^ets avec un coefficient non nul, et donc $x$ n'est pas primitif : contradiction.

Par suite, $F$ est de la forme $F=Gt\bullet^i$, avec $G \in \sforets$, $t \in \strees$, $t \neq \bullet$. Parmi tous les \'el\'ements de ${\cal P}(x)$ de cette forme, choisissons $F$ de sorte que $i$ soit minimal. Supposons $i\geq 1$. $Gt \otimes \bullet^i$ apparait dans l'\'ecriture de $\Delta_{Fr}(F)$. Comme $x$ est primitif, il existe $F' \in {\cal F}(x)$, $F' \neq F$, poss\'edant une coupe admissible gauche $c$ telle que $P^c(F')\otimes R^c(F')=Gt\otimes \bullet^i$. Or une telle for\^et est de la forme
$Ht\bullet^j$, avec $j <i$. On aboutit \`a une contradiction avec le choix de $F$, et donc $i=0$.

On a donc trouv\'e $F=t_1 \ldots t_n\in {\cal P}(x)$, telle que $t_n \neq \bullet$. Parmi toutes les $F$ de cette sorte, choisissons $F$ de sorte que :
\begin{enumerate}
\item si $F'=t'_1 \ldots t'_m \in {\cal P}(x)$, $t'_m \neq \bullet$, alors $poids(t'_m)\geq poids(t_n)$ ;
\item si de plus $poids(t'_m)=poids(t_n)$, alors $r(F')\geq r(F)$.
\end{enumerate}
Supposons $r(F)>1$ : $t_n =B^+(s_1\ldots s_m)$, $m\geq 2$. Alors $s_2 \ldots s_m \otimes t_1 \ldots t_{n-1} B^+(s_1)$ apparait dans l'\'ecriture de $\Delta_{Fr}(F)$.
Comme $x$ est primitif, il existe $F' \in {\cal F}(x)$, $F'\neq F$, poss\'edant une coupe admissible gauche $c$ telle que $P^c(F')\otimes p^c(F')=s_2 \ldots s_m \otimes t_1 \ldots t_{n-1} B^+(s_1)$. $F'$ est donc obtenue en greffant les diff\'erents arbres de $s_2 \ldots s_m$ sur les diff\'erents arbres de   $t_1 \ldots t_{n-1} B^+(s_1)$, ou en les intercalant entre ces arbres. Comme $c$ est admissible gauche, on ne peut pas greffer $s_1, \ldots, s_m$, sur des arbres de $t_1 \ldots t_{n-1} B^+(s_1)$ \'egaux \`a $\bullet$, donc n\'ecessairement, $p(F')\geq p(F)$, et donc $F' \in {\cal P}(x)$. Posons $F'=t'_1 \ldots t'_m$.  Trois cas se pr\'esentent :
\begin{enumerate}
\item soit on ne greffe pas tous les $s_i$ sur $B^+(s_1)$ : alors $t'_m$ est l'un des $s_i$ ou $B^+(s_1)$, et donc $poids(t'_m)<poids(t_n)$. 
Si $t'_m=\bullet$, alors $p(F')>p(F)=p(x)$, et donc $a_{F'}=0$. Sinon, par choix de $F$ (condition 1), on a alors $a_{F'}=0$, et donc dans les deux cas $F' \notin {\cal F}(x)$.
\item soit on  greffe tous les $s_i$ sur $B^+(s_1)$, mais pas tous sur la racine de $B^+(s_1)$ : alors  $poids(t'_m)=poids(t_n)$, et $r(F')<r(F)$, donc $a_{F'}=0$ (condition 2).
\item soit on greffe tous les $s_i$ sur la racine de $B^+(s_1)$ : alors on doit  n\'ecessairement greffer $s_2\ldots  s_m$  \`a gauche de $s_1$  (car $c $ est admissible gauche) et dans cet ordre, et donc $F'=F$.
\end{enumerate}
Dans les trois cas on  aboutit \`a une contradiction, et donc $r(F)=1$, ce que d\'emontre le lemme.  $\Box$
\\

Soit $\beta:\Hfr \longmapsto \Hfr$ d\'efinie par :
\begin{eqnarray*}
\beta(1)&=&0,\\
\beta(t_1 \ldots t_{n-1} \bullet)&=&  t_1 \ldots t_{n-1}, \\
\beta(t_1 \ldots t_n)&=& t_1 \ldots t_{n-1} B^-(t_n) \mbox{ si $r(t_1 \ldots t_n) = 1$, }\\
\beta(t_1 \ldots t_n)&=& 0 \mbox{ si $r(t_1 \ldots t_n) \geq 2$. }
\end{eqnarray*}

\begin{prop}
\begin{enumerate}
\item $\beta$ v\'erifie les propri\'et\'es suivantes :
\begin{enumerate}
\item $\beta$ est homog\`ene de degr\'e $-1$.
\item $\forall x,y \in \Hfr$, $\beta(xy)=\beta(x) \varepsilon(y)+x \beta(y)$.
\item $\beta_{\mid prim(\Hfr)}$ est injective.
\end{enumerate}
\item $\beta^{*g}: \Hfr^{*g} \longmapsto \Hfr^{*g} $ v\'erifie les propri\'et\'es suivantes :
\begin{enumerate}
\item $\beta^{*g}$ est homog\`ene de degr\'e $+1$.
\item $\forall f \in \Hfr^{*g}$, $\Delta(\beta^{*g}(f))= \beta^{*g}(f) \otimes 1+( Id \otimes B^{*g})\circ \Delta(f)$ o\`u $\Delta$ d\'esigne le coproduit de $\Hfr^{*g}$.
\item Soit $M$ l'id\'eal d'augmentation de $\Hfr^{*g}$ ; alors $\Hfr^{*g}=Im(\beta^{*g})+(1 \oplus M^2)$.
\end{enumerate}
\end{enumerate}
\end{prop}
{\it Preuve :} \begin{enumerate}
\item \begin{enumerate}
\item Evident, avec la d\'efinition de $\beta$.
\item On remarque que $\beta(t_1\ldots t_n)=t_1\ldots t_{n-1} \beta(t_n)$ ; le r\'esultat est alors imm\'ediat.
\item Soit $x  \in prim(\Hfr)$, non nul. Supposons $\beta(x)=0$. Si $\bullet \in {\cal P}(x)$, alors $\varepsilon(\beta(x))=a_{\bullet} \neq 0$ : contradiction.
Donc d'apr\`es le lemme pr\'ec\'edent,  il existe $F=t_1 \ldots t_n \in {\cal P}(x)$, $r(F)=1$. Alors $\beta(F)=t_1 \ldots t_{n-1} B^-(t_n)$, et $B^-(t_n) \in \strees$. Comme $\beta(x)=0$, il existe $F' \in {\cal F}(x)$, $F' \neq F$, telle que $\beta(F')=\beta(F)$. Alors n\'ecessairement $F'=t_1 \ldots t_{n-1} B^-(t_n)\bullet$, et donc $p(F')\geq p(F)+1>p(x)$ : on ne peut donc avoir $F' \in {\cal F}(x)$ : contradiction. Donc $\beta(x)\neq 0$.
\end{enumerate}
\item
\begin{enumerate}
\item  D\'ecoule du lemme \ref{lemme2}.
\item Soient $f\in \Hfr^{*g}$, $x,y \in \Hfr$.
\begin{eqnarray*}
(\Delta(\beta^{*g}(f)), x \otimes y)&=&(\beta^{*g}(f), xy)\\
&=&(f, \beta(xy))\\
&=&(f, \beta(x) \varepsilon(y)+x \beta(y))\\
&=&(\beta^{*g}(f),x)(1,y)+((Id \otimes \beta^{*g})\circ \Delta(f),x \otimes y).
\end{eqnarray*}
Le r\'esultat en d\'ecoule imm\'ediatement.
\item A l'aide de la proposition \ref{prop6}-2,  on identifie $Prim(\Hfr)^{*g}$ et $\Hfr^{*g}/ ((1) \oplus M^2)).$ Alors si $i:Prim(\Hfr)^{*g} \longmapsto \Hfr^{*g}$ est l'injection canonique,  $i^{*g}:\Hfr^{*g} \longmapsto \Hfr^{*g}/((1)\oplus M^2)$ est la surjection canonique. D'apr\`es 1.(c), $\beta \circ i$ est injective. Donc $i^{*g} \circ \beta^{*g}$ est surjective. On a donc :
\begin{eqnarray*}
Im(i^{*g} \circ \beta^{*g})&=& \frac{Im(\beta^{*g})}{(1)\oplus M^2}\\
&=& \frac{\Hfr^{*g} }{(1)\oplus M^2},
\end{eqnarray*}
d'o\`u le r\'esultat. $\Box$
\end{enumerate}
\end{enumerate}

Soit $\Phi: \sHr \longmapsto \Hfr^{*g}$ l'unique morphisme d'alg\`ebres de Hopf v\'erifiant $\Phi \circ B^+ = \beta^{*g} \circ \Phi$ (propri\'et\'e universelle de $\sHr$).
Montrons que $\Phi$ est homog\`ene de degr\'e z\'ero : soit $F\in \sforets$, de poids $n$, montrons que $\Phi(F)$ est homog\`ene de poids $n$ par r\'ecurrence sur $n$.
Si $n=0$, alors $F=1$, $\Phi(F)=1$. Supposons la propri\'et\'e vraie pour toute for\^ et de poids $n'<n$ : si $F \notin \strees$, alors il existe $F_1,F_2$, de poids strictement inf\'erieur \`a $n$, telles que $F=F_1F_2$. Par suite $\Phi(F)=\Phi(F_1)\Phi(F_2)$, et donc  $\Phi(F)$ est homog\`ene de poids $n$. Sinon, il existe $F'$ de poids $n-1$, telle que $F=B^+(F')$. Alors $\Phi(F)=\beta^{*g}(\Phi(F'))$ : comme $\beta^{*g}$ est homog\`ene de degr\'e 1, $\Phi(F)$ est homog\`ene de poids $n-1+1=n$.

Montrons que $\Phi$ est surjective : soit $y \in \Hfr^{*g}$, homog\`ene de poids $n$, montrons que $y \in Im(\Phi)$ par r\'ecurrence sur $n$. C'est \'evident si $n=0$. Sinon, 
$y \in (Im(\beta^{*g}) + M^2)\cap {\cal H}_n^* = \beta^{*g}({\cal H}_{n-1}^*) + M^2 \cap {\cal H}_n^*$. On peut donc se ramener \`a $y$ de la forme $\beta^{*g}(y')$, $y'$ homog\`ene de poids $n-1$, ou $y=y_1y_2$, $y_i \in M$, homog\`enes de poids strictement inf\'erieur \`a $n$. Dans le premier cas, il existe $x' \in \sHr $, $\Phi(x')=y'$, et alors $\Phi(B^+(x'))=\beta^{*g}\circ \Phi(x')= \beta^{*g}(y')=y$. Dans le deuxi\`eme cas, il existe $x_1,x_2 \in \sHr$, $\Phi(x_i)=y_i$, et donc $y=\Phi(x_1x_2)$.

Par suite, $\Phi$ \'etant homog\`ene de degr\'e z\'ero et surjectif, et les composantes homog\`enes de $\sHr$ et $\Hfr$ ayant les m\^ emes dimensions finies, $\Phi$ est un isomorphisme.

\begin{prop}
l'unique morphisme d'alg\`ebres de Hopf $\Phi:\sHr \longmapsto \Hfr^{*g}$ tel que $\Phi \circ B^+ =\beta^{*g} \circ \Phi$ est un isomorphisme d'alg\`ebres de Hopf gradu\'ees.
\end{prop}

\begin{theo}
 $\Hfr$ et $\sHr$ sont des alg\`ebres de Hopf gradu\'ees isomorphes.
\end{theo}
{\it Preuve :} $\Phi^{*g}:(\Hfr^{*g})^{*g}\longmapsto \sHr^{*g}$ est un isomorphisme d'alg\`ebres de Hopf gradu\'ees. De plus, $(\Hfr^{*g})^{*g}$ est canoniquement isomorphe \`a $\Hfr$ comme alg\`ebre de Hopf gradu\'ee, et $\sHr$ est isomorphe \`a $\sHr^{*g}$ comme alg\`ebre de Hopf gradu\'ee d'apr\`es le th\'eor\`eme \ref{theo27}.

\section{Sous-alg\`ebres des diff\'eomorphismes formels} 
\spa il s'agit dans cette  partie de mettre en \'evidence une  sous-alg\`ebre de Hopf de $\sHr$ et de $\Hfr$ dont l'abelianis\'ee est isomorphe \`a l'alg\`ebre de Connes-Moscovici ${\cal H}_{CM}$.

\subsection{Rappels et compl\'ements sur l'alg\`ebre de Hopf ${\cal H}_{CM}$}

\spa (Voir \cite{Moscovici,Connes}). Soit ${\cal H}_T$ l'alg\`ebre d\'efinie par les g\'en\'erateurs $X,Y,\delta_n$, $n \geq 1$, et les relations :
$$[Y,X]=X, \esp [Y,\delta_n]=n \delta_n,\esp  [\delta_n,\delta_m]=0,\esp [X,\delta_n]=\delta_{n+1}.$$
On la munit d'un coproduit donn\'e par :
\begin{eqnarray*}
\Delta(Y)&=&Y\otimes 1 +1 \otimes Y,\\
\Delta(X)&=&X\otimes 1+1 \otimes X+\delta_1 \otimes Y,\\
\Delta(\delta_1)&=&\delta_1 \otimes 1+1\otimes \delta_1.
\end{eqnarray*}
La sous-alg\`ebre de ${\cal H}_{T}$ engendr\'ee par les $\delta_n$, $n\geq 1$, est une sous-alg\`ebre de Hopf not\'ee ${\cal H}_{CM}$.
Elle est gradu\'ee en posant $poids(\delta_n)=n$.\\

${\cal H}_{CM}$ est une alg\`ebre de Hopf gradu\'ee commutative, et donc son dual gradu\'e $({\cal H}_{CM})^{*g}$  est l'alg\`ebre enveloppante d'une alg\`ebre de Lie $\g_{CM}$.
Par la proposition  \ref{prop6}, une base de $\g_{CM}$ est $(Z_n)_{n \geq 1}$, d\'efinie par :
\begin{eqnarray*}
Z_n(\delta_{i_1} \ldots \delta_{i_n})&=&0 \mbox{ si } n\neq 1,\\
Z_n(\delta_m)&=&\delta_{n,m}.
\end{eqnarray*}
On peut montrer la formule suivante :
\begin{equation}
\label{eq25}
[Z_n,Z_m]=(m-n)Z_{n+m}.
\end{equation}
\begin{prop}
Une base de $Prim({\cal H}_{CM})$ est donn\'ee par $(\delta_1,2\delta_2-\delta_1^2)$.
\end{prop}
{\it Preuve :}  soit $(Z_1^{\alpha_1} \ldots Z_k ^{\alpha_k})_{\alpha_1,\ldots \alpha_k}$ une base de Poincar\'e-Birkhoff-Witt de $({\cal H}_{CM})^{*g}$.
Si $\alpha_1+\ldots +\alpha_n\neq 1$, alors $Z_1^{\alpha_1} \ldots Z_k ^{\alpha_k} \in (1)+M_*^2$, o\`u $M_*$ est l'id\'eal d'augmentation de $({\cal H}_{CM})^{*g}$.
De plus, si $n\geq 3$, on a alors :
$$ Z_n = \frac{1}{n-2} [Z_1,Z_{n-1}],$$
et donc $Z_n \in M_*^2$. Par suite, $vect (Z_1,Z_2)$ est un compl\'ementaire de $(1)+M_*^2$ dans  $({\cal H}_{CM})^{*g}$, et donc $dim(({\cal H}_{CM})^{*g}/((1)+M_*^2))) \leq 2$.

Par suite, $Prim({\cal H}_{CM})=(({\cal H}_{CM})^{*g}/((1)+M_*^2)))^{*g}$ est de dimension au plus \'egale \`a 2.
On v\'erifie facilement que les deux \'el\'ements $\delta_1$ et $2\delta_2 -\delta_1^2$ sont primitifs ; ils sont lin\'eairement ind\'ependant car de poids diff\'erents, et donc forment une base de   $Prim({\cal H}_{CM})$. $\Box$

\subsection{Angles, greffes et coupes}

\spa On rappelle la notion d'angle d'un arbre enracin\'e d\'efinie par Kontsevich dans \cite{Konts} et utilis\'ee par Chapoton dans \cite{Chapoton} : 
\begin{defi}
\textnormal{Soit $t \in \strees$. Supposons $t$ dessin\'e dans le demi-disque sup\'erieur ouvert
$$D_+=\{(x,y)\in \mathbb{R}^2/y>0, \esp x^2+y^2<1\},$$
sauf la racine plac\'ee en $(0,0)$. On appelle angle de $t$ un couple $(s,\alpha)$ o\`u $s$ est un sommet de $t$ et $\alpha $ une composante connexe de $B_{\epsilon}(s) \cap( D_+-t)$ o\`u $B_{\epsilon}(s)$ est un petit disque de centre $s$. On note $Angles(t)$ l'ensemble des angles de $t$.}

\textnormal{$Angles(t)$ est muni d'une relation d'ordre totale de gauche \`a droite de la mani\`ere suivante : en consid\'erant chaque angle de $t$ comme une direction issue d'un sommet , on peut tracer un chemin de chaque angle vers un point du cercle unit\'e, de sorte que ces chemins ne se coupent pas. On obtient alors un point du demi-cercle associ\'e \`a chaque angle. L'ordre de ces points de gauche \`a droite d\'etermine l'ordre total sur les angles.}
\end{defi}


\begin{defi}
\textnormal{Soit $F=t_1 \ldots t_m\in \sforets$, $t \in \strees$. Une greffe de $F$ sur $t$ est une suite croissante de $m$ angles de $t$. Le r\'esultat d'une greffe $g = ((s_1,\alpha_1), \ldots ,(s_m, \alpha_m))$ de $F$ sur $t$ est l'arbre not\'e $R_g(F,t)$ obtenu en greffant $t_i$ sur le sommet $s_i$ dans la composante connexe $\alpha_i$ de $B_{\epsilon}(s_i) \cap( D_+-t)$, avec la condition suivante : si $(s_i, \alpha_i)=(s_{i+1}, \alpha_{i+1})$, alors $t_{i}$ est greff\'e \`a gauche de $t_{i+1}$.} 
\end{defi}


\begin{prop}
Soient $F \in \sforets$, $t\in \strees$. On pose :
\begin{eqnarray*}
G_{F,t}&=&\{ \makebox{greffes de $F$ sur $t$}\} ,\\
C_{F,t}&=& \bigcup_{t' \in \strees} \{c \makebox{ coupe admissible de $t'$ telle que $P^{c}(t')=F$ et $R^c(t')=t$}\}.
\end{eqnarray*}
Soit $f_1:C_{F,t} \longmapsto G_{F,t}$, qui \`a une coupe $c$ de l'arbre $t'$ associe l'unique greffe $g$ de $F$ sur $t$, telle que $R_g(F,t)=t'$, les ar\^etes cr\'e\'ees lors de la greffe \'etant les ar\^etes de $t'$ sur lesquelles agit $c$.\\
Soit $f_2:G_{F,t} \longmapsto C_{F,t}$, qui \`a une greffe $g$ de $F$ sur $t$ associe la coupe $c$ de $t'=R_g(F,t)$, portant sur les ar\^etes cr\'e\'ees lors de la greffe.\\
Alors $f_1$ et $f_2$ sont des bijections r\'eciproques l'une de l'autre. 
\end{prop}
{\it Preuve :} par construction de $R_g(F,t)$,  la coupe $f_2(g)$ est bien admissible. Elle v\'erifie de plus $P^{f_2(g)}(R_g(F,t))=F$,  $R^{f_2(g)}(R_g(F,t))=t$.  Donc $f_2$ est bien d\'efinie. 
Le reste est imm\'ediat. $\Box$
\\

Soient $F,G,H \in \sforets$. Rappelons que $n(F,G;H)$ est le nombre de coupes admissibles $c$ de $H$ telles que $P^c(H)=F$, $R^c(H)=G$. On pose  $n_G(F,G;H)$  le nombre de coupes admissibles gauches $c$ de $H$ telles que $P^c(H)=F$, $R^c(H)=G$.

\begin{cor}
\label{cor52} Soit $t_1,t_2\in \trees$. On a dans $\sHr$  :
$$ [ e_{t_1}, e_{t_2}]= \sum_{\makebox{$g$ greffe de $t_2$ sur $t_1$}} e_{R_{g}(t_1,t_2)} -  \sum_{\makebox{$g$ greffe de $t_1$ sur $t_2$}} e_{R_{g}(t_2,t_1)}.$$
\end{cor}
{\it Preuve :}
\begin{eqnarray*}
[ e_{t_1}, e_{t_2}]&=& \sum_{t' \in \strees} n(t_1,t_2; t') e_{t'} -  \sum_{t' \in \strees} n(t_2,t_1; t') e_{t'}\\
&=& \sum_{\stackrel{c \in C_{t_1,t_2}}{c \makebox{\footnotesize{ coupe de }}t'}} e_{t'} - \sum_{\stackrel{c \in C_{t_2,t_1}}{c \makebox{\footnotesize{ coupe de }}t'}} e_{t'}\\
&=& \sum_{g \in G_{t_1,t_2}} e_{R_g(t_1,t_2)}- \sum_{g \in G_{t_2,t_1}} e_{R_g(t_2,t_1)}.
\end{eqnarray*}  
(On a utilis\'e la proposition \ref{prop29} pour la premi\`ere \'egalit\'e, et les propri\'et\'es de $f_1$ pour la troisi\`eme). $\Box$
\begin{cor}
\label{cor60} Soit $F=t_1 \ldots t_m \in \sforets$, $t \in \strees$ de poids $n$.
\begin{eqnarray*}
\sum_{t' \in \strees} n(F,t;t') &=& \binom{2n+m-2}{m}\\
\sum_{t' \in \strees} n_G(F,t;t') &=& \binom{n+m-2}{m}.
\end{eqnarray*}
\end{cor}
{\it Preuve :} on a $\sum_{t' \in \strees} n(F,t;t') = card(C_{F,t})=card(G_{F,t})$, car $f_1$ est bijective. Or en notant $a=card(Angles(t))$, le cardinal de $G_{F,t}$ est le nombre de suites croissantes de $m$ \'el\'ements choisis parmi $a$ : d'apr\`es un lemme classique,
$$\sum_{t' \in \strees} n(F,t;t') = \binom{m+a-1}{m}.$$

Calculons $a$. Chaque sommet $s$ de $t$ est le sommet de $f(s)+1$ angles, o\`u $f(s)$ est le nombre d'ar\^ etes issues de $s$. Donc :
\begin{eqnarray*}
a&=&\sum_{s \in som(t)} f(s)+1\\
&=& \left(\sum_{s \in som(t)} f(s) \right)+n\\
&=&\left(\makebox{nombre d'ar\^ etes de $t$} \right)+n\\
&=&n-1+n\\
&=&2n-1
\end{eqnarray*}
ce qui prouve le premier r\'esutat.

Consid\'erons : 
$$C_{F,t}^G=\bigcup_{t' \in \strees} \{c \makebox{ coupe admissible gauche de $t'$ telle que $P^{c}(t')=F$ et $R^c(t')=t$}\} \subset C_{F,t}.$$
On dira que $(s, \alpha) \in Angles(t)$ est un angle gauche si il est le plus petit parmi les angles ayant le m\^ eme sommet $s$, et on notera $Angles_G(t)$ l'ensemble des angles gauches de $t$. Alors l'image de $C_{F,t}^G$ par $f_1$ est l'ensemble des greffes $((s_1,\alpha_1), \ldots, (s_m,\alpha_m))$ telles que pour tout $i$, $(s_i, \alpha_i)$ ne soit pas un angle gauche. D'o\`u :
\begin{eqnarray*}
\sum_{t' \in \strees} n_G(F,t;t') &=& card(C_{F,t}^G)\\
&=&\binom{m+a-card(angles_G(t))-1}{m}.
\end{eqnarray*}
Or il y a exactement $n$ angles gauches parmi les angles de $t$ (un par sommet), ce qui donne le deuxi\`eme r\'esultat. $\Box$

\begin{cor}
Soit $F=t_1 \ldots t_m \in \sforets$, $G \in \sforets$ de poids $n$.
\begin{eqnarray*}
\sum_{H \in \sforets} n(F,G;H) &=& \binom{2n+m}{m}\\
\sum_{H \in \sforets} n_G(F,G;H) &=& \binom{n+m}{m}.
\end{eqnarray*}
\end{cor}
{\it Preuve :}
par la bijection $\alpha$ de la preuve de la proposition \ref{pro7}, pour tout $F,G,H \in \sforets$, $n(F,G;H)=n(F,B^+(G);B^+(H))$. Donc :
\begin{eqnarray*}
\sum_{H \in \sforets} n(F,G;H) &=& \sum_{H \in \sforets} n(F,B^+(G),B^+(H))\\
&=& \sum_{t' \in \strees} n(F,B^+(G),t')\\
&=& \binom{2(n+1)+m-2}{m}
\end{eqnarray*}
car $B^+(G)$ est un arbre de poids $n+1$.\\

Posons $H=s_1 \ldots s_k$, et soit $c \in {\cal A}d^G(H)$, telle que $P^c(H)=F$ et $R^c(H)=G$. Consid\'erons $\alpha(c)$. Si $s_1\neq t_1$, $c_{\mid s_1}$ n'est pas la coupe totale de $s_1$, et donc $\alpha(c)$ est une coupe admissible gauche de $B^+(H)$ telle que $P^{\alpha(c)}(B^+(H))=F$ et $R^{\alpha(c)}(B^+(H))=B^+(G)$ ; de plus, toute coupe de cette forme est obtenue ; comme $\alpha$ est injective :  
$$n_G(t_1 \ldots t_m, G;s_1 \ldots s_k)= n_G(t_1\ldots t_m,B^+(G);B^+(s_1 \ldots s_k)), \mbox{ si }t_1 \neq s_1.$$

Si $s_1=t_1$, on obtient toujours toutes les coupes admissibles gauches de $B^+(H)$ telles que $P^{c}(B^+(H))=F$ et $R^{c}(B^+(H))=B^+(G)$ ; on obtient \'egalement les coupes admissibles de $B^+(s_1 \ldots s_k)$ portant sur l'ar\^ete menant \`a $s_1$ et v\'erifiant cette condition : $c_{\mid s_2 \ldots s_k}$  est admissible gauche, avec $P^{c}(s_2 \ldots s_k)=t_2 \ldots t_m$, $R^{c}(s_2 \ldots s_k)=G$. On a donc :
\begin{eqnarray*}
n_G(t_1 \ldots t_m, G;s_1 \ldots s_k)&=& n_G(t_1\ldots t_m,B^+(G);B^+(s_1 \ldots s_k))\\
&&+n_G(t_2 \ldots t_m,G;s_2 \ldots s_k), \mbox{ si }t_1 = s_1.
\end{eqnarray*}
Donc :
\begin{eqnarray*}
\sum_{H \in \sforets} n_G(F,G;H)&=&\sum_{H \in \sforets} n_G(F,B^+(G);B^+(H))+\sum_{H \in \sforets} n_G(t_2\ldots t_m,G;H)\\
&=& \binom{n+m-1}{m} +\sum_{H\in \sforets} n_G(t_2\ldots t_m,G;H).
\end{eqnarray*}
Terminons par une r\'ecurrence sur $m$ :
si $m=1$, on a :
\begin{eqnarray*}
\sum_{H \in \sforets} n_G(F,G;H)&=&\binom{n+1-1}{1} + \sum_{H \in \sforets} n_G(1,G;H)\\
&=&n+1\\
&=&\binom{n+m}{n}.
\end{eqnarray*}
Supposons la formule vraie au rang $m-1$ :
\begin{eqnarray*}
\sum_{H \in \sforets} n_G(F,G;H)&=&\binom{n+m-1}{m} + \binom{n+m-1}{m-1}\\
&=& \binom{n+m}{m}. \esp \Box
\end{eqnarray*}

\subsection{Construction des sous-alg\`ebres}

\spa Dans $\sHr$ ou $\Hfr$, on consid\`ere les \'el\'ements suivants pour tout $n\geq 1$ :
\begin{eqnarray*}
u_n &=& \sum_{\stackrel{ F \in \sforets}{poids(F)=n}}F,\\
v_n&=& \sum_{\stackrel{t \in \strees}{poids(t)=n}}t.
\end{eqnarray*}
On a facilement :
$$u_n =\sum_{l>0}\esp \sum _{\stackrel{a_1+ \ldots a_l=n}{a_i>0}} v_{a_1} \ldots v_{a_l},$$
et donc les $u_n$ et les $v_n$ engendrent (librement) la m\^eme sous-alg\`ebre de $\sHr$ ou de $\Hfr$. On la note $\cal H$.

\begin{theo}
\label{theo61}
 Pour tout $n \geq 1$, on a :
\begin{eqnarray*}
\tdelta(v_n)&=&\sum_{k=1}^{n-1} \sum_{l>0} \binom{2n-2k+l-2}{l} \left(\sum_{\stackrel{a_1+\ldots +a_l=k}{a_i>0}} v_{a_1} \ldots v_{a_l} \right) \otimes v_{n-k},\\
\tdelta_{Fr}(v_n)&=&\sum_{k=1}^{n-1} \sum_{l>0} \binom{n-k+l-2}{l} \left(\sum_{\stackrel{a_1+\ldots +a_l=k}{a_i>0}} v_{a_1} \ldots v_{a_l} \right) \otimes v_{n-k}\esp ;\\
\tdelta(u_n)&=&\sum_{k=1}^{n-1} \sum_{l>0} \binom{2n-2k+l}{l} \left(\sum_{\stackrel{a_1+\ldots +a_l=k}{a_i>0}} v_{a_1} \ldots v_{a_l} \right) \otimes u_{n-k},\\
\tdelta_{Fr}(u_n)&=&\sum_{k=1}^{n-1} \sum_{l>0} \binom{n-k+l}{l} \left(\sum_{\stackrel{a_1+\ldots +a_l=k}{a_i>0}} v_{a_1} \ldots v_{a_l} \right) \otimes u_{n-k}.
\end{eqnarray*}
\end{theo}
{\it Preuve : } $v_n$ est une combinaison lin\'eaire d'arbres de poids $n$ ; par d\'efinition de $\Delta$ :
\begin{eqnarray*}
\tdelta(v_n) &=& \sum_{poids(t')=n} \esp \sum_{F \in \sforets }\esp \sum_{t \in \strees} n(F,t;t') F \otimes t\\
&=& \sum_{poids(t')=n} \esp \sum_{k=1}^{n-1} \esp \sum_{l>0} \esp \sum_{poids(t_1\ldots t_l)=k}\esp \sum_{poids(t)=n-k} n(t_1 \ldots t_l,t;t') F \otimes t\\
&=& \sum_{k=1}^{n-1} \sum_{l>0} \esp \sum_{poids(t_1 \ldots t_l)=k}\esp \sum_{poids(t)=n-k} \left(\sum_{poids(t')=n} n(t_1 \ldots t_l,t;t')\right) t_1 \ldots t_l \otimes t\\
&=& \sum_{k=1}^{n-1} \sum_{l>0} \esp \sum_{poids(t_1 \ldots t_l)=k}\esp \sum_{poids(t)=n-k} \binom{2n-2k+l-2}{l} t_1 \ldots t_l\otimes t\\
\end{eqnarray*}
\begin{eqnarray*}
&=& \sum_{k=1}^{n-1} \sum_{l>0} \esp \esp  \binom{2n-2k+l-2}{l}  \left(\sum_{poids(t_1 \ldots t_l)=k}t_1 \ldots t_l\right)\otimes \left(\sum_{poids(t)=n-k} t\right)\\
&=& \sum_{k=1}^{n-1} \sum_{l>0} \esp \esp  \binom{2n-2k+l-2}{l}  \left(\sum_{\stackrel{a_1+\ldots +a_l=k}{a_i>0}} v_{a_1} \ldots v_{a_l} \right)\otimes v_{n-k}.
\end{eqnarray*}
(On a utilis\'e le corollaire  \ref{cor60} ainsi que le fait que $n(F,t;t')=0$ si $poids(t') \neq poids(F)+poids(t)$ pour la quatri\`eme \'egalit\'e.)

 Les trois autres calculs sont identiques. $\Box$

\begin{cor}
$\cal H$ est une sous alg\`ebre de Hopf de $\Hfr$ et de $\sHr$ ; de plus les abelianis\'ees de $({\cal H}, \Delta)$ et de  $({\cal H}, \Delta_{Fr})$ sont isomorphes \`a l'alg\`ebre de Hopf de Connes-Moscovici comme alg\`ebres gradu\'ees.
\end{cor}
{\it Preuve :} consid\'erons le cas de $({\cal H}, \Delta_{Fr})$. Il suffit de montrer que $ {\cal H}_{ab}^{*g}$ est isomorphe \`a ${\cal H}_{CM}^{*g}$ ; comme il s'agit de deux alg\`ebres gradu\'ees cocommutatives (et donc d'alg\`ebres enveloppantes), il s'agit donc de montrer que leurs alg\`ebres de Lie des \'el\'ements primitifs sont des alg\`ebres de Lie isomorphes. Dans le cas de ${\cal H}_{ab}^{*g}$, d'apr\`es la proposition \ref{prop6}, une base de l'alg\`ebre de Lie des \'el\'ements primitifs $\g_{Fr}$ est donn\'ee par $(L_i)_{i\in \mathbb{N}^*}$
d\'efinie par :
$$L_i(v_{i_1} \ldots v_{i_n})= \delta_{v_i, v_{i_1} \ldots v_{i_n}}.$$
De plus, $L_i$ est homog\`ene de poids $i$.
On a alors, par homog\'en\'eit\'e, $[L_i, L_j]_{Fr}= \alpha_{i,j} L_{i+j}$, $\alpha_{i,j} \in \mathbb{Q}$. Par dualit\'e :
\begin{eqnarray*}
\alpha_{i,j}&=& ([L_i,L_j]_{Fr}, v_{i+j})\\
&=&(L_i \otimes L_j - L_j \otimes L_i, \Delta_{Fr}(v_{i+j}))\\
&=&(L_i \otimes L_j, \sum_{l>0} \binom{j+l-2}{l} \left(\sum_{\stackrel{a_1+\ldots +a_l=i}{a_k>0}} v_{a_1} \ldots v_{a_l} \right) \otimes v_j)\\
&&-(L_j \otimes L_i, \sum_{l>0} \binom{i+l-2}{l} \left(\sum_{\stackrel{a_1+\ldots +a_l=j}{a_k>0}} v_{a_1} \ldots v_{a_l} \right) \otimes v_i)+0\\
&=&(L_i \otimes L_j, \binom{j+1-2}{1} v_i \otimes v_j) -(L_j \otimes L_i, \binom{i+1-2}{1} v_j \otimes v_i) +0\\
&=& j-i.
\end{eqnarray*}
(On a utilis\'e l'expression de $\Delta_{Fr}(v_{i+j})$ ainsi que l'homog\'en\'eit\'e de $L_i$ et $L_j$ pour la troisi\`eme \'egalit\'e).\\
Donc $Z_i \longmapsto L_i$ est un isomorphisme de $\g_{CM}$ sur $\g_{Fr}$ d'apr\`es (\ref{eq25}).

Dans le cas de $({\cal H},\Delta)$,   on note $\g $ l'alg\`ebre de Lie des primitifs de ${\cal H}_{ab}^{*g}$. $\g$ a la m\^eme base que $\g_{Fr}$, et un calcul semblable au pr\'ec\'edent montre que  :
$$[L_i, L_j]= 2(j-i) L_{i+j}.$$
Donc $Z_i \longmapsto \frac{1}{2} L_i$ est un isomorphisme de $\g_{CM}$ sur $\g$. $\Box$\\

{\it Remarque :} on retrouve ainsi les r\'esultats de \cite{Frabetti}.\\
  
\subsection{Isomorphisme entre $({\cal H}, \Delta_{Fr})$ et $({\cal H}, \Delta)$}

\spa Soit $V=\bigoplus V_n$ un espace gradu\'e. On d\'efinit $val(v)$ pour tout $v \in V$ par :
$$val(v)=\max\{n/v \in V_n \oplus V_{n+1}\oplus \ldots \}.$$
$V$ est alors muni d'une distance donn\'ee par $d(v,v')=2^{-val(v-v')}$.
Un compl\'et\'e de $V$ pour cette distance est donn\'e par :
$$\overline{V}=\prod_{n=0}^{+\infty} V_n.$$
Les \'el\'ements de $V$ seront not\'es $\sum v_n$, $v_n \in V_n$ pour tout $n \in \mathbb{N}$.
Soit $V$, $V'$ deux espaces gradu\'es, et soit $f:V \longmapsto V'$ homog\`ene de degr\'e $i$. On v\'erifie alors que $f$ est lipschitzienne de rapport $2^{-i}$, et donc se prolonge de mani\`ere unique en un application $\overline{f}:\overline{V} \longmapsto \overline{V}'$.

En particulier, si $(A,m)$ est une alg\`ebre gradu\'ee, $m$ se prolonge en $\overline{m} : \overline{A\otimes A} \longmapsto \overline{A}$. On a une injection naturelle :
\begin{eqnarray*}
\overline{A}\otimes \overline{A} &\longmapsto &\overline{A\otimes A}\\
\left(\sum_i a_i \right)\otimes\left(\sum_j b_j\right)& \longmapsto & \sum_{n} \left(\sum_{i+j=n}a_i \otimes b_j\right).
\end{eqnarray*}
On peut donc consid\'erer $\overline{m}: \overline{A}\otimes \overline{A}\longmapsto \overline{A}$ ; un simple raisonnement par densit\'e montre que $(\overline{A}, \overline{m})$ est  une alg\`ebre. Son produit est donn\'e par :
$$\left(\sum_i a_i \right)\left(\sum_j b_j\right)=\sum_{n} \left(\sum_{i+j=n}a_i  b_j\right).$$

Si $(A,m,\Delta)$ est une big\`ebre gradu\'ee, on peut prolonger $\Delta$ en $\overline{\Delta}:\overline{A} \longmapsto \overline{A\otimes A}$. Un raisonnement par densit\'e montre que :
$$\overline{\Delta}(ab)=\overline{\Delta}(a) \overline{\Delta}(b), \esp \forall a,b \in \overline{A}.$$ 

Appliquons ceci \`a l'alg\`ebre $\cal H$. On consid\`ere :
$$ U=1+\sum_{n=1}^{+\infty} u_n \in \overline{\cal H}, \esp V=\sum_{n=1}^{+\infty} v_n \in \overline{\cal H}.$$

\begin{prop}
\begin{eqnarray*}
\overline{\Delta}_{Fr}(U)&=&\sum_{j=0}^{+\infty} U^{j+1} \otimes u_j,\\
\overline{\Delta}(U)&=&\sum_{j=0}^{+\infty} U^{2j+1} \otimes u_j.
\end{eqnarray*}
\end{prop}
{\it Preuve :}
\begin{eqnarray}
\nonumber U&=&1+\sum_{n=1}^{+\infty} \sum_{l>0} \esp \sum_{\stackrel{a_1+\ldots +a_l=n}{a_i>0}} v_{a_1} \ldots v_{a_l}\\
\nonumber &=&1+\sum_{l>0} \sum_{a_i>0} v_{a_1} \ldots v_{a_l}\\
\nonumber &=&1+ \sum_{l>0} V^l\\
\label{eqnA} &=&\frac{1}{1-V}.
\end{eqnarray}
De plus :
\begin{eqnarray*}
{\Delta}_{Fr}(u_n)&=&  \sum_{k=1}^{n-1} \sum_{l>0} \binom{n-k+l}{l} \left(\sum_{\stackrel{a_1+\ldots +a_l=k}{a_i>0}} v_{a_1} \ldots v_{a_l} \right) \otimes u_{n-k}\\
&&+1\otimes u_n +u_n\otimes 1\\
&=&  \sum_{k=1}^{n} \sum_{l>0} \binom{n-k+l}{l} \left(\sum_{\stackrel{a_1+\ldots +a_l=k}{a_i>0}} v_{a_1} \ldots v_{a_l} \right) \otimes u_{n-k}\\
&&+1\otimes u_n\\
&=&  \sum_{j=0}^{n-1} \sum_{l>0} \binom{j+l}{l} \left(\sum_{\stackrel{a_1+\ldots +a_l=n-j}{a_i>0}} v_{a_1} \ldots v_{a_l} \right) \otimes u_{j}\\
&&+1\otimes u_n,
\end{eqnarray*}
cette formule restant vraie pour $n=0$. Par suite :
\begin{eqnarray*}
\overline{\Delta}_{Fr}(U)&=&\sum_{0}^{+\infty} \Delta_{Fr}(u_n)\\
&=&\sum_{j=0}^{+\infty} \sum_{l>0}  \binom{j+l}{j} \left(\sum_{a_1, \ldots, a_l>0} v_{a_1}\ldots v_{a_l}\right)\otimes u_j +1 \otimes U\\
&=&\sum_{j=0}^{+\infty} \left(\sum_{l>0} \binom{j+l}{j} V^l \right)\otimes u_j +1 \otimes U.
\end{eqnarray*}
Or dans $\mathbb{Q}[[X]]$, on a :
$$\sum_{l=0}^{+\infty} \binom{j+l}{j} X^l=\frac{1}{(1-X)^{j+1}}.$$
Donc :
\begin{eqnarray*}
\overline{\Delta}_{Fr}(U)&=&\sum_{j=0}^{+\infty} \frac{1}{(1-V)^{j+1}} \otimes u_j -\sum_{j=0}^{+\infty}1\otimes u_j  +1 \otimes U\\
&=&\sum_{j=0}^{+\infty} {U^{j+1}} \otimes u_j .
\end{eqnarray*}
(On  a utilis\'e (\ref{eqnA}) pour la deuxi\`eme \'egalit\'e).\\
Le calcul de $\overline{\Delta}(U)$ est similaire. $\Box$

\begin{theo}
On consid\`ere les \'el\'ements de $\cal H$ d\'efinis de la mani\`ere suivante :
\begin{eqnarray*}
z_n&=&2u_n +\sum_{k=1}^{n-1} u_k u_{n-k}\esp ;\\ 
w_n&=&\frac{1}{2} u_n -\frac{1}{2} \sum_{i=1}^{n-1} w_i  w_{n-i}.
\end{eqnarray*}
Soit $\Phi :({\cal H}, \Delta_{Fr}) \longmapsto ({\cal H}, \Delta)$ l'unique morphisme d'alg\`ebres envoyant $u_n$ sur $z_n$ pour tout $n \in \mathbb{N}^*$. Alors $\Phi$ est un isomorphisme d'alg\`ebres de Hopf gradu\'ees ; son inverse est donn\'e par $\Phi^{-1}(u_n) = w_n$.
\end{theo}
{\it Preuve :} une r\'ecurrence simple montre que les $z_n$ engendrent librement $\cal H$ ; par suite, $\Phi$ est bijectif. $z_n$ est homog\`ene de poids $n$, et donc $\Phi$ est homog\`ene de degr\'e z\'ero. 

Soit $Z=1+\sum z_n \in \overline{\cal H}$.
\begin{eqnarray*}
Z&=&1+\sum_{n=1}^{+\infty} \left(2u_n +\sum_{k=1}^{n-1}u_k u_{n-k}\right)\\
&=&1+2U+ (U-1)(U-1)\\
&=&U^2.
\end{eqnarray*}
On a alors :
\begin{eqnarray*}
\overline{\Delta}(Z)&=&\overline{\Delta}(U)^2\\
&=&\left(\sum_{j=0}^{+\infty} U^{2j+1} \otimes u_j\right) \left(\sum_{k=0}^{+\infty} U^{2k+1} \otimes u_k\right)\\ 
&=&\sum_{n=0}^{+\infty}\sum_{j+k=n} U^{2(j+k+1)}\otimes u_ju_k\\
&=&\sum_{n=0}^{+\infty} Z^{n+1} \otimes \left(2 u_n +\sum_{j=1}^{n-1} u_j u_{n-j}\right)\\
&=& \sum_{n=0}^{+\infty} Z^{n+1} \otimes z_n.
\end{eqnarray*}
 Consid\'erons $\overline{\Phi}:(\overline{\cal H}, \overline{\Delta}_{Fr}) \longmapsto (\overline{\cal H}, \overline{\Delta})$. Le calcul pr\'ec\'edent montre que $\overline{\Delta}(\overline{\Phi}(U))= (\overline{\Phi}\otimes \overline{\Phi})\overline{\Delta}_{Fr}(U).$ Par suite, en consid\'erant chaque composante homog\`ene, on en d\'eduit que $\Phi$ est un morphisme d'alg\`ebres de Hopf.

Enfin, montrons par r\'ecurrence que $\Phi(w_n)=u_n$ : si $n=1$, alors $z_1=2 u_1$, $w_1=\frac{1}{2} u_1$. Par suite, $\Phi(w_1)=\frac{1}{2}z_1=u_1$. Supposons la propri\'et\'e vrai jusqu'au rang $n-1$ :
\begin{eqnarray*}
\Phi(w_n)&=& \frac{1}{2} \Phi(u_n)-\frac{1}{2} \sum_{i=1}^{n-1} \Phi(w_i)\Phi(w_{n-i})\\
&=& u_n +\frac{1}{2} \sum_{i=1}^{n-1} u_iu_{n-i}-\frac{1}{2} \sum_{i=1}^{n-1} u_iu_{n-i}\\
&=&u_n. \esp \Box
\end{eqnarray*}

{\it Remarque :} on montre facilement que dans $\overline{\cal H}$ :
$$ \left(1+\sum_{n=1}^{+\infty}w_n \right)^2=U.$$

\section{Cog\`ebre tensorielle d'un espace vectoriel}
\spa Dans toute cette section, $V$ d\'esigne un espace vectoriel quelconque sur un corps commutatif $K$.

\subsection{ Construction et caract\'erisation}
\label{part4}
\spa Soit $T(V)=\bigoplus_{n=0}^{\infty} V^{\otimes n}.$
Pour $v_1, \ldots v_n \in V$, on note leur produit tensoriel dans $V^{\otimes n}$ $v_1 \top \ldots \top v_n$ plut\^ot que $v_1 \otimes \ldots \otimes v_n$ ou $v_1 \ldots v_n$.
On munit cet espace d'une structure de cog\`ebre en posant :
\begin{eqnarray}
\label{eqn13}
\Delta(v_1 \top \ldots \top v_n)=\sum_{k=0}^{k=n} (v_1 \top \ldots \top v_k) \otimes (v_{k+1} \top \ldots \top v_n).
\end{eqnarray}
D'apr\`es \cite{Bou1}, chapitre III, \S 11, $(T(V),\Delta)$ est bien une cog\`ebre appel\'ee cog\`ebre tensorielle de $V$, et sa counit\'e est donn\'ee par :
$$\varepsilon(1)=1, \esp \varepsilon(v_1 \top \ldots \top v_n)=0 \mbox{ si $n\geq 1$}.$$
$(V^{\otimes n})_{n \in \mathbb{N}}$ est une graduation de la cog\`ebre $T(V)$ v\'erifiant $(C_1)$. $T(V)$ est donc filtr\'ee par $deg_p$ comme on l'a vu dans la section \ref{partB}.\\

{\it Remarque :} d'apr\`es la propri\'et\'e 7 du th\'eor\`eme \ref{theo28},  $\Hr$ est une cog\`ebre tensorielle, $\top$ \'etant donn\'e par :
$$ e_{t_1} \top \ldots \top e_{t_n} = e_{t_1\ldots t_n}.$$ 
\label{section72}

\begin{figure}[h]
\framebox(450,100){
\begin{picture}(0,0)(80,-20)
$p=-2$
\begin{picture}(8,20)(0,0)
\put(0,0){\circle*{3}}
\put(0,0){\line(0,1){10}}
\put(0,10){\circle*{3}}
\end{picture}
$+$
\begin{picture}(20,20)(-2,0)
\put(0,0){\circle*{3}}
\put(10,0){\circle*{3}}
\end{picture}
$=e$
\begin{picture}(8,20)(0,10)
\put(0,0){\circle*{3}}
\put(0,0){\line(0,1){10}}
\put(0,10){\circle*{3}}
\end{picture}
 ;
\end{picture}
\begin{picture}(0,0)(-170,-20)
$q=$
\begin{picture}(8,10)(0,0)
\put(3,0){\circle*{3}}
\end{picture}
$=e$
\begin{picture}(8,10)(0,0)
\put(0,0){\circle*{3}}
\end{picture}
 ;
\end{picture}
\\
\begin{picture}(0,0)(180,30)
$q\top p =e$
\begin{picture}(20,20)(0,8)
\put(0,0){\circle*{3}}
\put(10,0){\circle*{3}}
\put(10,10){\circle*{3}}
\put(10,0){\line(0,1){10}}
\end{picture}
$=-$
\begin{picture}(10,30)(0,0)
\put(0,0){\circle*{3}}
\put(0,10){\circle*{3}}
\put(0,20){\circle*{3}}
\put(0,0){\line(0,1){10}}
\put(0,10){\line(0,1){10}}
\end{picture}
$-$
\begin{picture}(25,20)(0,0)
\put(10,0){\circle*{3}}
\put(0,10){\circle*{3}}
\put(20,10){\circle*{3}}
\put(10,0){\line(1,1){10}}
\put(10,0){\line(-1,1){10}}
\end{picture}
$+$
\begin{picture}(20,20)(0,0)
\put(0,0){\circle*{3}}
\put(10,0){\circle*{3}}
\put(10,10){\circle*{3}}
\put(10,0){\line(0,1){10}}
\end{picture}
.
\end{picture}
}
\caption{{\it l'application bilin\'eaire $\top$ dans ${\cal H}_{P,R}$.}}
\end{figure}

On rappelle que $\tilde{\Delta}^k$ est d\'efini dans la section \ref{partB}.
\begin{lemme}
\label{lemme13}
Soient $v_1, \ldots, v_n \in V$. 
\begin{eqnarray*}
\tilde{\Delta}^{n-1}(v_1 \top \ldots \top v_n)&=&v_1 \otimes \ldots \otimes v_n \mbox{ si }n\geq 2 \esp ;\\
\tilde{\Delta}^{k}(v_1 \top \ldots \top v_n)&=&0 \mbox{ si }k\geq n.
\end{eqnarray*}
\end{lemme}
{\it Preuve :} r\'ecurrence facile sur $n$. $\Box$

\begin{prop}
\label{prop14}
\begin{enumerate}
\item $1$ est le seul \'el\'ement non nul de $T(V)$ tel que $\Delta(e)=e \otimes e$.
\item Soit $v\in V$. On d\'efinit $L_v :T(V) \longmapsto T(V)$ par $L_v(v_1 \top \ldots \top v_k)= v_1 \top \ldots \top v_k \top v$. Alors
$ \forall x \in T(V), \esp \Delta(L_v(x))=L_v(x) \otimes 1 + (Id \otimes L_v)\circ \Delta(x).$
\item $Prim(T(V))=V$.
\item $Ker(\tilde{\Delta}^n)=V \oplus \ldots \oplus V^{\otimes n}.$
\end{enumerate}
\end{prop}
{\it Preuve : }

1. Voir \cite{Bou1}.

2. D\'ecoule imm\'ediatement de (\ref{eqn13}).

3. Soit $x$ un primitif de $T(V)$. Alors $\varepsilon(x)=0$, donc $x \in V \oplus \ldots \oplus V^{\otimes n}$ pour un certain $n \geq 1$. Posons $x=x_1 + \ldots +x_n$, avec $x_i \in V^{\otimes i}$. Supposons $n \geq 2$, et $x_n \neq 0$.
D'apr\`es le lemme \ref{lemme13}, $\tilde{\Delta}^{n-1}(x)=\tilde{\Delta}^{n-1}(x_n)=0$, car $\tilde{\Delta}(x)=0$.
Or, toujours d'apr\`es le lemme \ref{lemme13}, $\tilde{\Delta}^{n-1} :V^{\otimes n}\longmapsto V^{\otimes n}$ est l'identit\'e :
par suite, $x_n=0$ : on aboutit \`a une contradiction. Donc $n=1$, et $Prim(T(V)) \subseteq V$. La r\'eciproque est triviale.

4. On vient de le montrer pour $n=1$. Soit $n>1$ ; supposons le r\'esultat acquis pour tout $k<n$, et soit $x\in T(V)$, $\tilde{\Delta}^n(x)=0$.
Posons $\tilde{\Delta}^{n-1}(x)=\sum x^{(1)} \otimes \ldots \otimes x^{(n)}$. D'apr\`es le lemme \ref{lemme11}, on peut supposer que les $x^{(i)}$ sont primitifs, donc sont des \'el\'ements de $V$.
Alors  $\tilde{\Delta}^{n-1}(x-\sum x^{(1)} \top \ldots \top x^{(n)})=0$, donc $x \in Ker(\tilde{\Delta}^{n-1})+ V^{\otimes n}$,
d'o\`u $ Ker(\tilde{\Delta}^{n}) \subseteq V \oplus \ldots \oplus V^{\otimes n}$ d'apr\`es l'hypoth\`ese de r\'ecurrence.
L'inclusion r\'eciproque d\'ecoule imm\'ediatement du lemme \ref{lemme13}. $\Box$\\

On constate alors que $T(V)_{deg_p\leq n}=K \oplus V \oplus \ldots \oplus V^{\otimes n}$. La filtration par $deg_p$ provient donc de la graduation dont les composantes  homog\`enes sont donn\'ees par 
$T(V)_{deg_p=n}=V^{\otimes n}$.

\begin{theo}
\label{theo13}
Soit $(C, \Delta_C, \varepsilon)$ une cog\`ebre v\'erifiant :
\begin{enumerate}
\item $\exists \esp e \in C-\{0\}, \esp \Delta_C(e)=e \otimes e$.
\item Soit $Prim(C)=\{x \in C/\Delta_C(x) =x \otimes e + e \otimes x\}$ ; alors
    $\exists \esp L :Prim(C) \longmapsto {\cal L}(C)$, v\'erifiant :
  \begin{description}
    \item[\textnormal{a)}] $Lp(e)=p, \esp \forall p\in Prim(C),$
    \item[\textnormal{b)}] $\Delta_C\left(L_p(x)\right)=L_p(x) \otimes e + (Id \otimes L_p)\circ \Delta_C(x), \esp \forall x \in C$.
 \end{description}
\item On pose $\tilde{\Delta}_C(x)=\tilde{\Delta}_C^1(x)=\Delta_C(x)-x \otimes e - e \otimes x$, et par r\'ecurrence on d\'efinit
            $\tilde{\Delta}_C^n=\left( \tilde{\Delta}_C^{n-1} \otimes Id\right) \circ \tilde{\Delta}_C^1$ ; alors 
pour tout $x\in Ker(\varepsilon)$, il existe $n\geq 1$, $\tilde{\Delta}_C^n(x)=0$.
\end{enumerate}
Alors $C$ et $T(Prim(C))$ sont isomorphes.
\end{theo}
{\it Preuve : }soient $p_1, \ldots ,p_n$ des \'el\'ements primitifs de $C$.\\
On d\'efinit par r\'ecurrence $p_1 \ttop \ldots \ttop p_n=L_{p_n}(p_1\ttop \ldots \ttop p_{n-1})$ (on utilise la convention $p_1 \ttop \ldots \ttop p_n=e$ si $n=0$).
Montrons par r\'ecurrence sur $n$ que : 
\begin{eqnarray}
\label{eqn14}
\Delta_C(p_1 \ttop \ldots \ttop p_n)&=&\sum_{k=0}^{k=n}(p_1 \ttop \ldots \ttop p_k) \otimes (p_{k+1} \ttop \ldots \ttop p_n).
\end{eqnarray}
C'est vrai pour $n=0$. Supposons le r\'esultat vrai au rang $n-1$ ; alors
\begin{eqnarray*}
\Delta_C(p_1 \ttop \ldots \ttop p_n)&=&\Delta_C\left(L_{p_n}(p_1 \ttop \ldots \ttop p_{n-1})\right)\\
                                  &=&(p_1 \ttop \ldots \ttop p_n)  \otimes e + (Id \otimes L_{p_n})\circ\Delta_C(p_1 \ttop \ldots \ttop p_{n-1})\\
                                  &=& (p_1 \ttop \ldots \ttop p_n) \otimes e + \sum_{k=0}^{k=n-1} (p_1 \ttop \ldots \ttop p_k) \otimes (p_{k+1} \ttop \ldots \ttop p_n)\\
                                  &=& \sum_{k=0}^{k=n}(p_1 \ttop \ldots \ttop p_k) \otimes (p_{k+1} \ttop \ldots \ttop p_n).
\end{eqnarray*}

Comme $L$ est lin\'eaire, $(p_1, \ldots, p_n) \longmapsto p_1 \ttop \ldots \ttop p_n$ est $n$-lin\'eaire. On peut donc d\'efinir :
\begin{eqnarray}
\nonumber
F :T(Prim(C))& \longmapsto & C\\
  p_1 \top \ldots \top p_n & \longmapsto & p_1 \ttop \ldots \ttop p_n.
\end{eqnarray}
D'apr\`es (\ref{eqn13}) et (\ref{eqn14}), $F$ est un morphisme de cog\`ebres.\\

Montrons que $F$ est injective : soit $x=x_0 1+x_1 + \ldots +x_n$, avec $x_0 \in K$, $x_i \in Prim(C)^{\otimes i}$, tel que $F(x)=0$. Supposons $x_n\neq 0$.
Comme $\varepsilon(F(x))=\varepsilon(x)=x_0$, on a $x_0=0$ : $n\geq 1$. De plus, $F(p)=p$, $\forall p\in Prim(C)$, donc $n \geq 2$.
Comme $F$ est un morphisme de cog\`ebres et que $F(1)=e$, on a $\tilde{\Delta}_C^k\circ F=F^{\otimes (k+1)}\circ \tilde{\Delta}^k$ pour tout $k\geq 1 $.
Par suite, dans $Prim(C)^{\otimes n}$ :
\begin{eqnarray*}
\tilde{\Delta}_C^{n-1}(F(x))&=&0\\
&=& F^{\otimes n} \circ \tilde{\Delta}^{n-1}(x)\\
&=& F^{\otimes n} \circ \tilde{\Delta}^{n-1}(x_n)\\
&=& \tilde{\Delta}^{n-1}(x_n)\\
&=& x_n.
\end{eqnarray*}
(On a utilis\'e le lemme \ref{lemme13} pour la troisi\`eme et la derni\`ere \'egalit\'e et le fait que $F(p)=p$ $\forall p\in Prim(C)$ pour la quatri\`eme).\\
On aboutit \`a une contradiction, donc  $F$ est injective. \\

Il reste \`a montrer que $F$ est surjective.

 Montrons que $Ker(\tilde{\Delta}_C^{n})\subseteq F(Prim(C) \oplus \ldots \oplus Prim(C)^{\otimes n})$.
En utilisant l'hypoth\`ese 3 du th\'eor\`eme, on pourra conclure.
Proc\'edons par r\'ecurrence sur $n$ : c'est vrai pour $n=1$, $Ker(\tilde{\Delta}_C^{1})=Prim(C)\subseteq F(Prim(C))$. Supposons l'hypoth\`ese vraie au rang $n-1$.
Soit $x \in C$ tel que $\tilde{\Delta}_C^{n}(x)=0$, posons $\tilde{\Delta}_C^{n-1}(x)=\sum x^{(1)}\otimes \ldots \otimes x^{(n)}$.
D'apr\`es le lemme \ref{lemme11}, on peut supposer les $x^{(i)}$ primitifs.
Alors $\tilde{\Delta}_C^{n-1}(x-F(\sum x^{(1)} \top \ldots \top x^{(n)}))=0$ ; 
donc $x \in Ker(\tilde{\Delta}_C^{n-1})+F(Prim(C)^{\otimes n})$ ; par l'hypoth\`ese de r\'ecurrence,
$Ker(\tilde{\Delta}_C^{n})\subseteq F(Prim(C) \oplus \ldots \oplus Prim(C)^{\otimes n})$. $\Box$

\subsection{Cas des alg\`ebres de Hopf $\nH$}
\spa On note ${\cal H}_R$ l'alg\`ebre de Hopf des arbres enracin\'es (non plans) de \cite{Connes}.
On note $\nH$ l'alg\`ebre de Hopf des arbres enracin\'es d\'ecor\'es par $\cal D$, o\`u $\cal D$ est un ensemble fini non vide (voir  \cite{Connes,Kreimer2}).
On note  $\ntrees$ l'ensemble des arbres enracin\'es d\'ecor\'es par $\cal D$, et $\nforets$ l'ensemble des mon\^omes en les \'el\'ements de $\ntrees$ de $\nH$.\\

Soient $F,G \in \nforets$. On pose :
\begin{eqnarray*}
F \otop G&=& \frac{1}{poids(G)} \sum_{s \in som(G)}  \mbox{(greffe de $F$ sur le sommet $s$), si $G \neq 1$ ;}\\
&=&0 \mbox{, si $G=1$.}
\end{eqnarray*}
On prolonge $\otop$ en une application bilin\'eaire de $\nH \times \nH$ dans $\nH$. On consid\`ere :
\begin{eqnarray*}
L \esp : Prim(\nH)&\longmapsto & {\cal L}(\nH)\\
p&\longmapsto& L_p \esp : 
\left\{ \begin{array}{rcl} \nH &\longmapsto &\nH\\
x &\longmapsto& x \otop p
\end{array} \right.
\end{eqnarray*}

\begin{figure}[h]
\begin{picture}(15,30)(-10,0)
\put(0,0){\circle*{3}}
\put(0,0){\line(0,1){10}}
\put(0,10){\circle*{3}}
\end{picture}
$\otop$
\begin{picture}(30,30)(-10,0)
\put(0,0){\circle*{3}}
\put(0,0){\line(1,1){10}}
\put(0,0){\line(-1,1){10}}
\put(-10,10){\circle*{3}}
\put(10,10){\circle*{3}}
\end{picture}
= $\frac{1}{3} \left(
\begin{picture}(30,30)(-13,15)
\put(0,0){\circle*{3}}
\put(0,0){\line(1,1){10}}
\put(0,0){\line(-1,1){10}}
\put(-10,10){\circle*{3}}
\put(10,10){\circle*{3}}
\put(-10,10){\line(0,1){10}}
\put(-10,20){\line(0,1){10}}
\put(-10,20){\circle*{3}}
\put(-10,30){\circle*{3}}
\end{picture}
+
\begin{picture}(30,30)(-13,15)
\put(0,0){\circle*{3}}
\put(0,0){\line(1,1){10}}
\put(0,0){\line(-1,1){10}}
\put(-10,10){\circle*{3}}
\put(10,10){\circle*{3}}
\put(10,10){\line(0,1){10}}
\put(10,20){\line(0,1){10}}
\put(10,20){\circle*{3}}
\put(10,30){\circle*{3}}
\end{picture}
+
\begin{picture}(30,30)(-13,15)
\put(0,0){\circle*{3}}
\put(0,0){\line(1,1){10}}
\put(0,0){\line(-1,1){10}}
\put(-10,10){\circle*{3}}
\put(10,10){\circle*{3}}
\put(0,0){\line(0,1){10}}
\put(-10,10){\line(0,1){10}}
\put(0,10){\circle*{3}}
\put(-10,20){\circle*{3}}
\end{picture}
\right)$
 ; 
\begin{picture}(20,10)
\end{picture}
\begin{picture}(20,30)(0,0)
\put(0,0){\circle*{3}}
\put(10,0){\circle*{3}}
\end{picture}
$\otop$
\begin{picture}(10,30)(-10,10)
\put(0,0){\circle*{3}}
\put(0,0){\line(0,1){10}}
\put(0,10){\circle*{3}}
\put(0,10){\line(0,1){10}}
\put(0,20){\circle*{3}}
\end{picture}
$=\frac{1}{3}
\left(
\begin{picture}(25,30)(-10,10)
\put(0,0){\circle*{3}}
\put(0,0){\line(0,1){10}}
\put(0,10){\circle*{3}}
\put(10,10){\line(0,1){10}}
\put(10,20){\circle*{3}}
\put(10,10){\circle*{3}}
\put(-10,10){\circle*{3}}
\put(0,0){\line(1,1){10}}
\put(0,0){\line(-1,1){10}}
\end{picture}
+\begin{picture}(20,30)(-10,10)
\put(0,0){\circle*{3}}
\put(0,0){\line(0,1){10}}
\put(0,10){\circle*{3}}
\put(0,10){\line(0,1){10}}
\put(0,20){\circle*{3}}
\put(10,20){\circle*{3}}
\put(-10,20){\circle*{3}}
\put(0,10){\line(1,1){10}}
\put(0,10){\line(-1,1){10}}
\end{picture}
+\begin{picture}(20,40)(-10,10)
\put(0,0){\circle*{3}}
\put(0,0){\line(0,1){10}}
\put(0,10){\circle*{3}}
\put(0,10){\line(0,1){10}}
\put(0,20){\circle*{3}}
\put(10,30){\circle*{3}}
\put(-10,30){\circle*{3}}
\put(0,20){\line(1,1){10}}
\put(0,20){\line(-1,1){10}}
\end{picture}
\right)$
 ;

\begin{picture}(130,0)
\end{picture}
\begin{picture}(10,30)(0,0)
\put(0,0){\circle*{3}}
\put(0,0){\line(0,1){10}}
\put(0,10){\circle*{3}}
\end{picture}
$\otop$
\begin{picture}(20,30)(-5,0)
\put(0,0){\circle*{3}}
\put(10,0){\circle*{3}}
\end{picture}
$= \frac{1}{2} \left(
\begin{picture}(20,30)(-5,10)
\put(0,0){\circle*{3}}
\put(10,0){\circle*{3}}
\put(0,10){\circle*{3}}
\put(0,10){\line(0,1){10}}
\put(0,20){\circle*{3}}
\put(0,0){\line(0,1){10}}
\end{picture}
+
\begin{picture}(20,30)(-5,10)
\put(0,0){\circle*{3}}
\put(10,0){\circle*{3}}
\put(10,10){\circle*{3}}
\put(10,10){\line(0,1){10}}
\put(10,20){\circle*{3}}
\put(10,0){\line(0,1){10}}
\end{picture}
\right)$

\begin{picture}(300,-200)(0,5)
\put(0,0){\line(0,1){135}}
\put(0,0){\line(1,0){460}}
\put(460,0){\line(0,1){135}}
\put(0,135){\line(1,0){460}}
\end{picture}

\caption{\it L'application bilin\'eaire  $\otop$.}

\end{figure}

\begin{prop}
  $L$ v\'erifie la condition 2 du th\'eor\`eme \ref{theo13}, avec $e=1$.
\end{prop}
{\it Preuve :} 
soient $F,G \in \nforets$. Soit $s$ un sommet de $G$, et soit $H$ la for\^et obtenue en greffant $F$ sur le sommet $s$.
Soit $c$ une  coupe admissible non vide et non totale de $H$.
\begin{enumerate}
\item $c$ coupe les ar\^etes reliant les racines de $F$ \`a $s$. Alors :
\begin{enumerate}
\item soit $c_{\mid G}$ est vide, et alors $P^c(H)=F$, $R^c(H)=G$.
\item soit $c' =c_{\mid G}$ est admissible, non vide, et alors $P^c(H)=FP^{c'}(G)$, $R^c(H)=R^{c'}(G)$. De plus, comme $c$ est admissible, $s$ est n\'ecessairement l'un des sommets de $R^{c'}(G)$.
\end{enumerate}
\item $c$ coupe au moins une ar\^ete de $F$ ou une ar\^ ete de $s$ vers une racine de $F$, et ne coupe pas toutes les ar\^etes de $s$ vers une racine de $F$. Alors $c'=c_{\mid G}$ et $c''=c_{\mid F}$ sont admissibles, et $c''$ est non vide.
\begin{enumerate}
\item Si $c'$ est vide, alors $P^c(H)=P^{c''}(F)$ ; de plus   $R^{c}(H)$ est la greffe de $R^{c''}(F)$ sur le sommet $s$ de $G$.
\item Si $c'$ est non vide, alors $P^c(H)=P^{c'}(F) P^{c''}(G)$ ; de plus, $s$ est un sommet de $R^{c''}(G),$  et $R^c(H)$ est la greffe de $R^{c'}(F)$ sur $R^{c''}(G)$.
\end{enumerate}
\item $c$ ne coupe que des ar\^etes de $G$ : posons $c'=c_{\mid G}$.
\begin{enumerate}
\item  Si $s$ est un sommet de $P^{c'}(G)$,  $P^{c}(H)$ est la greffe de $F$ sur le sommet $s$ de $P^{c'}(G)$ et $R^{c}(H)=R^{c'}(G)$.
\item Si $s$ est un sommet de $R^{c'}(G)$, alors $P^c(H)=P^{c'}(G)$, $R^{c}(H)$ est la greffe de $F$ sur le sommet $s$ de $R^{c'}(G)$.
\end{enumerate}
\end{enumerate}
En sommant sur $s$, et en posant $\tdelta(F)=\sum F'\otimes F''$, $\tdelta(G)= \sum G' \otimes G''$ :
\begin{eqnarray*}
n\tdelta( F \otop G)&=&n F \otimes G + \sum n'' FG'\otimes G''+\sum n F' \otimes (F'' \otop G)\\
&&+ \sum \sum n'' F'G' \otimes (F'' \otop G'') + \sum n' (F \otop G') \otimes G'' +\sum n'' G'\otimes (F \otop G'').
 \end{eqnarray*}
($n$ d\'esigne le poids de $G$, $n'$ le poids de $G'$, $n''$ le poids de $G''$.)\\

Soit $h_F:\nH\otimes \nH \longmapsto \nH\otimes \nH$ d\'efinie par :
\begin{eqnarray*}
h_F(X \otimes Y) &=& \frac{y}{x+y} FX \otimes Y +\frac{y}{x+y} \sum F' X \otimes (F'' \otop Y)\\
&&+ \frac{x}{x+y} (F \otop X)\otimes Y + \frac{y}{x+y} X \otimes (F \otop Y),
\end{eqnarray*}
o\`u $X,Y$ sont des \'el\'ements de $\nforets$ de poids respectifs $x,y$.  En remarquant que $n'+n''=n$, on obtient :
\begin{eqnarray*}
\tdelta(F \otop G) &=& F \otimes G + \sum F' \otimes (F'' \otop G) + h_F( \tdelta(G)).
\end{eqnarray*}
Par lin\'earit\'e, pour tout $p$ primitif de $\nH$ et $F \in \nforets$ :
$$ \tdelta(F \otop p)= F \otimes p + \sum F' \otimes (F'' \otop p),$$ 
ce qui est \'equivalent \`a la condition 2 du th\'eor\`eme \ref{theo13}. $\Box$
 \\

 La condition 3 d\'ecoule du lemme \ref{lemme12}. On en d\'eduit que $\nH$ est une cog\`ebre isomorphe \`a $T(Prim(\nH))$. Plus pr\'ecis\'ement, on pose :
\begin{eqnarray*}
\overline{F}_n\esp : Prim(\nH)^{\otimes n}& \longmapsto & \nH\\
p_1 \otimes \ldots \otimes p_n &\longmapsto &p_1\otop \ldots \otop p_n.
\end{eqnarray*}
Les $\overline{F}_n$ sont alors injectives, et $\nH=\bigoplus Im(\overline{F}_n)$.
\label{part10.2}
\subsection{Cas de l'abelianis\'ee de $\Hr$}

\spa  On d\'efinit $\breve{\top}:\Hr \otimes \Hr \longmapsto \Hr$ par r\'ecurrence sur $poids(F)$ de la mani\`ere suivante : 
\begin{eqnarray*}
t_1 \ldots t_n \breve{\top} 1  &=&0 \\
t_1 \ldots t_n \breve{\top} \bullet_d &=& \sum_{\sigma \in S_n} B_d^+( t_{\sigma(1)} \ldots t_{\sigma(n)}),\\
G \breve{\top} t'_1 \ldots t'_m &=& \sum_{i=1}^m  t'_1 \ldots (G \breve{\top} t'_i) \ldots t'_m,\\
t_1 \ldots t_n \breve{\top} B_d^+(F)&=&B^+_d(t_1 \ldots t_n\breve{\top}F) + \sum_{\sigma \in S_n} B^+_d (t_{\sigma(1)} \ldots t_{\sigma(n)} F).
\end{eqnarray*}
(C'est-\`a-dire qu'on greffe les $t_{\sigma(1)} \ldots t_{\sigma(n)}$ "le plus \`a gauche possible" sur chaque sommet de $F$).\\

On pose :
$$t_1 \ldots t_n \mtop F = \frac{1}{n! poids(F)} t_1\ldots t_n \breve{\top} F.$$

Soient $F$ et $G$ $\in \forets$, toutes deux diff\'erentes de $1$. On pose $\tdelta(F)=\sum F' \otimes F''$ et 
$\tdelta(G)=\sum G' \otimes G''$. Une \'etude simple des coupes admissibles de chaque for\^et de  $G \breve{\top} F$ montre que :
\begin{eqnarray}
\nonumber \Delta(F \mtop G)&=& (F\mtop G) \otimes 1 +1 \otimes (F \mtop G) +  F \otimes G + \sum G' \otimes (F'' \mtop G) \\
\label{eq29} &&+  h_F(\tdelta(G))+I \otimes \Hr,
\end{eqnarray}
o\`u $h_F:\Hr \otimes \Hr \longmapsto \Hr \otimes \Hr$ est une certaine application lin\'eaire
et $I$ le sous-espace de $\Hr$ engendr\'e par les $t_1 \ldots t_n -t_{\sigma(1)}\ldots t_{\sigma(n)}$, $t_1 \ldots t_n \in \forets$, $\sigma \in S_n$.
Or $I$ est l'id\'eal bilat\`ere engendr\'e par les $xy-yx$, $x,y \in \Hr$. De plus, on a facilement :
\begin{eqnarray*}
( t_1 \ldots t_n -t_{\sigma(1)}\ldots t_{\sigma(n)}) \mtop F&=&0,\\
G \mtop (t'_1 \ldots t'_m-t'_{\sigma(1)} \ldots t'_{\sigma(m)})&\in & I. 
\end{eqnarray*}
Par suite, $\mtop$ passe au quotient dans $(\Hr)_{ab} =\Hr/I$. Comme dans le cas de $\nH$, la condition 2 du th\'eor\`eme \ref{theo13} est v\'erifi\'ee d'apr\`es (\ref{eq29}). On en d\'eduit que la cog\`ebre $(\Hr)_{ab}$ est isomorphe \`a la cog\`ebre $T(Prim((\Hr)_{ab}))$.
\label{part10.3}

\subsection{Bigraduation de $T(V)$}
\spa Dans ce paragraphe, on suppose que $V$ est muni d'une graduation $(V_i)_{i\geq 1}$, telle que les $V_i$ soient de dimension finie ; on pose $v_i=dim(V_i)$ et $P(X)=\sum v_i X^i$ la s\'erie g\'en\'eratrice des $v_i$.
On suppose $v_0=0$. Cette graduation induit une graduation de la cog\`ebre $T(V)$. On notera $T(V)_{poids=n}$ la composante homog\`ene de degr\'e $n$ de $T(V)$, $r_n$ sa dimension, et $R(X)=\sum r_n X^n$ la s\'erie g\'en\'eratrice des $r_n$.\\
On notera $T(V)_{poids=n, deg_p=m}=T(V)_{poids=n} \cap V^{\otimes m}$, et $h_{n,m}=dim(T(V)_{poids=n,deg_p=m})$. Enfin, on introduit les s\'eries formelles suivantes :
$$H_m(X)=\sum_{n\geq 0} h_{n,m} X^{n}, \esp H(X,Y)=\sum_{n\geq 0, m\geq 0} h_{n,m} X^n Y^m.$$
On remarque que $H_0(X)=1$ et que $H_1(X)=P(X)$.
\begin{theo}
\begin{eqnarray*}
H_m(X)&=&P(X)^m, \esp \forall m \in \mathbb{N},\\
R(X)&=&\frac{1}{1-P(X)},\\
H(X,Y)&=&\frac{R(X)}{(1-Y)R(X)+Y}.
\end{eqnarray*}
\end{theo}
{\it Preuve :}
\begin{eqnarray*}
r_n&=&\sum_{k\geq 1} \sum_{\stackrel{\mbox{\scriptsize{  $a_1,\ldots,a_k$ tous non nuls,}}}{a_1+\ldots+a_k=n}}
v_{a_1} \ldots v_{a_k}\\
&=&\sum_{b_1+2b_2+\ldots+nb_n=n} \frac{(b_1+\ldots + b_n)!}{b_1!\ldots b_n!} v_1^{b_1} \ldots v_n^{b_n}.
\end{eqnarray*}
D'o\`u $R(X)=\sum_{i\geq 0} P(X)^i=(1-P(X))^{-1}$.
\begin{eqnarray*}
h_{n,m}&=&\sum_{\stackrel{\mbox{\scriptsize{ $a_1,\ldots,a_m$ tous non nuls,}}}{a_1+\ldots+a_m=n} }
v_{a_1} \ldots v_{a_m}\\
&=&\sum_{\stackrel{b_1+2b_2+\ldots+nb_n=n,}{b_1+b_2+\ldots+b_n=m}}  \frac{m!}{b_1!\ldots b_n!} v_1^{b_1} \ldots v_n^{b_n}.
\end{eqnarray*}
Par suite, $H_m(X)=H_1(X)^m=P(X)^m$.

On a de plus :
\begin{eqnarray*}
H(X,Y)&=&\sum_{m\geq 0} H_m(X)Y^m\\
  &=&\sum_{m\geq 0} \left(P(X)Y\right)^m\\
&=&\frac{1}{1-P(X)Y}\\
&=&\frac{1}{1-Y + \frac{Y}{R(X)}} \\
&=&\frac{R(X)}{(1-Y)R(X)+Y}. \esp \Box
\end{eqnarray*}

Supposons $D$ fini de cardinal $\cal D$. Les r\'esultats pr\'ec\'edents s'appliquent \`a $\nH$, $\Hr$ et $(\Hr)_{ab}$.  Dans le cas de $\nH$, on retrouve les r\'esultats de \cite{Kreimer,Foissy}. 
Dans le cas de$ \Hr$, on pose $p_n=dim(Prim(\Hr)\cap {\cal H}_n)$. On obtient :
\begin{eqnarray*}
P(X)&=&1-\frac{1}{R(X)}\\
    &=& 1- (1-T(DX))  \mbox{ d'apr\`es (\ref{eqnR})}\\
    &=& T(DX).
\end{eqnarray*}

\begin{prop}
\label{pro16}
Soit $p_n=dim(Prim(\Hr)\cap {\cal H}_n)$, et $P(X)=\sum p_nX^n.$
Alors :
$$ P(X)=\frac{1-\sqrt{1-4DX}}{2}, \esp  p_n = \frac{(2n-2)}{n!(n-1)!}D^n=D^n \tau_n.$$
\end{prop}

\section{Endomorphismes de $\Hr$ et comodules sur $\Hr$}
\spa On a vu dans la section \ref{part4} que pour $p_1,\ldots,p_n \in Prim(\Hr)$ :
$$\Delta(p_1\top \ldots \top p_n)=\sum_{i=0}^n (p_1\top \ldots \top p_i)\otimes(p_{i+1}\top \ldots \top p_n).$$
Pour $n \in \mathbb{N}$, on consid\`ere:
\begin{eqnarray*}
F_n:(Prim(\Hr))^{\otimes n}&\longmapsto&\Hr\\
p_1\otimes \ldots \otimes p_n &\longmapsto & p_1\top \ldots \top p_n.
\end{eqnarray*}
Les $F_n$ sont injectives, et $\Hr=\bigoplus Im(F_n)$.

\subsection{Endomorphismes de cog\`ebre}
\spa {\it Notations :} 
\begin{enumerate}
\item Soit $u:Prim(\Hr)^{\otimes i}\longmapsto Prim(\Hr)^{\otimes j}$. On d\'efinit $\overline{u}:Im(F_i)\longmapsto Im(F_j)$ par $\overline{u}=F_j\circ u \circ F_i^{-1}$.
\item $\pi_1$ est la projection sur $Im(F_1)=Prim(\Hr)$ dans la somme directe $\Hr=\bigoplus Im(F_i)$.
\end{enumerate}

\begin{theo}
\begin{enumerate}
\item  Pour tout $i\in \mathbb{N}^*$, soit $u_i: Prim(\Hr)^{\otimes i}\longmapsto Prim(\Hr)$. On d\'efinit $\Phi_{(u_i)}:\Hr \longmapsto \Hr$ par:
\begin{eqnarray*}
\Phi_{(u_i)}(1)&=&1, \\
\Phi_{(u_i)}(p_1 \top \ldots \top p_n)&=&\sum_{k=1}^n \esp \sum_{\stackrel{a_1+\ldots +a_k=n,}{a_i>0}} (\overline{u_{a_1} \otimes \ldots \otimes u_{a_k}})(p_1\top \ldots \top p_n).
\end{eqnarray*}
Alors $\Phi_{(u_i)}$ est un endomorphisme de cog\`ebre.
\item  Soit $\Phi$ un endomorphisme de cog\`ebre de $\Hr$. Alors il existe une unique famille $(u_i)_{i \in \mathbb{N}^*}$, $u_i:Prim(\Hr)^{\otimes i} \longmapsto Prim(\Hr)$,
telle que $\Phi=\Phi_{(u_i)}$.
\end{enumerate}
\end{theo}
{\it Preuve :}

1. On a imm\'ediatement $\varepsilon \circ \Phi_{(u_i)}=\varepsilon$. Comme $\Phi_{(u_i)}(1)=1$, il suffit de montrer:
$$\tdelta(\Phi_{(u_i)}(p_1 \top \ldots \top p_n))=\Phi_{(u_i)} \otimes \Phi_{(u_i)}(\tdelta(p_1\top \ldots \top p_n)).$$
On a:
$$\Phi_{(u_i)} \otimes \Phi_{(u_i)} (\tilde{\Delta}(p_1 \top \ldots \top p_n))=$$
$$\hspace*{-0.5cm} \sum_{j=0}^n \esp \sum_{a_1+ \ldots +a_k  =  j} \esp
\sum_{b_1 + \ldots +b_l  = n-j} \left[ (\overline{u_{a_1} \otimes \ldots \otimes u_{a_k}}\,) \otimes (\overline{u_{b_1} \otimes \ldots \otimes u_{b_l}}\,) \right]
                \left[ (p_1 \top \ldots \top p_{j} ) \otimes ( p_{j+1} \top \ldots \top p_n ) \right] $$
$$= \tilde{\Delta}\left( \sum_{d_1 + \ldots + d_m = n} (\overline{ u_{d_1} \otimes \ldots \otimes u_{d_m} }\,) ( p_1 \top \ldots \top p_n) 
          - \overline{u_n}(p_1 \top \ldots \top p_n) \right)$$
$$= \tilde{\Delta}\left( \sum_{d_1 + \ldots + d_m = n} (\overline{ u_{d_1} \otimes \ldots \otimes u_{d_m} }\,) ( p_1 \top \ldots \top p_n) 
         \right)-0$$
$$= \tilde{\Delta}(\Phi_{(u_i)}(p_1 \top \ldots \top p_n)).$$

2. Montrons par r\'ecurrence sur $n$ qu'il existe $u_i:Prim(\Hr)^{\otimes i} \longmapsto Prim(\Hr)$ pour tout $i\leq n$, tel que si on pose
$u_i^{(n)}=u_i $ si $i \leq n$ et $u_i^{(n)}=0$ si $i>n$, alors $\Phi=\Phi_{(u_i^{(n)})}$ sur $Im(F_0)\oplus \ldots \oplus Im(F_n)$. 

Comme $\Phi(1)$ est un \'el\'ement groupo\"idal, n\'ecessairement $\Phi(1)=1$, ce qui prouve la propri\'et\'e au rang $0$. Supposons la propri\'et\'e vraie au rang $n-1$.
On pose $\Phi^{(n-1)}=\Phi_{(u_i^{(n-1)})}$. On a alors:
\begin{eqnarray*}
\tdelta(\Phi(p_1\top \ldots \top p_n))&=&\Phi \otimes \Phi (\tdelta(p_1\top \ldots \top p_n))\\
&=&\Phi^{(n-1)} \otimes \Phi^{(n-1)} (\tdelta(p_1\top \ldots \top p_n))\\
&=&\tdelta(\Phi^{(n-1)}(p_1\top \ldots \top p_n)).
\end{eqnarray*}
(On a utilis\'e le fait que $\tdelta(p_1\top \ldots \top p_n) \in Im(F_1) \otimes Im(F_{n-1}) + \ldots + Im(F_{n-1}) \otimes Im(F_1)$ pour la deuxi\`eme \'egalit\'e.)\\

Donc $(\Phi-\Phi^{(n-1)})(p_1\top \ldots \top p_n)$ est primitif.
On d\'efinit alors $u_n: Prim(\Hr)^{\otimes n}\longmapsto Prim(\Hr)$ par: $$\overline{u}_n(p_1\top \ldots \top p_n)=(\Phi-\Phi^{(n-1)})(p_1\top \ldots \top p_n).$$
Comme $\Phi_{(u_i^{(n)})}=\Phi_{(u_i^{(n-1)})}$ sur $Im(F_0) \oplus \ldots \oplus Im(F_{n-1})$, et que $\Phi_{(u_i^{(n)})}=\Phi_{(u_i^{(n-1)})}+\overline{u}_n$ sur $Im(F_n)$, on a $\Phi=\Phi_{(u_i^{(n)})}$ sur $Im(F_0) \oplus \ldots \oplus Im(F_n)$.

On conclut en remarquant que $\Phi_{(u_i^{(n+m)})} = \Phi_{(u_i^{(n)})}$ sur $Im(F_0)\oplus \ldots \oplus Im(F_n)$, et donc  $\Phi=\Phi_{(u_i)}$ sur $\bigoplus Im(F_i)=\Hr$.\\

Unicit\'e des $u_i$ : on a $\overline{u}_i=\pi_1 \circ \Phi_{\mid Im(F_i)}$. $\Box$

\begin{cor}
\begin{eqnarray*}
End_{\mbox{\scriptsize{cog\`ebre}}}(\Hr)&\longmapsto & {\cal L}(\Hr, Prim(\Hr)) \mbox{ est une bijection.}\\
\Phi&\longmapsto & \pi_1\circ \Phi
\end{eqnarray*}
\end{cor}
{\it Preuve :} injectivit\'e : supposons que $\pi_1\circ \Phi=\pi_1 \circ \Phi'$.
Montrons que $\Phi(p_1\top \ldots \top p_n)=\Phi'(p_1\top \ldots \top p_n)$ par r\'ecurrence sur $n$.
Pour $n=0$, $\Phi(1)=\Phi'(1)=1$. Supposons la propri\'et\'e vraie pour tout $k<n$.
\begin{eqnarray*}
\tdelta \circ \Phi(p_1\top \ldots \top p_n)&=& \Phi \otimes \Phi(\tdelta(p_1\top \ldots \top p_n))\\
&=& \Phi' \otimes \Phi'(\tdelta(p_1\top \ldots \top p_n))\\
&=&\tdelta \circ \Phi'(p_1\top \ldots \top p_n).
\end{eqnarray*}
Donc $(\Phi-\Phi')(p_1\top \ldots \top p_n)$ est primitif, par suite :
$$(\Phi-\Phi')(p_1\top \ldots \top p_n)=\pi_1\circ(\Phi-\Phi')(p_1\top \ldots \top p_n)=0.$$

Surjectivit\'e : soit $u\in {\cal L}(\Hr, Prim(\Hr))$. Soit $u_i$ telle  que $\overline{u}_i=u_{\mid Im(F_i)}$. Alors $\pi_1 \circ \Phi_{(u_i)}=u$. $\Box$

\begin{cor}
Soit $\Phi \in End_{\mbox{\scriptsize{cog\`ebre}}}(\Hr)$. Alors $\Phi(Prim(\Hr))\subseteq Prim(\Hr)$ ; soit $\phi: Prim(\Hr)\longmapsto Prim(\Hr)$ d\'efini par $\phi(p)=\Phi(p)$, $\forall p \in Prim(\Hr)$. Alors $\Phi$ est injectif (respectivement bijectif) si et seulement si $\phi$ est injectif (respectivement  bijectif). Si $\phi$ est surjectif, alors $\Phi$ est surjectif.
\end{cor}
{\it Preuve :} on peut supposer $\Phi$ de la forme  $\Phi_{(u_i)}$. Alors $\phi=u_1$.

Injectivit\'e : $\Rightarrow$ : \'evident.

$\Leftarrow$ : soit $x \in Ker(\Phi)$, $x \neq 0$. On peut \'ecrire $x=x_n+y$, avec $x_n \in Im(F_n)$, non nul, et $deg_p(y)<n$. Alors:
$$\Phi(x)=(\overline{u_1\otimes \ldots \otimes u_1})(x_n)+ Im(F_0)\oplus \ldots \oplus Im(F_{n-1}) =0.$$ 
Donc comme $(\overline{u_1\otimes \ldots \otimes u_1})(x_n) \in Im(F_n)$, on a $(\overline{u_1\otimes \ldots \otimes u_1})(x_n)=0$. Or $u_1$ est injectif, donc
$u_1 \otimes \ldots \otimes u_1$ est injectif, et donc $\overline{u_1\otimes \ldots \otimes u_1}$ est injectif : par suite $x_n=0$ : on aboutit \`a une contradiction. 
Donc $\Phi$ est injectif.\\

Surjectivit\'e : $\Leftarrow$ : supposons $u_1$ surjectif. Montrons que $Im(F_0)\oplus \ldots \oplus Im(F_n) \subset Im(\Phi)$ $\forall n$. Si $n=0$, c'est \'evident car $\Phi(1)=1$. Supposons la propri\'et\'e vraie au rang $n-1$. Soient $p_1, \ldots, p_n \in Prim(\Hr)$. Il existe $q_1, \ldots, q_n \in prim(\Hr)$, tels que $p_i=u_1(q_i)$.
Alors:
$$\Phi(q_1\top \ldots \top q_n)=p_1\top \ldots \top p_n + Im(F_0)\oplus\ldots \oplus Im(F_{n-1}).$$
Donc $p_1 \top \ldots \top p_n\in Im(\Phi)$, ce qui prouve la propri\'et\'e au rang $n$.\\

Bijectivit\'e : $\Leftarrow$ : d\'ecoule imm\'ediatement de ce qui pr\'ec\`ede.

$\Rightarrow$ : supposons $\Phi=\Phi_{(u_i)}$  bijectif. Son inverse est un morphisme de cog\`ebres, donc de la forme $\Phi_{(v_i)}$. On a alors :
\begin{eqnarray*}
(\Phi \circ \Phi^{-1})_{\mid Prim(\Hr)}&=& Id_{Prim(\Hr)}\\
&=&\Phi \circ v_1\\
&=&u_1 \circ v_1.
\end{eqnarray*}
De m\^eme, $v_1 \circ u_1=Id_{Prim(\Hr)}$, donc $u_1$ est bijectif. $\Box$ \\

{\it Remarque : } 
on peut avoir $\Phi$ sujectif sans que $\phi$ le soit. Par exemple, pour $u_1,$ $u_3$, $u_4\ldots$ nuls, et $u_2:Prim(\Hr)^{\otimes 2} \longmapsto Prim(\Hr)$ surjectif, on a:
\begin{eqnarray*}
\Phi_{(u_i)}(p_1 \top \ldots \top p_{2n+1})&=&0,\\
\Phi_{(u_i)}(p_1 \top \ldots \top p_{2n})&=&\underbrace{\overline{u_2\otimes \ldots \otimes u_2}}_{ \mbox{$n$ fois}}(p_1 \top \ldots \top p_{2n}).
\end{eqnarray*}
Comme $u_2$ est surjectif, $u_2^{\otimes i}:Prim(\Hr)^{\otimes 2i}\longmapsto Prim(\Hr)^{\otimes i}$ est surjectif, et donc
$\overline{u_2^{\otimes i}}:Im(F_{2i})\longmapsto Im(F_i)$ est surjectif. Donc $\Phi_{(u_i)}$ est surjectif. Cependant, $\phi_{(u_i)}=u_1=0$ n'est pas surjectif.
   
\subsection{Endomorphismes de big\`ebre}
 
 \begin{theo}
\begin{enumerate}
\item Soit $(P_t)_{t \in \trees}$ une famille d'\'el\'ements primitifs de  $\Hr$ index\'ee par $\trees$.
Soit $\Phi_{(P_t)}$ l'unique endomorphisme d'alg\`ebre de $\Hr$ d\'efini par r\'ecurrence sur le poids par:
\begin{eqnarray*}
 \Phi_{(P_t)}(\bullet_d)&=&P_{\bullet_d}, \\
 \Phi_{(P_t)}(t)&=&\left(\sum_{(t)} \Phi_{(P_t)}(t') \top P_{t''}\right)+P_t \esp \mbox{ pour tout }t \in \trees \mbox{, avec } \tilde{\Delta}(t)=\sum_{(t)} t' \otimes t'' .
\end{eqnarray*}
Alors $\Phi_{(P_t)}$ est un endomorphisme de big\`ebre.
\item Soit $\Phi$ un endomorphisme de big\`ebre de $\Hr$ ; alors il existe une unique famille $(P_t)_{t\in \trees}$ telle que $\Phi=\Phi_{(P_t)}$.
\end{enumerate}
\end{theo}
{\it Preuve :} remarquons que $\Phi_{(P_t)}$ est bien d\'efini, car $\Hr$ est librement engendr\'ee par $\trees$.

1. Il s'agit de montrer que $\tdelta \circ \Phi_{(P_t)}(t)=\Phi_{(P_t)}\circ \Phi_{(P_t)}(\tdelta(t))$ $\forall t \in \trees$.
Proc\'edons par r\'ecurrence sur le poids de $t$. Si $t$ est de poids $1$, alors $t$ est de la forme $\bullet_d$, et donc $\Phi_{(P_t)}(t)=P_t$ ; comme $t$ est $P_t$ sont tous les deux primitifs, la propri\'et\'e est v\'erifi\'ee. Supposons la propri\'et\'e vraie pour tout arbre de poids strictement inf\'erieur \`a $n$. Comme $\Phi_{(P_t)}$ est un morphisme d'alg\`ebres, la propri\'et\'e est vraie pour tout \'el\'ement de $\Hr$ de poids strictement inf\'erieur \`a $n$. On a $\tdelta(t)=\sum t'\otimes t''$, avec $t''\in \trees$.

\begin{eqnarray*}
\tilde{\Delta}(\Phi_{(P_t)}(t))&=&\sum_{(t)} \tilde{\Delta}\left(\Phi_{(P_t)}(t') \top P_{t''}\right)\\
& =&\sum_{(t)} \Phi_{(P_t)}(t') \otimes P_{t''}+ \sum_{(t)} \sum_{(\Phi(t'))} \Phi_{(P_t)}(t')' \otimes (\Phi_{(P_t)}(t')'' \top P_{t''})\\
& =&\sum_{(t)} \Phi_{(P_t)}(t') \otimes P_{t''}+ \sum_{(t)} \sum_{(t')} \Phi_{(P_t)}((t')') \otimes (\Phi_{(P_t)}((t')'') \top P_{t''})\\
& =&\sum_{(t)} \Phi_{(P_t)}(t') \otimes P_{t''}+ \sum_{(t)} \sum_{(t'')} \Phi_{(P_t)}(t') \otimes (\Phi_{(P_t)}((t'')') \top P_{(t'')''})\\
&=&\sum_{(t)} \Phi_{(P_t)}(t') \otimes \left[ \sum_{(t'')}\left(\Phi_{(P_t)}{((t'')}') \top P_{{(t'')}''}\right) + P_t'' \right]\\
& =&\sum_{(t)} \Phi_{(P_t)}(t') \otimes \Phi_{(P_t)}(t'').
\end{eqnarray*}
(On a utilis\'e le fait que $x \longmapsto x\top p$ soit un $1-$cocycle pour tout primitif $p$ pour la deuxi\`eme \'egalit\'e, l'hypoth\`ese de r\'ecurrence pour la troisi\`eme \'egalit\'e et la coassociativit\'e de $\tdelta$ pour la quatri\`eme.)

2. Montrons qu'il existe $(P_t)_{poids(t)\leq n}$ telle que si on pose $P_t^{(n)}=P_t$ si  $ poids(t)\leq n$, et $P_t^{(n)}=0$ si $poids(t)>n$, alors $\Phi(t)=\Phi_{(P_t^{(n)})}(t)$  pour tout $t$ de poids inf\'erieur ou \'egal \`a $n$. Pour $n=1$, alors on prend $P_{\bullet_d}=\Phi(\bullet_d)$. Supposons la propri\'et\'e vraie au rang $n-1$. Posons $\Phi^{(n)}=\Phi_{(P_t^{(n)})}$. Soit $t\in \trees$, $poids(t)=n$.
\begin{eqnarray*}
\tdelta(\Phi(t))&=&\Phi\otimes \Phi(\tdelta(t))\\
&=&\Phi^{(n-1)} \otimes \Phi^{(n-1)}(\tdelta(t))\\
&=& \tdelta(\Phi^{(n-1)}(t)).
\end{eqnarray*}
On prend alors $P_t=\Phi(t)-\Phi^{(n)}(t)$. On a bien $P_t \in prim(\Hr)$, et $\Phi(t)=\Phi^{(n)}(t)$ pour tout $t\in \trees$ de poids inf\'erieur ou \'egal \`a $n$.\\
Comme $\Phi^{(n+m)}(t)=\Phi^{(n)}(t)$ pour tout $t\in \trees$ de poids inf\'erieur ou \'egal \`a $n$, on a $\Phi=\Phi_{(P_t)}$.

Unicit\'e : on a $P_t= \pi_1 \circ \Phi_{(P_t)}(t)$ pour tout $t \in \trees$. $\Box$

\begin{cor}
Soit $\cal T$ le sous-espace de $\Hr$ engendr\'e par les \'el\'ements de $\trees$. Alors :
\begin{eqnarray*}
End_{\mbox{\scriptsize{big\`ebre}}}(\Hr)&\longmapsto & {\cal L}({\cal T}, Prim(\Hr)) \mbox{ est une bijection.}\\
\Phi&\longmapsto & \pi_1\circ \Phi_{\mid {\cal T}}
\end{eqnarray*}
\end{cor}
{\it Preuve :} 

Surjectivit\'e : soit $u\in  {\cal L}({\cal T}, Prim(\Hr))$. On pose $P_t=u(t)$ pour tout $t \in \trees$. Alors $\pi_1\circ \Phi_{(P_t)}=u$ sur $\cal T$.

Injectivit\'e : si $\pi_1 \circ \Phi_{(P_t)}=\pi_1 \circ \Phi_{(P'_t)}$ sur $\cal T$, alors pour tout $t \in \trees$, 
\begin{eqnarray*}
\pi_1 \circ \Phi_{(P_t)}(t)&=&P_t\\
&=&\pi_1 \circ \Phi_{(P'_t)}(t)\\
&=&P'_t,
\end{eqnarray*}
et donc $\Phi_{(P_t)}=\Phi_{(P'_t)}$. $\Box$

\begin{cor}
Soit $\Phi \in End_{\mbox{\scriptsize{big\`ebre}}}(\Hr)$ ; alors $\Phi$ est un endomorphisme d'alg\`ebre de Hopf, c'est-\`a-dire $\Phi \circ S=S\circ \Phi$.
\end{cor}
{\it Preuve :} on suppose $\Phi=\Phi_{(P_t)}$.
Il faut montrer $S \circ \Phi(t)=\Phi \circ S(t)$, $\forall t \in \trees$. Si $t=\bullet_d$, alors 
$S \circ \Phi(t)=\Phi \circ S(t)=-\Phi(t)$ car $t$ est primitif. Supposons cette propri\'et\'e vraie pour tout $t$ de poids inf\'erieur ou \'egal \`a $n$. \\

Soit $x \in \Hr$ ; posons $\tdelta(x)=\sum x'\otimes x''$. Pour tout $p$ primitif, on a alors :
\begin{eqnarray*}
m \circ (S \otimes Id) \circ \Delta(x\top p )&=& S(x\top p) + x\top p +S(x) p  +\sum S(x') (x''\top p)\\
&=& \varepsilon(x \top p)1\\
&=&0.
\end{eqnarray*}
Donc $S(x \top p)=-x \top p-\sum S(x') (x''\top p) -S(x) p .$ 
\\
En particulier, si $x = \bullet_d$ : 
$S \circ B_d^+ (x)=-B_d^+(x)-\sum S(x') B^+_d(x'')-S(x)\bullet_d $.\\

Soit $t\in \trees$, $poids(t)=n$. Posons $t=B_d^+(F)$. On a alors :
\begin{eqnarray*}
  \Phi \circ S(B_d^+(F))&=&-\Phi(F) \top P_{\bullet_d} - \sum \Phi(F') \top P_{B_d^+(F'')}-P_{B_d^+(F)} \\
&&- \sum \Phi \circ S (F') (\Phi(F'') \top P_{\bullet_d})
- \sum \Phi \circ S(F') [\Phi(F'') \top P_{B_d^+(F''')}]\\
&& -\sum \Phi \circ S(F') P_{B_d^+(F')}-\Phi \circ S(F)P_{\bullet_d}   \esp ;\\
 S \circ \Phi(B_d^+(F))&=&S(\Phi(F) \top P_{\bullet_d} ) + \sum S( \Phi(F') \top P_{B_d^+(F'')}) + S(P_{B_d^+(F)})\\
&=& -\Phi(F) \top P_{\bullet_d} - \sum S \circ \Phi(F') ( \Phi(F'') \top P_{\bullet_d})  - S \circ \Phi(F) P_{\bullet_d} \\
&& -\sum \Phi(F') \top P_{B_d^+(F'')} - \sum S \circ \Phi(F') [ \Phi(F'') \top P_{B_d^+(F''')}]\\
&&  - \sum_{(F)} S \circ \Phi(F') P_{B_d^+(F'')}- P_{B_d^+(F)} . 
\end{eqnarray*}
L'hypoth\`ese de r\'ecurrence permet imm\'ediatement de conclure. $\Box$

\subsection{Comodules sur $\Hr$}
\spa $\Hr$ et ${\cal H}_R$ \'etant toutes deux isomorphes \`a une cog\`ebre tensorielle, le th\'eor\`eme 3.8 et les r\'esultat des sections 5 et 6 de \cite{Foissy} restent vrais dans le cas de $\Hr$.

\section{Compl\'ements sur la cog\`ebre $T(V)$}
\spa On suppose maintenant qu'on a dot\'e $T(V)$ d'un produit, de sorte que $(T(V),m,\eta,\Delta,\varepsilon)$ soit une big\`ebre.
Le but de cette partie est de montrer que si $m$ est commutatif, alors cette big\`ebre est isomorphe \`a l'alg\`ebre de battage introduite dans la section \ref{part12.3}.

L'\'el\'ement neutre pour l'alg\`ebre v\'erifie $\Delta(e)=e \otimes e$ ; il s'agit donc de $1 \in V^{\otimes 0}$ d'apr\`es la proposition \ref{prop14}-1.
\subsection{Existence d'une antipode} 
\spa Soient $p_1 ,\ldots, p_n \in V$. Soit $c$ une partie de $\{1 ,\ldots, n-1\}$. on note $(p_1 \top \ldots \top p_n)_c$ l'\'el\'ement 
de $T(V)$ obtenu en remplaceant le $i$-\`eme $\top$ de $p_1 \top \ldots \top p_n$ par un produit pour tout $i\in c$,
les $\top$ \'etant index\'es de la gauche vers la droite.\\

{\it Exemple :} \begin{eqnarray*}
(p_1 \top p_2 \top p_3)_{\emptyset}&=&p_1 \top p_2 \top p_3 \esp ;\\
(p_1 \top p_2 \top p_3)_{\{1\}}&=&p_1 (p_2 \top p_3) \esp ;\\
(p_1 \top p_2 \top p_3)_{\{2\}}&=&(p_1 \top p_2) p_3 \esp ;\\
(p_1 \top p_2 \top p_3)_{\{1,2\}}&=&p_1  p_2  p_3.
\end{eqnarray*}
  
\begin{theo}
\label{theo38}
pour $p_1, \ldots, p_n \in V$, $n\geq 1$ on d\'efinit :
\begin{eqnarray*}
S(p_1 \top \ldots \top p_n)&=&-\sum_{c \subseteq \{1, \ldots ,n-1\}} (-1)^{card(c)} (p_1 \top \ldots \top p_n)_c \esp ;\\
S(1)&=&1.
\end{eqnarray*}
Alors $S$ est une antipode pour $(T(V),m,\eta,\Delta,\varepsilon)$.
 \end{theo}
{\it Preuve :} soit $c$ une partie non vide de $\{1, \ldots, n-1\}$ ; soit $j=\max (c)$.
Alors $(p_1 \top \ldots \top p_n)_c=(p_1 \top \ldots \top p_{j})_{c-\{j\}}(p_{j+1} \top\ldots \top j_n).$ 
Posons $c'=c-\{j\}$ ; alors $c' \subseteq \{1, \ldots,j-1\}$.
D'o\`u :
\begin{eqnarray*}
 S(p_1\top \ldots \top p_n)&=&-(p_1 \top \ldots \top p_n)_{\emptyset} \\
 &&-\sum_{j=1}^{j=n-1} \sum_{c' \subseteq \{1, \ldots , j-1 \}}
      (-1)^{card(c')+1} (p_1 \top \ldots \top p_j)_{c'}(p_{j+1} \top \ldots \top p_n)\\
 &=&-(p_1 \top \ldots \top p_n) -\sum_{j=1}^{j=n} S(p_1 \top \ldots p_j)(p_{j+1} \top \ldots \top p_n).
\end{eqnarray*}
Alors :
\begin{eqnarray*}
m \circ (S \otimes Id) \circ \Delta(p_1 \top \ldots \top p_n)
&=& S(p_1 \top \ldots \top p_n) + p_1\top \ldots \top p_n \\
&&+ \sum_{j=1}^{n} S(p_1\top \ldots \top p_j) (p_{j+1} \top \ldots \top p_n)\\
&=&0\\
&=&\varepsilon(p_1\top \ldots \top p_n)1.
\end{eqnarray*}
On montre de m\^eme que $m \circ (Id\otimes S) \circ \Delta= \eta \circ \varepsilon$, et donc $S$ est une antipode. $\Box$\\

\label{partie71}

\subsection{1-cocycles de $T(V)$}
\begin{prop}
Soit $u : T(V) \longmapsto V$ une application lin\'eaire quelconque. On d\'efinit $L_u :T(V) \longmapsto T(V)$ par :
\begin{eqnarray*}
L_u(1)&=&u(1) \esp ;\\
L_u(v_1 \top \ldots \top v_n)&=&\sum_{j=1}^{n-1} v_1 \top \ldots \top v_j \top u(v_{j+1} \top \ldots \top v_n)\\
&&   + v_1 \top \ldots \top v_n \top u(1) + u(v_1 \top \ldots \top v_n).
\end{eqnarray*}
Alors $L_u$ est un 1-cocycle de $(T(V),m,\eta,\Delta,\varepsilon)$.
\end{prop}
{\it Preuve :}
\begin{eqnarray*}
\Delta(L_u(1))&=&\Delta(u(1))\\
&=&u(1) \otimes 1 + 1 \otimes u(1)\\
&=&L_u(1) \otimes 1 + (Id \otimes L_u) \circ \Delta(1).
\end{eqnarray*}
Si $n \geq 1$ :
\begin{eqnarray*}
\Delta(L_u(v_1 \top \ldots \top v_n))&=& L_u(v_1 \top \ldots \top v_n)\otimes 1 + 1 \otimes L_u(v_1 \top \ldots \top v_n)\\
&&+\sum_{j=1}^{n-1} \sum_{k=1}^{j} (v_1 \top \ldots \top v_k) \otimes (v_{k+1} \top \ldots \top v_j \top u(v_{j+1} \top \ldots \top v_n))\\
&&+ \sum_{k=1}^{n} (v_1 \top \ldots \top v_k) \otimes (v_{k+1} \top \ldots \top v_n \top u(1))\\
&=&L_u(v_1 \top \ldots \top v_n)\otimes 1 + 1 \otimes L_u(v_1 \top \ldots \top v_n)  \\ 
&&+\sum_{k=1}^{n-1} (v_1 \top \ldots \top v_k) \otimes [\sum_{j=k+1}^{n-1}(v_{k+1} \top \ldots \top v_j \top u(v_{j+1} \top \ldots \top v_n)) \\
&& + u(v_{k+1} \top \ldots \top v_n)+ v_{k+1}\top \ldots \top v_n \top u(1)]+(v_1 \top \ldots \top v_n) \otimes u(1)  \\   
&=&L_u(v_1 \top \ldots \top v_n)\otimes 1+ (Id \otimes L_u)\left[ (v_1\top \ldots \top v_n) \otimes 1\right.\\
&& + 1\otimes (v_1\top \ldots \top v_n) + \sum_{k=1}^{n-1} (v_1 \top \ldots \top v_k) \otimes (v_{k+1} \top \ldots \top v_n)]\\
&=&L_u(v_1 \top \ldots \top v_n)\otimes 1+ (Id \otimes L_u) \circ \Delta(v_1\top \ldots \top v_n). \esp \Box
\end{eqnarray*}

Soit $\pi_1$ la projection de $T(V)$ sur $V$ dans la somme directe $K \oplus V \oplus V^2\oplus \ldots$

\begin{theo}
$$\begin{tabular}{rcl}
\hspace{-5cm}Soit $\Phi$ : $\cocy(T(V)^*)$ & $\longmapsto$ & ${\cal L}(T(V),V)$\\
         $L$ & $\longmapsto$ & $\pi_1 \circ L.$
\end{tabular}$$
Alors $\Phi$ est un isomorphisme d'espaces vectoriels. son inverse est donn\'e par $u \longmapsto L_u$. 
\end{theo}
{\it Preuve :} clairement, $\pi_1 \circ L_u=u$, $\forall u \in {\cal L}(T(V),V)$. Donc $\Phi$ est surjectif. Soit $L \in  \cocy(T(V)^*)$, tel que $\pi_1 \circ L=0$.
$L(1)$ est primitif, donc dans $V$ ; par suite, $L(1)=\pi_1 \circ L(1)=0$. Supposons que $L$ s'annule sur $K \oplus \ldots \oplus V^{\otimes n-1}$.
\begin{eqnarray*}
\Delta(L(v_1 \top \ldots \top v_n))&=&L(v_1 \top \ldots \top v_n) \otimes 1+1 \otimes L(v_1 \top \ldots \top v_n) \\
&&+\sum_{i=1}^{n} (v_1\top \ldots \top v_i) \otimes L(v_{i+1} \top \ldots \top v_n).
\end{eqnarray*}
D'apr\`es l'hypoth\`ese de r\'ecurrence, on a $L(v_1 \top \ldots \top v_n)$ primitif, donc dans $V$. Comme $\pi_1 \circ L=0$,
 $L(v_1 \top \ldots \top v_n)=0$. $\Box$ 
\begin{cor}
\label{cor41}
Soit $L :K \oplus \ldots \oplus V^{\otimes n} \longmapsto T(V)$, v\'erifiant  (\ref{eqn9}). Alors il existe $\overline{L} \in \cocy(T(V)^*)$, 
tel que $\overline{L}(x) = L(x)$ pour tout $x \in K \oplus \ldots \oplus V^{\otimes n} $.
 \end{cor}
{\it Preuve :} soit $u \in {\cal L}(T(V),V)$, tel que $u(x)=\pi_1 \circ L(x)$ pour tout $x \in K \oplus \ldots \oplus V^{\otimes n} $.
En reprenant la preuve du th\'eor\`eme pr\'ec\'edent, on montre que $L_u$ convient. $\Box$

\subsection{Produit de battage $*$ sur $T(V)$}
\label{part12.3}
\spa On d\'efinit un produit sur $T(V)$ en posant :
\begin{eqnarray}
\label{eqn29}
(v_1\top \ldots \top v_p)*(v_{p+1} \top \ldots \top v_{p+q})&=& \sum_{\sigma \in bat(p,q)} v_{\sigma(1)} \top \ldots \top v_{\sigma(p+q)}.
\end{eqnarray}
Ce produit a un \'el\'ement neutre : $1\in V^{\otimes 0}$.
\begin{theo}
$(T(V),*,1, \Delta, \varepsilon)$ est une alg\`ebre de Hopf commutative. Son antipode est donn\'ee par :
$$S_*(v_1\top \ldots \top v_n)=(-1)^n v_n \top \ldots \top v_1.$$
$(V^{\otimes n})_{n \in \mathbb{N}}$ est une graduation de l'alg\`ebre de Hopf $(T(V),*,1,\Delta,\varepsilon,S_*)$.
\end{theo} 
{\it Preuve :}

Associativit\'e de $*$ : on note $bat(p,q,r)$ l'ensemble des $\sigma \in S_{p+q+r}$, tels que $\sigma$ soit croissant sur
$\{1, \ldots, p\}$, $\{p+1, \ldots, p+q\}$ et $\{p+q+1, \ldots, p+q+r\}$.
Alors :
\begin{eqnarray*} 
&&\left[ (v_1 \top \ldots \top v_p)*(v_{p+1} \top \ldots \top v_{p+q})\right]*(v_{p+q+1} \top \ldots \top v_{p+q+r})\\
&=&\sum_{\sigma \in bat(p,q,r)} v_{\sigma(1)} \top \ldots \top  v_{\sigma(p+q+r)}\\
&=&(v_1 \top \ldots \top v_p)*\left[ (v_{p+1} \top \ldots \top v_{p+q})*(v_{p+q+1} \top \ldots \top v_{p+q+r})\right].
\end{eqnarray*}

$*$ est commutatif : on le montre de la m\^eme mani\`ere que la commutativit\'e du gradu\'e associ\'e \`a la filtration par $deg_p$ dans la preuve du corollaire \ref{corM}. 

$\varepsilon$ morphisme d'alg\`ebres : imm\'ediat.

$\Delta$ morphisme d'alg\`ebres :
Remarquons que :
\begin{eqnarray*}
(v_1\top \ldots \top v_p)*(v_{p+1} \top \ldots \top v_{p+q})&=& [(v_1 \top \ldots v_{p})*(v_{p+1} \top \ldots \top v_{p+q-1})]\top v_{p+q}\\
&&+[(v_1 \top \ldots v_{p-1})*(v_{p+1} \top \ldots \top v_{p+q})]\top v_p.
\end{eqnarray*}
D'o\`u, si $v,w \in V$, $x,y \in T(V)$, avec les notations de la proposition \ref{prop14}-2 :
\begin{eqnarray}
\label{eqn30} 
L_v(x)*L_w(y)&=&L_w\left(L_v(x)*y\right)+L_v\left(x *L_w(y)\right).
\end{eqnarray}
Montrons que $\Delta(a*b)=\Delta(a)*\Delta(b)$ par r\'ecurrence sur $n=deg_p(a)+deg_p(b)$.
Si $n=0$, alors $a$ et $b$ sont des constantes, et le r\'esultat est imm\'ediat.
Supposons la propri\'et\'e vraie au rang $n-1$. C'est \'evident si $deg_p(a)=0$ ou $deg_p(b)=0$.
Sinon, on peut se ramener au cas o\`u $a=L_v(x)$, $b=L_w(y)$. Notons que $deg_p(x)=deg_p(a)-1$, et $deg_p(y)=deg_p(b)-1$.
On pose 
$\Delta(x)=\sum x'\otimes x'', \esp \Delta(y)=\sum y' \otimes y''.$
Alors :
\begin{eqnarray*}
\Delta(L_v(x)*L_w(y))&=&(a*b)\otimes 1\\
&&+(Id\otimes L_w)\left( \sum (L_v(x)*y')\otimes y''+\sum \sum (x'*y')\otimes (L_v(x'')*y'')\right)\\
&&+(Id\otimes L_v)\left(\sum (x'*L_w(y))\otimes x''+\sum \sum (x'*y')\otimes (x'' *L_w(y''))\right)\\
&=&(a*b)\otimes 1+ \sum (L_v(x)*y')\otimes L_w(y'')+\sum (x'*L_w(y))\otimes L_v(x'')\\
&&+\sum \sum (x'*y')\otimes \left[L_w(L_v(x'')*y'')+L_v(x'' *L_w(y''))\right]\\
&=&(a*b)\otimes 1+ \sum (L_v(x)*y')\otimes L_w(y'')+\sum (x'*L_w(y))\otimes L_v(x'')\\  
&&+\sum \sum (x'*y') \otimes \left[L_v(x'')*L_w(y'')\right]\\
&=&\left(L_v(x)\otimes 1+\sum x' \otimes L_v(x'')\right)*\left(L_v(y)\otimes 1+\sum y' \otimes L_w(y'')\right)\\
&=&\Delta(a)*\Delta(b).
\end{eqnarray*}
(On a utilis\'e le fait que $L_v$ et $L_w$ sont des 1-cocycles pour la premi\`ere  et la cinqui\`eme \'egalit\'e, 
(\ref{eqn30}) pour la troisi\`eme, l'hypoth\`ese de r\'ecurrence appliqu\'ee \`a $x*L_w(y)$ et $L_v(x)*y$ pour la premi\`ere.) 

Donc $(T(V),*,1,\Delta,\varepsilon)$ est une big\`ebre ; d'apr\`es le th\'eor\`eme \ref{theo38}, $T(V)$ a une antipode $S_*$.
Montrons par r\'ecurrence sur $n$ que $S_*(v_1 \top \ldots \top  v_n)=(-1)^n v_n \top \ldots \top v_1$. C'est imm\'ediat pour $n=0$.
Comme $v \in V$ est primitif, on a $S_*(v)=-v$, et donc c'est vrai pour $n=1$.
\\
Si $n \geq 2$ :
\begin{eqnarray*}
m_* \circ (Id \otimes S_*) \circ \Delta(v_1 \top \ldots \top v_n)&=&0\\
&=& (v_1 \top \ldots \top v_n) + S_*(v_1 \top \ldots \top v_n)\\
&&+ \sum_{j=1}^{n-1} (-1)^{n-j} (v_1 \top \ldots \top v_j) * (v_{n}\top \ldots \top v_{j+1})\\ 
&=& (v_1 \top \ldots \top v_n) + S_*(v_1 \top \ldots  \top v_n)\\
&&+ \sum_{j=1}^{n-1} (-1)^{n-j} [(v_1 \top \ldots \top v_{j-1}) * (v_{n}\top \ldots \top v_{j+1})]\top v_j\\
&&+ \sum_{j=1}^{n-1} (-1)^{n-j} [(v_1 \top \ldots \top v_j) * (v_{n}\top \ldots \top v_{j+2})]\top v_{j+1}\\
&=& (v_1 \top \ldots \top v_n) + S_*(v_1 \top \ldots  \top v_n)\\
&&+ \sum_{j=1}^{n-1} (-1)^{n-j} [(v_1 \top \ldots \top v_{j-1}) * (v_{n}\top \ldots \top v_{j+1})]\top v_j\\
&&- \sum_{j=2}^{n} (-1)^{n-j} [(v_1 \top \ldots \top v_{j-1}) * (v_{n}\top \ldots \top v_{j+1})]\top v_{j}\\
&=&( v_1 \top \ldots \top v_n) + S_*(v_1 \top \ldots  \top v_n)\\
&&+ (-1)^{n-1} (v_n \top\ldots \top v_1) -(v_1 \top \ldots \top v_n).
\end{eqnarray*}
(On a utilis\'e (\ref{eqn30}) pour la troisi\`eme \'egalit\'e.)\\
Le r\'esultat est alors imm\'ediat. $\Box$\\

Il est imm\'ediat que $(V^{\otimes n})_{n \geq 0}$ est une graduation de l'alg\`ebre de Hopf $(T(V),*,1,\Delta,\varepsilon,S_*)$.
 
\subsection{Propri\'et\'es de $(T(V),*,1,\Delta,\varepsilon,S_*)$}

\spa On suppose maintenant que $K$ est de caract\'eristique nulle, et que $V$ est muni d'une graduation $(V_i)_{i \geq 0}$ telle que les $V_i$ soient de dimension finie, avec $V_0=(0)$. Cette graduation induit une graduation de  $T(V)$, et il est facile de voir que c'est une graduation d'alg\`ebre de Hopf v\'erifiant $(C_1)$ et $(C_2)$.

\begin{cor}
\label{cor44}
\begin{enumerate}
\item Soit $M=Ker(\varepsilon)$. Soit $G$ un suppl\'ementaire gradu\'e de $M^2$ dans $M$.
Alors $\xi_G :S(G) \longmapsto T(V)$, fixant chaque \'el\'ement de $G$, est un isomorphisme d'alg\`ebres.
\item $\forall x,y \in T(V)$, $deg_p(xy)=deg_p(x)+deg_p(y)$.\\
($deg_p$ est d\'efini dans la section \ref{partB}.)
\end{enumerate}
\end{cor}
{\it Preuve :} 
d\'ecoule imm\'ediatement de la proposition \ref{proL} et du corollaire \ref{corM}. $\Box$\\

Soit $(v_i)_{i \in I}$ une base de $V$ form\'ee d'\'el\'ements homog\`enes. Alors $(v_{i_1} \top \ldots \top v_{i_k})_{k\geq 0,\esp i_j \in I}$ est une base de $T(V)$ form\'ee d'\'el\'ements homog\`enes.

Soit $f_{i_1,\ldots,i_k} \in T(V)^*$, d\'efinie par $f_{i_1,\ldots,i_k} (v_{j_1} \top \ldots \top v_{j_n})=\delta_{(i_1,\ldots, i_k),(j_1,\ldots,j_n)}$.
D'apr\`es la proposition \ref{lemme1}, $(f_{i_1,\ldots, i_k})_{k\geq 0,\esp i_j \in I}$  est une base de $T(V)^{*g}$. 
 \begin{eqnarray*}
f_{i_1,\ldots,i_k}f_{i'_1,\ldots,i'_l}(v_{j_1} \top \ldots \top v_{j_n})&=&
f_{i_1,\ldots,i_k}\otimes f_{i'_1,\ldots,i'_l}\left(\Delta(v_{j_1} \top \ldots \top v_{j_n})\right)\\
&=&f_{i_1,\ldots,i_k}\otimes f_{i'_1,\ldots,i'_l}\left(\sum_{s=0}^{s=n}(v_{j_1} \top \ldots \top v_{j_s})\otimes (v_{j_{s+1}}\top \ldots \top v_{j_n})\right)\\
&=& \delta_{(i_1,\ldots,i_k,i'_1,\ldots,i'_l),(j_1,\ldots,j_n)}.
\end{eqnarray*}
Et donc
$ f_{i_1,\ldots,i_k}f_{i'_1,\ldots,i'_l}=f_{i_1,\ldots,i_k,i'_1,\ldots,i'_l}.$
On en d\'eduit que $T(V)^{*g}$ est isomorphe \`a l'alg\`ebre libre g\'en\'er\'ee par les $f_i$, $i\in I$. De plus :
\begin{eqnarray*}
\Delta(f_i)((v_{i_1}\top \ldots \top v_{i_n}) \otimes  (v_{j_1}\top \ldots \top v_{j_m})) 
&=& f_i((v_{i_1}\top \ldots \top v_{i_n}) * (v_{j_1}\top \ldots \top v_{j_m}))
\end{eqnarray*}
 Si $i_n + j_m >1$, alors $(v_{i_1}\top \ldots \top v_{i_n}) * (v_{j_1}\top \ldots \top v_{j_m}) \in V^{\otimes n}$ pour un certain $n\geq2$, et donc
 $\Delta(f_i)((v_{i_1}\top \ldots \top v_{i_n}) \otimes  (v_{j_1}\top \ldots \top v_{j_m})) =0$.
De plus, $\Delta(f_i)(e_j\otimes 1)=\Delta(f_i)(1 \otimes e_j)=f_i(e_j)=\delta_{i,j}$. Donc :
$$\Delta(f_i)(x \otimes y)=(f_i \otimes 1+ 1 \otimes f_i)(x \otimes y) \esp \forall x,y \in T(V).$$
Donc $f_i$ est primitif. On a d\'emontr\'e :

\begin{theo}
\label{theo58}
Le dual gradu\'e de $(T(V),*,1,\Delta,\varepsilon,S_*)$ est l'alg\`ebre gradu\'ee $T(V^{*g})$ muni du coproduit donn\'e par :
$$\Delta(f)=f\otimes 1 +1 \otimes f \esp \forall f \in   V^{*g}.$$
 \end{theo}

$(T(V)^{*g},*,1,\Delta,\varepsilon,S_*)$  est donc l'alg\`ebre enveloppante de l'alg\`ebre de Lie libre g\'en\'er\'ee par $V^{*g}$ (voir \cite{Bou2}, ch II, $\S 3$).
Par suite, l'alg\`ebre de Lie $Prim(T(V)^{*g})$ est une alg\`ebre de Lie libre.

\subsection{Filtration par $deg_p$}
\begin{lemme}
\label{lem44}
Soit $m$ un produit donnant une structure de big\`ebre sur $T(V)$. Soit $p,q$ des entiers sup\'erieurs \`a $1$.
On a :
$$\tilde{\Delta}^{p+q-1}((v_1\top \ldots \top v_p).(v_{p+1} \top \ldots \top v_{p+q}))
= \sum_{\sigma \in bat(p,q)} v_{\sigma(1)}\otimes \ldots \otimes v_{\sigma(p+q)}.$$
\end{lemme}
{\it Preuve :} d\'ecoule des lemmes \ref{lemI} et \ref{lemme13}. $\Box$

\begin{prop}
Soit $m$ un produit donnant une structure de big\`ebre sur $T(V)$.
L'alg\`ebre de Hopf gradu\'ee associ\'ee \`a la filtration par $deg_p$ est isomorphe \`a $(T(V),*,\eta,\Delta,\varepsilon)$.
\end{prop}
{\it Preuve :} on identifie naturellement $T(V)_{deg_p\leq n}/ T(V)_{deg_p\leq n-1}$ et $V^{\otimes n}$. On a alors un isomorphisme
naturel d'espaces vectoriels  $\Upsilon :gr(T(V)) \longmapsto T(V)$. 
Comme la filtration de la cog\`ebre $T(V)$ provient d'une graduation de cog\`ebre, c'est un isomorphisme de cog\`ebres.\\
Soit $v_1 \top \ldots \top v_p\in V^{\otimes p}$, $v_{p+1} \top \ldots \top v_{p+q}\in V^{\otimes q}$. Par le lemme \ref{lem44},
$$(v_1 \top \ldots \top v_p)(v_{p+1} \top \ldots \top v_{p+q})-(v_1 \top \ldots \top v_p)*(v_{p+1} \top \ldots \top v_{p+q})
\in Ker(\tilde{\Delta}^{p+q-1}).$$
Donc 
$$(v_1 \top \ldots \top v_p)(v_{p+1} \top \ldots \top v_{p+q})-(v_1 \top \ldots \top v_p)*(v_{p+1} \top \ldots \top v_{p+q})
\in (T(V))_{\leq p+q-1}.$$
Comme $(v_1\top \ldots \top v_p)*(v_{p+1}\top \ldots \top v_{p+q}) \in V^{\otimes p+q}$, 
$\Upsilon$ est un isomorphisme d'alg\`ebres. $\Box$

\subsection{Produits commutatifs sur $T(V)$}
\begin{theo}
\label{theo60}
 On suppose $K$ de caract\'eristique nulle,  et $V$ de dimension finie ou infinie d\'enombrable.
Soit $m$ un produit commutatif sur $T(V)$, tel que $(T(V),m,\eta,\Delta,\varepsilon)$
soit une big\`ebre (et donc une alg\`ebre de Hopf par le th\'eor\`eme \ref{theo38}).
Alors $(T(V),m,\eta,\Delta,\varepsilon,S)$ et $(T(V),*,1,\Delta,\varepsilon,S_*)$
sont isomorphes comme alg\`ebres de Hopf.
\end{theo}
{\it Preuve :} soit $M$ l'id\'eal d'augmentation de $T(V)$, et $M^{*2}$ son carr\'e
dans $(T(V),*,1)$. On note $M^{*2}_i=M^{*2} \cap V^{\otimes i}$. On a $M^{*2}=\bigoplus  M^{*2}_i$.
Soit $G_i$ sous-espace de $T(V)$, tel que $V^{\otimes i}=M^{*2}_i \oplus G_i$ ($i \geq 1$).
On pose $G=\bigoplus G_i$, de sorte que $M=G \oplus M^{*2}$.\\
Soit $(v_n)_{n \in I}$ une base de $V$. Par le th\'eor\`eme de la base incompl\`ete, 
on peut choisir $G_i$ ayant une base de la forme $(v_{j_1}\top \ldots \top v_{j_i})_{(j_1,\ldots,j_i) \in K_i}$ 
avec $K_i$ une partie de $I^i$.\\

On construit par r\'ecurrence $\phi_n :T(V)_{deg_p\leq n} \longmapsto T(V)_{deg_p\leq n}$ v\'erifiant :
\begin{enumerate}
\item ${\phi_n}_{\mid T(V)_{deg_p\leq n-1}}=\phi_{n-1}$ ;
\item $\phi_n$ morphisme de cog\`ebres ;
\item $\forall x,y \in T(V)$, tel que $deg_p(x)+deg_p(y) \leq n$, 
$\phi_n(x*y)=\phi_n(x) \phi_n(y)$.
\end{enumerate}

On prend $\phi_0=Id_K$. Les conditions 2 et 3 sont v\'erifi\'ees  pour $n=0$.
On prend $\phi_1=Id_{K\oplus V}$. Les conditions 1,2 et 3 sont v\'erifi\'ees pour $n=1$.
Supposons $\phi_n$ construit, et construisons $\phi_{n+1}$.
D'apr\`es la condition 1, il reste \`a d\'efinir $\phi_{n+1}$ sur $V^{\otimes n+1}=G_{n+1} \oplus M^{*2}_{n+1}$.

Comme $V$ est de dimension au plus infinie d\'enombrable, $V$ peut \^etre muni d'une graduation $(V_i)_{i\in \mathbb{N}}$ telle que $V_0=(0)$, et $dim(V_i)$ soit finie pour tout $i$. En choisissant la base $(v_n)$ form\'ee d'\'el\'ements homog\`enes, $G$ est alors un sous-espace gradu\'e de $T(V)$. D'apr\`es la proposition \ref{cor44}, $\xi_G :S(G) \longmapsto (T(V),*,1)$ est un isomorphisme d'alg\`ebres. 
On peut donc d\'efinir $\phi_{n+1}$ sur $M^{*2}_{n+1}$  par :
\begin{eqnarray}
\label{eqn50}
\phi_{n+1}(x_1 *\ldots * x_k)&=&\prod_{i=1}^{k} \phi_n(x_i)
\end{eqnarray}
pour $x_i,\ldots, x_k \in G_1 \oplus \ldots \oplus G_n$, $\sum deg_p(x_i)=n+1$.

Il reste \`a d\'efinir $\phi_{n+1}(v_{j_1} \top \ldots \top v_{j_{n+1}})$ pour $(j_1,\ldots,j_{n+1}) \in K_{n+1}$.\\

Soient $w_1,\ldots,w_k$ dans $V$, $k\leq n$.
On pose $w_1 \mtop \ldots \mtop w_k=\phi_n(w_1 \top \ldots \top w_k)$.
Comme $\phi_n$ est un morphisme de cog\`ebres, et que $\phi_n(1)=1$, on a
$\tilde{\Delta} \circ \phi_n=(\phi_n \otimes \phi_n) \circ \tilde{\Delta}$, 
et par it\'eration, $\tilde{\Delta}^{k-1} \circ \phi_n=\phi_n^{\otimes k} \circ \tilde{\Delta}$, si $k\leq n$.
Donc :
\begin{eqnarray}
\nonumber \tilde{\Delta}(w_1 \mtop \ldots \mtop w_k)&=&
\phi_n(w_1) \otimes \ldots \otimes \phi_n(w_k)\\
\label{eqn51} &=&w_1 \otimes \ldots \otimes w_k.
\end{eqnarray}
(Car ${\phi_k}_{\mid V}=Id_V$ si $k\geq 1$.)

Montrons par r\'ecurrence que $\phi_k$ est bijectif si $k\leq n$.
C'est imm\'ediat pour $k=0,1$. Soit $k \geq 2$, et supposons $\phi_{k-1}$ bijectif.

Surjectivit\'e : il faut montrer que $w_1 \top \ldots \top w_k \in Im(\phi_k)$, $\forall w_1,\ldots,w_k \in V$.
D'apr\`es (\ref{eqn51}),  $\tilde{\Delta}^{k-1}((w_1 \top \ldots \top w_k)-(w_1 \mtop \ldots \mtop w_k))=0$.
On en d\'eduit 
 $deg_p((w_1 \top \ldots \top w_k)-(w_1 \mtop \ldots \mtop w_k))<k$, 
et donc : $$ (w_1 \top \ldots \top w_k)-(w_1 \mtop \ldots \mtop w_k) \in Im(\phi_{k-1}) \subset Im(\phi_k).$$
Comme par d\'efinition, $w_1 \mtop \ldots \mtop w_k$ est dans $Im(\phi_k)$, on a le r\'esultat.

Injectivit\'e : soit $x=x_0+\ldots+x_l$, $x_i \in V^{\otimes i}$, $l\leq k$, $x_l \neq 0$,
tel que $\phi_k(x)=0$. Si $l=0$ ou $1$, alors $\phi_k(x)=x$ : contradiction.
Donc $l\geq 2$.
De plus, $\varepsilon \circ \phi_k(x)=0=\varepsilon(x)=x_0$ : $x_0=0$. D'apr\`es le lemme \ref{lemme13} :
$$  \tilde{\Delta}^{l-1}(x)=\tilde{\Delta}^{l-1}(x_l) \neq 0.$$
De plus, $\tilde{\Delta}^{l-1}(x_l) \in V^{\otimes l}$, donc : 
$$ \phi_k^{\otimes l} \circ \tilde{\Delta}^{l-1}(x)=\tilde{\Delta}^{l-1}(x_l)\neq 0.$$
Or 
$$ \phi_k^{\otimes l} \circ \tilde{\Delta}^{l-1}(x)=\tilde{\Delta}^{l-1}(\phi_k(x))=0.$$
On aboutit \`a une contradiction, et donc $\phi_k$ est injectif $\forall k \leq n$.\\

On peut donc d\'efinir, pour $v \in V$ :
\begin{eqnarray*}
L_v :T(V)_{deg_p\leq n-1} & \longmapsto & T(V)\\
\phi_n(v_1 \top \ldots \top v_k) & \longmapsto & \phi_n(v_1 \top \ldots \top v_k \top v)\\
=v_1 \mtop \ldots \mtop v_k &  & =v_1 \mtop \ldots \mtop v_k \mtop v.\\
(k \leq n-1)
\end{eqnarray*}
Comme $\phi_n$ est un morphisme de cog\`ebres, on a :
$$\forall x \in T(V)_{deg_p\leq n-1},\esp \Delta(L_v(x))=L_v(x) \otimes 1 + (Id \otimes L_v) \circ \Delta(x).$$
D'apr\`es le corollaire \ref{cor41}, il existe $\overline{L}_v$ un 1-cocycle de $T(V)$, tel que
${\overline{L}_v}_{\mid T(V)_{\leq n-1}}=L_v.$\\

On pose alors : 
$$\phi_{n+1}(v_{j_1} \top \ldots \top v_{j_{n+1}})=\overline{L}_{v_{j_{n+1}}}
(v_{j_1} \mtop \ldots \mtop v_{j_n}), \esp \forall (j_1, \ldots , j_{n+1}) \in K_{n+1}.$$
Alors $\phi_{n+1}$ v\'erifie 1 par construction, et 3 d'apr\`es (\ref{eqn50}). 
Reste \`a montrer 2. Soit $x \in T(V)_{deg_p\leq n+1}$. On veut monter que 
$\Delta \circ \phi_{n+1}(x)=(\phi_{n+1} \otimes \phi_{n+1}) \circ \Delta (x)$. 
Si $x \in T(V)_{deg_p\leq n}$,
on a :
\begin{eqnarray*}
\Delta\circ \phi_{n+1}(x)&=&\Delta \circ \phi_{n}(x)\\
&=& (\phi_n \otimes \phi_n) \circ \Delta (x)\\
&=& (\phi_{n+1} \otimes \phi_{n+1}) \circ \Delta (x).
\end{eqnarray*}

Si $x \in M^{*2}_{n+1}$, on peut se ramener \`a $x=x_1*x_2$, $x_1$ et $x_2$ dans $M$.
D'apr\`es la  proposition \ref{cor44}-2, $deg_p(x)=deg_p(x_1)+deg_p(x_2)$, et donc $deg_p(x_i)\leq n$ pour $i=1,2$.
Donc on a le r\'esultat pour $x_1$ et $x_2$ ; en utilisant (\ref{eqn50}) et le fait que $\Delta$
soit un morphisme d'alg\`ebres, on a le r\'esultat pour $x$.

Si $x \notin M^{*2}_{n+1}$, on peut se ram\`ener \`a $x=v_{j_1} \top \ldots \top v_{j_{n+1}}$, $(j_1,\ldots,j_{n+1})\in K_{n+1}$.
\begin{eqnarray*}
\Delta\circ \phi_{n+1}(x)&=&
\Delta\circ \overline{L}_{v_{j_{n+1}}}(v_{j_1} \mtop \ldots \mtop v_{j_n})\\
&=& x \otimes 1 +1 \otimes x\\
&&+\sum_{k=1}^{n} (v_{j_1} \mtop \ldots \mtop v_{j_k})
\otimes \overline{L}_{v_{j_{n+1}}}(v_{j_{k+1}} \mtop \ldots \mtop v_{j_n})\\ 
&=& x \otimes 1 +1 \otimes x\\
&&+\sum_{k=1}^{n} (v_{j_1} \mtop \ldots \mtop v_{j_k})
\otimes L_{v_{j_{n+1}}}(v_{j_{k+1}} \mtop \ldots \mtop v_{j_n})\\ 
&=& x \otimes 1 +1 \otimes x\\
&&+\sum_{k=1}^{n} \phi_n(v_{j_1} \top \ldots \top v_{j_k})
\otimes \phi_n(v_{j_{k+1}} \top \ldots \top v_{j_{n+1}})\\ 
&=& (\phi_{n+1}\otimes \phi_{n+1})\circ \Delta(x).
\end{eqnarray*}
(On a utilis\'e le fait que $\overline{L}_{v_{j_{n+1}}}$ soit un 1-cocycle pour la deuxi\`eme \'egalit\'e.)\\

On d\'efinit alors $\phi :T(V) \longmapsto T(V)$ par $\phi_{\mid T(V)_{deg_p\leq n}} =\phi_n$ ;  
par la condition 1, $\phi$ est bien d\'efini. Par la condition 2, $\phi$ est un morphisme de cog\`ebres ;
par la condition 3, $\phi$ est un morphisme d'alg\`ebres. De plus, comme les $\phi_n$ sont bijectifs,
$\phi$ est un isomorphisme d'alg\`ebres de Hopf. $\Box$

\section{Applications des r\'esultats aux alg\`ebres de Hopf $\nH$}
\subsection{Lien avec $\Hr$}
\spa On a une surjection $\varpi :\trees \longmapsto \ntrees$ qui \`a un arbre enracin\'e plan d\'ecor\'e associe le m\^eme graphe, muni des m\^emes d\'ecorations, en oubliant la donn\'ee du plongement dans le plan.
On prolonge $\varpi$ de $\forets$ vers $\nforets$ en posant $\varpi(t_1 \ldots t_n)=\varpi(t_1) \ldots \varpi(t_n)$.

\begin{prop} 
On consid\`ere l'application suivante :
$$\Phi :\left\{
\begin{array}{rcl}
\Hr & \longmapsto & \nH\\
F \in \forets & \longmapsto &  \varpi(F).
\end{array}\right. $$
Alors $\Phi$ est un morphisme surjectif d'alg\`ebres de Hopf gradu\'ees. 
\end{prop}
{\it Preuve :}
on note ${\cal B}_d^+$ l'op\'erateur de $\nH$ qui envoie un \'el\'ement $t_1 \ldots t_n \in \nforets$ \`a l'arbre obtenu en greffant les racines de $t_1, \ldots, t_n$ \`a une racine commune d\'ecor\'ee par $d$.
D'apr\`es \cite{Connes}, il s'agit d'un 1-cocycle de $\nH$.

Par la propri\'et\'e universelle de $\Hr$, il existe un unique morphisme d'alg\`ebres de Hopf  $\Phi' :\Hr \longmapsto \nH$ tel que $\Phi' \circ B_d^+={\cal B}_d^+ \circ \Phi'.$
Montrons par r\'ecurrence sur $n=poids(F)$ que $\Phi'(F)=\varpi(F)$. C'est vrai si $n=0$. Suppposons ce r\'esultat vrai pour toute for\^et de poids inf\'erieur ou \'egal \`a $n$, et soit $F$ de poids $n+1$.
Si $F\notin \trees$, il existe $F_1,F_2 \in \forets$, non vides, telles que $F=F_1F_2$. Alors $\Phi'(F)=\Phi'(F_1)\Phi'(F_2)=\varpi(F_1)\varpi(F_2)=\varpi(F)$.
Sinon, il existe $d \in \cal D$, $F_1 \in \forets,$ tels que $F=B_d^+(F_1)$ ; alors $\Phi'(F)={\cal B}_d^+(\varpi(F_1))=\varpi(F)$ ; par suite, $\Phi=\Phi'$ est un morphisme d'alg\`ebres de Hopf. Enfin, il est imm\'ediat que $\Phi$ est homog\`ene de degr\'e 0. $\Box$

\subsection{Primitifs de $\nH$}

\spa On utilise les notations des parties \ref{part10.2} et \ref{part10.3}. 

Pour $p_1, \ldots , p_n \in Prim(\Hr)$, on d\'efinit par r\'ecurrence :
$$ p_1 \mtop \ldots \mtop p_n = (p_1 \mtop \ldots \mtop p_{n-1})\mtop p_n.$$

\begin{lemme}
\label{lemme81}
Soient $x, y\in \Hr$. Alors $\Phi(x \mtop y)= \Phi(x) \otop \Phi(y)$.
\end{lemme}
{\it Preuve :} c'est imm\'ediat si $x,y \in \forets$. $\Box$
\\

A l'aide de (\ref{eq29}), et en remarquant que $I$ est un coid\'eal inclus dans $Ker(\Phi)$, on montre par une r\'ecurrence simple le lemme suivant :

\begin{lemme}
\label{lemme82}
Soient $p_1,\ldots, p_n \in Prim(\Hr)$. On a :
$$\tdelta^{n-1}(p_1 \mtop \ldots \mtop p_n) = p_1 \otimes \ldots \otimes p_n + \sum_{i=1}^{n} \Hr \otimes \ldots \otimes \underbrace{Ker(\Phi)}_{\mbox{$i^{\mbox{\`eme}}$ position}} \otimes \ldots \otimes \Hr.$$
\end{lemme}

On pose donc :
\begin{eqnarray*}
\tilde{F}_i : Prim(\Hr)^{\otimes i} &\longmapsto & \Hr\\
p_1 \otimes \ldots \otimes p_i &\longmapsto &p_1 \mtop \ldots \mtop p_i.
\end{eqnarray*}
On a $\Hr =Ker(\Phi)+ \sum Im(\tilde{F}_i)$.  De plus, d'apr\`es le lemme \ref{lemme81}, 
\begin{eqnarray}
\nonumber \tilde{F}_i(p_1 \otimes \ldots \otimes p_i)&=&\Phi(p_1) \otop \ldots \otop \Phi(p_i)\\
\label{eqn33} &=& \overline{F}_i (\Phi(p_1) \otimes \ldots \otimes \Phi(p_i)).
\end{eqnarray}

\begin{theo}
$\Phi :Prim(\Hr)\longmapsto Prim(\nH)$ est surjectif.
\end{theo}
{\it Preuve} : comme $\Phi$ est un morphisme d'alg\`ebres de Hopf, $\Phi(Prim(\Hr))\subseteq Prim(\nH)$. 
Soit $p \in Prim(\nH)$. Comme $\Phi$ est surjectif, il existe $x \in \Hr$, $\Phi(x)=p$. On peut supposer que $x = x_0 1+\tilde{F}_1(y_1)+ \ldots +\tilde{F}_n(y_n)$, $x_0 \in \mathbb{Q}$,
$y_i \in Prim(\Hr)^{\otimes i}$. $\varepsilon(x)=\varepsilon(\Phi(x))=\varepsilon(p)=0$, donc $x_0=0$. Supposons $n >1$. On a alors :
\begin{eqnarray*}
\tdelta^{n-1}(\Phi(x))&=&0\\
&=&\Phi^{\otimes n} \circ \tdelta^{n-1}(x)\\
&=& \Phi^{\otimes n} \circ \tdelta^{n-1}(\tilde{F}_n(y_n))+0\\
&=&\Phi^{\otimes n} (y_n)+0.
\end{eqnarray*}
(on a utilis\'e le lemme \ref{lemme82} pour la troisi\`eme et la quatri\`eme \'egalit\'e).\\

D'apr\`es (\ref{eqn33}), on a alors :
\begin{eqnarray*}
\Phi (\tilde{F}_n(y_n))&=& \overline{F}_n ( \Phi^{\otimes n}(y_n))\\
&=&0.
\end{eqnarray*}
Donc $p=\Phi(x)=\Phi(\tilde{F}_1(y_1))+\ldots +\Phi(\tilde{F}_{n-1}(y_{n-1}))$. Par une r\'ecurrence descendante, on d\'emontre que
$p=\Phi(x)=\Phi(\tilde{F}_1(y_1))$. Or $\tilde{F}_1(y_1) \in Prim(\Hr)$, et donc $\Phi(prim(\Hr))=Prim(\nH)$. $\Box$

\subsection{Dual gradu\'e de $\nH$}
On va chercher $(Ker(\Phi))^{\perp}$ dans la dualit\'e entre $\Hr$ et elle-m\^eme de la section \ref{part55}. 
\begin{prop}
\label{pro63}
\begin{enumerate}
\item $(Ker(\Phi))^{\perp}$ est une sous-alg\`ebre de Hopf cocommutative de $\Hr$, stable par $\gamma_d$, $\forall d \in \cal D$.
\item Pour $\overline{F} \in \nforets$, on pose : $$e_{\overline{F}}=\sum_{\varpi(F')=\overline{F}} \esp e_{F'}.$$
Alors $(e_{\overline{F}})_{\overline{F} \in \nforets}$ est une base de $(Ker(\Phi))^{\perp}$.
\item Pour $\overline{F}_1,\overline{F}_2,\overline{F} \in \nforets$, soit $n(\overline{F}_1,\overline{F}_2 ;\overline{F})$ le  nombre de
coupes admissibles de $\overline{F}$ telles que $P^c(\overline{F})=\overline{F}_1$ et $R^c(\overline{F})=\overline{F}_2$. Alors :
$$e_{\overline{F}_1} e_{\overline{F}_2}=\sum_{\overline{F} \in \nforets} n(\overline{F}_1,\overline{F}_2 ;\overline{F}) e_{\overline{F}}.$$
\end{enumerate}
\end{prop}
{\it Preuve :}

1. Comme $\Phi$ est un morphisme d'alg\`ebres de Hopf, $Ker(\Phi)$ est un id\'eal bilat\`ere, un coid\'eal, et est stable par $S$.
De plus, $\Phi \circ B_d^+={\cal B}_d^+ \circ \Phi$, donc $Ker(\Phi)$ est stable par $B_d^+$. 

$1 \in (Ker(\Phi))^{\perp}$ car $Ker(\Phi) \subset Ker(\varepsilon)$.

Soient $x,y \in (Ker(\Phi))^{\perp}$, $z \in Ker(\Phi)$ ; $(xy, z) =(x \otimes y, \Delta(z))$ ; or $\Delta(z) \in
Ker(\Phi)\otimes \Hr + \Hr \otimes Ker(\Phi)$ ; donc $(xy,z)=0$ : $xy \in (Ker(\Phi))^{\perp}$.

Soit $x \in (Ker(\Phi))^{\perp}$, $y \otimes z \in \Hr \otimes Ker(\Phi) + Ker(\Phi)\otimes \Hr$.
Alors $(\Delta(x),y\otimes z)=(x,yz)=0$ car $Ker(\Phi)$ est un id\'eal bilat\`ere ; par suite, $\Delta(x) \in (\Hr \otimes Ker(\Phi) + Ker(\Phi)\otimes \Hr)^{\perp}=(Ker(\Phi))^{\perp} \otimes(Ker(\Phi))^{\perp}$.

Comme $Ker(\Phi)$ est stable par $S$, $(Ker(\Phi))^{\perp}$ est stable par $S^{*g}=S$.

Comme $Ker(\Phi)$ est stable par $B_d^+$, $(Ker(\Phi))^{\perp}$ est stable par $(B_d^+)^{*g}=\gamma_d$.

Enfin, soient $x \in (Ker(\Phi))^{\perp}$, $y,z \in \Hr$.
\begin{eqnarray*}
(\Delta(x)-\Delta^{op}(x), y\otimes z)&=&(\Delta(x),y\otimes z -z \otimes y)\\
&=&(x, yz-zy).
\end{eqnarray*}
Comme $\nH$ est commutative, $yz-zy \in Ker(\Phi)$, et donc $(x,yz-zy)$ est nul ; par suite, $\Delta(x)=\Delta^{op}(x)$.\\

2. La famille $(G-G')_{G,G' \in \forets, \varpi(G)=\varpi(G')}$ g\'en\`ere lin\'eairement $Ker(\Phi)$. On montre  facilement que
$(e_{\overline{F}}, G-G')=0$ si $\varpi(G)=\varpi(G')$. Donc $e_{\overline{F}} \in (Ker(\Phi))^{\perp}$. Il est imm\'ediat qu'il s'agit d'une famille libre ; en comparant les dimensions des composantes homog\`enes
de $(Ker(\Phi))^{\perp}$ et de l'espace engendr\'e par les $e_{\overline{F}}$, on montre qu'il s'agit d'une base de $(Ker(\Phi))^{\perp}$.

3. La dualit\'e entre $\Hr$ et elle m\^eme engendre une dualit\'e entre $(Ker(\Phi))^{\perp}$ et $\Hr/Ker(\Phi) \approx \nH$. On a alors, pour $F,G \in \forets$ :
\begin{eqnarray*}
(e_{\varpi(F)},\varpi(G))&=&(\sum_{\varpi(F')=\varpi(F)} e_{F'}, G\esp )\\
&=&\delta_{\varpi(F), \varpi(G)}.
\end{eqnarray*}
Alors pour $\overline{F}_1,\overline{F}_2, \overline{F} \in \nforets :$
\begin{eqnarray*}
(e_{\overline{F}_1} e_{\overline{F}_2}, \overline{F})&=&(e_{\overline{F}_1} \otimes e_{\overline{F}_2}, \Delta(\overline{F}))\\
&=& n(\overline{F}_1,\overline{F}_2 ;\overline{F}).
\end{eqnarray*}
D'o\`u le r\'esultat annonc\'e. $\Box$

\begin{cor}
Soit ${\cal L}_1^{\cal D}$ l'alg\`ebre de Lie de base $(e_{\overline{t}})_{\overline{t} \in \ntrees}$, avec :
$$[e_{\overline{t}_1} ;e_{\overline{t}_2}]= \sum_{\overline{t}\in \ntrees} \left(n(\overline{t}_1,\overline{t}_2 ;\overline{t})
- n(\overline{t}_2,\overline{t}_1 ;\overline{t})\right) e_{\overline{t}}.$$
Alors le dual gradu\'e de $\nH$ est isomorphe \`a ${\cal U}({\cal L}_1^{\cal D})$ comme alg\`ebre de Hopf gradu\'ee.
\end{cor}
{\it Preuve :} $(\nH)^{*g} \approx (Ker(\Phi))^{\perp}$. Comme $(Ker(\Phi))^{\perp}$ est cocommutative, d'apr\`es le th\'eor\`eme \ref{theoK},
$(Ker(\Phi))^{\perp} \approx {\cal U}(Prim(Ker(\Phi))^{\perp} )$. De plus, une base de 
$Prim((Ker(\Phi)^{\perp}) $ est $(e_{\overline{t}})_{\overline{t} \in \ntrees}$ d'apr\`es le th\'eor\`eme \ref{theo28}-7.
La formule pour $ [e_{\overline{t}_1} ;e_{\overline{t}_2}]$ se d\'eduit 
de la proposition \ref{pro63}-3. $\Box$ \\

{\it Remarque :} on a retrouv\'e ainsi le r\'esultat de \cite{Connes,Panaite}.

\subsection{Applications du th\'eor\`eme 80}

\begin{theo}
Soit $\cal D$ un ensemble non vide, fini ou d\'enombrable. Soit $V$ un espace de dimension infinie d\'enombrable.
Alors les alg\`ebres de Hopf $\nH$ et $(T(V),*,1, \Delta, \varepsilon,S_*)$ sont isomorphes.
\end{theo}
{\it Preuve :} on a vu dans la partie \ref{part10.2} que $\nH$ et $T(Prim(\nH))$ sont des cog\`ebres isomorphes, $Prim(\nH)$ \'etant de dimension infinie d\'enombrable (donc isomorphe \`a $V$).
Comme $\nH$ est commutative, d'apr\`es le th\'eor\`eme \ref{theo60}, $\nH$ est isomorphe \`a l'alg\`ebre de Hopf $(T(V),*,1,\Delta,\varepsilon,S_*)$. $\Box$\\

{\it Remarque :} si $\cal D$ et $\cal D'$ sont deux ensembles non vides, finis ou d\'enombrables, ${\cal H}_R^{\cal D}$  et ${\cal H}_R^{\cal D'}$ sont donc isomorphes comme  alg\`ebres de Hopf.
Elles sont isomorphes comme alg\`ebres de Hopf gradu\'ees si et seulement si $\cal D$ et $\cal D'$ ont le m\^eme cardinal. 

\begin{cor}
Les alg\`ebres de Lie ${\cal L}_1^{\cal D}$ sont des alg\`ebres de Lie libres.
\end{cor}
{\it Preuve :} application directe de la proposition \ref{theo58}. $\Box$

\begin{prop}
Les alg\`ebres de Lie $Prim(\Hr)$ sont des alg\`ebres de Lie libres.
\end{prop}
{\it Preuve :} d'apr\`es la partie \ref{part10.3}, $(\Hr)_{ab}$ est une cog\`ebre tensorielle, ainsi qu'une alg\`ebre commutative. D'apr\`es le th\'eor\`eme \ref{theo60} et la proposition
\ref{theo58}, l'alg\`ebre de Lie $Prim((\Hr)_{ab}^{*g})$  est libre. Or $(\Hr)^{*g}_{ab} $ s'identifie naturellement \`a l'orthogonal de l'id\'eal bilat\`ere $I$ engendr\'e par les $xy-yx$.
Donc $Prim((\Hr)_{ab}^{*g})\approx Prim(I^{\perp})$, qui est une sous-alg\`ebre de Lie de $Prim(\Hr)$. Mais $I \subseteq M^2$, donc $Prim(\Hr) \subset  (1) + Prim( \Hr)  = (M^2)^{\perp} \subset I^{\perp}$. Donc $Prim(I^{\perp})=Prim(\Hr)$. $\Box$ 

\section{Appendice}

\subsection{Valeurs des $\tau_k$}
\label{part81}

$\begin{array}{|rcl|rcl|}
\hline
\tau_{1}&=&1&\tau_{13}&=&208012\\
\tau_{2}&=&1&\tau_{14}&=&742900\\
\tau_{3}&=&2&\tau_{15}&=&2674440\\
\tau_{4}&=&5&\tau_{16}&=&9694845\\
\tau_{5}&=&14&\tau_{17}&=&35357670\\
\tau_{6}&=&42&\tau_{18}&=&129644790\\
\tau_{7}&=&132&\tau_{19}&=&477638700\\
\tau_{8}&=&429&\tau_{20}&=&1767263190\\
\tau_{9}&=&1430&\tau_{21}&=&6564120420\\
\tau_{10}&=&4862&\tau_{22}&=&24466267020\\
\tau_{11}&=&16796&\tau_{23}&=&91482563640\\
\tau_{12}&=&58786&\tau_{24}&=&343059613650\\
\hline
\end{array}
$
 
\subsection{Forme bilin\'eaire $(\hspace{1mm},)$ dans ${\cal H}_{P,R}$}
\spa On se place dans ${\cal H}_{P,R}$.
On ordonne les for\^ets de poids $n$ par $\geq$ :
$F_1 \leq F_2 \leq \ldots \leq F_{r_n}$ 
 ${\cal B}'_n=(F_i)_{i\leq r_n}$ est une base de  ${\cal H}_n$ (voir la figure 12).

\begin{figure}[h]
\framebox(450,180){
\begin{picture}(0,0)(170,-50)
\large{${\cal B}'_2=$}\huge{(}
\begin{picture}(17,20)(-3,0)
\put(0,0){\circle*{5}}
\put(10,0){\circle*{5}}
\end{picture}
\large{,}
\begin{picture}(7,10)(-3,0)
\put(0,0){\circle*{5}}
\put(0,10){\circle*{5}}
\put(0,0){\line(0,1){10}}
\end{picture}
\huge{)}\large{ ; $\esp {\cal B}'_3=$}\huge{(}
\begin{picture}(25,10)(0,0)
\put(0,0){\circle*{5}}
\put(10,0){\circle*{5}}
\put(20,0){\circle*{5}}
\end{picture}
\large{,}
\begin{picture}(20,10)(-3,0)
\put(0,0){\circle*{5}}
\put(10,0){\circle*{5}}
\put(0,10){\circle*{5}}
\put(0,0){\line(0,1){10}}
\end{picture}
\large{,}
\begin{picture}(20,10)(-3,0)
\put(0,0){\circle*{5}}
\put(10,0){\circle*{5}}
\put(10,10){\circle*{5}}
\put(10,0){\line(0,1){10}}
\end{picture}
\large{,}
\begin{picture}(25,15)(-3,0)
\put(10,0){\circle*{5}}
\put(10,0){\line(-1,1){10}}
\put(10,0){\line(1,1){10}}
\put(0,10){\circle*{5}}
\put(20,10){\circle*{5}}
\end{picture}
\large{,}
\begin{picture}(7,30)(-3,0)
\put(0,0){\circle*{5}}
\put(0,10){\circle*{5}}
\put(0,0){\line(0,1){10}}
\put(0,10){\line(0,1){10}}
\put(0,20){\circle*{5}}
\end{picture}
\huge{)}\large{ ;}
\begin{picture}(0,0)(340,50)
\large{${\cal B}'_4=$}\huge{(}
\begin{picture}(35,10)(-3,0)
\put(0,0){\circle*{5}}
\put(10,0){\circle*{5}}
\put(20,0){\circle*{5}}
\put(30,0){\circle*{5}}
\end{picture}
\large{,}
\begin{picture}(30,10)(-3,0)
\put(0,0){\circle*{5}}
\put(10,0){\circle*{5}}
\put(0,10){\circle*{5}}
\put(00,0){\line(0,1){10}}
\put(20,0){\circle*{5}}
\end{picture}
\large{,}
\begin{picture}(30,10)(-3,0)
\put(0,0){\circle*{5}}
\put(10,0){\circle*{5}}
\put(10,10){\circle*{5}}
\put(10,0){\line(0,1){10}}
\put(20,0){\circle*{5}}
\end{picture}
\large{,}
\begin{picture}(35,15)(-3,0)
\put(10,0){\circle*{5}}
\put(10,0){\line(-1,1){10}}
\put(10,0){\line(1,1){10}}
\put(0,10){\circle*{5}}
\put(20,10){\circle*{5}}
\put(30,0){\circle*{5}}
\end{picture}
\large{,}
\begin{picture}(20,30)(-3,0)
\put(0,0){\circle*{5}}
\put(0,10){\circle*{5}}
\put(0,0){\line(0,1){10}}
\put(0,10){\line(0,1){10}}
\put(0,20){\circle*{5}}
\put(10,0){\circle*{5}}
\end{picture}
\large{,}
\begin{picture}(30,10)(-3,0)
\put(0,0){\circle*{5}}
\put(10,0){\circle*{5}}
\put(20,10){\circle*{5}}
\put(20,0){\line(0,1){10}}
\put(20,0){\circle*{5}}
\end{picture}
\large{,}
\begin{picture}(20,10)(-3,0)
\put(0,0){\circle*{5}}
\put(0,10){\circle*{5}}
\put(0,0){\line(0,1){10}}
\put(10,0){\circle*{5}}
\put(10,10){\circle*{5}}
\put(10,0){\line(0,1){10}}
\end{picture}
\large{,}
\end{picture}
\end{picture}
\begin{picture}(0,0)(100,60)
\begin{picture}(35,15)(-3,0)
\put(20,0){\circle*{5}}
\put(20,0){\line(-1,1){10}}
\put(20,0){\line(1,1){10}}
\put(10,10){\circle*{5}}
\put(30,10){\circle*{5}}
\put(0,0){\circle*{5}}
\end{picture}
\large{,}
\begin{picture}(20,30)(-3,0)
\put(10,0){\circle*{5}}
\put(10,10){\circle*{5}}
\put(10,0){\line(0,1){10}}
\put(10,10){\line(0,1){10}}
\put(10,20){\circle*{5}}
\put(0,0){\circle*{5}}
\end{picture}
\large{,}
\begin{picture}(22,10)(0,0)
\put(0,10){\circle*{5}}
\put(10,10){\circle*{5}}
\put(20,10){\circle*{5}}
\put(10,0){\circle*{5}}
\put(10,0){\line(-1,1){10}}
\put(10,0){\line(1,1){10}}
\put(10,0){\line(0,1){10}}
\end{picture}
\large{,}
\begin{picture}(22,10)(-3,0)
\put(0,10){\circle*{5}}
\put(20,10){\circle*{5}}
\put(0,20){\circle*{5}}
\put(0,10){\line(0,1){10}}
\put(10,0){\circle*{5}}
\put(10,0){\line(-1,1){10}}
\put(10,0){\line(1,1){10}}
\end{picture}
\large{,}
\begin{picture}(25,10)(-3,0)
\put(0,10){\circle*{5}}
\put(20,10){\circle*{5}}
\put(20,20){\circle*{5}}
\put(20,10){\line(0,1){10}}
\put(10,0){\circle*{5}}
\put(10,0){\line(-1,1){10}}
\put(10,0){\line(1,1){10}}
\end{picture}
\large{,}
\begin{picture}(25,15)(-3,0)
\put(10,10){\circle*{5}}
\put(10,10){\line(-1,1){10}}
\put(10,10){\line(1,1){10}}
\put(0,20){\circle*{5}}
\put(20,20){\circle*{5}}
\put(10,0){\circle*{5}}
\put(10,0){\line(0,1){10}}
\end{picture}
\large{,}
\begin{picture}(7,30)(-3,0)
\put(0,10){\circle*{5}}
\put(0,20){\circle*{5}}
\put(0,10){\line(0,1){10}}
\put(0,20){\line(0,1){10}}
\put(0,30){\circle*{5}}
\put(0,0){\circle*{5}}
\put(0,0){\line(0,1){10}}
\end{picture}
\huge{)}\large{.}
\end{picture}
}
\caption{{\it les bases ${\cal B}'_n$ pour $n \leq 4$.}}
\end{figure}
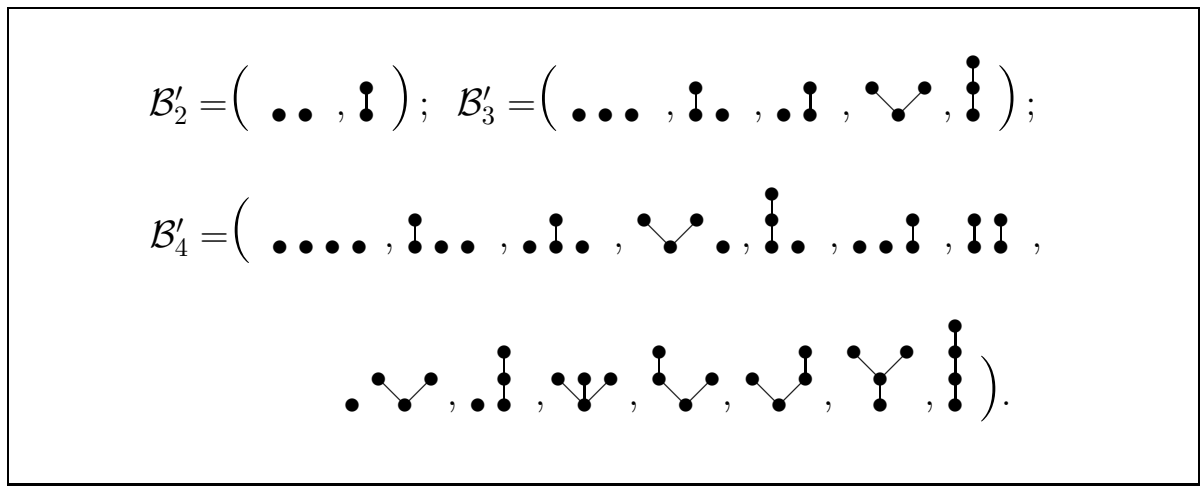

$A'_n$ est la matrice de la forme bilin\'eaire $(\hspace{1mm},)$ restreinte \`a ${\cal H}_n \times {\cal H}_n$ dans la base ${\cal B}'_n$. 
On obtient :\begin{eqnarray*}
 A'_1=\left[ \begin{array}{c}
           1
\end{array} \right] ; \hspace{8mm}
 A'_2=\left[ \begin{array}{cc}
2&1\\
1&0
\end{array} \right] ;&&
 A'_3=\left[ \begin{array}{ccccc}
6&3&3&2&1\\
3&1&1&1&0\\
3&1&1&0&0\\
2&1&0&0&0\\
1&0&0&0&0
\end{array} \right] ;
\end{eqnarray*}
$$ A'_4=\left[ \begin{array}{cccccccccccccc}
24&12&12&8&4&12&6&8&4&6&3&3&2&1\\
12&5&5&4&1&5&2&4&1&3&1&1&1&0\\
12&5&5&3&1&5&2&3&1&3&1&1&0&0\\
8&4&3&2&1&2&1&2&0&2&1&0&0&0\\
4&1&1&1&0&1&0&1&0&1&0&0&0&0\\
12&5&5&2&1&5&2&2&1&0&0&0&0&0\\
6&2&2&1&0&2&1&1&0&0&0&0&0&0\\
8&4&3&2&1&2&1&0&0&0&0&0&0&0\\
4&1&1&0&0&1&0&0&0&0&0&0&0&0\\
6&3&3&2&1&0&0&0&0&0&0&0&0&0\\
3&1&1&1&0&0&0&0&0&0&0&0&0&0\\
3&1&1&0&0&0&0&0&0&0&0&0&0&0\\
2&1&0&0&0&0&0&0&0&0&0&0&0&0\\
1&0&0&0&0&0&0&0&0&0&0&0&0&0\\
\end{array} \right]. $$
Soit $P'_n=Pass({\cal B'}_n,(e_{F_i})_{i\leq r_n}).$ $P'_n={A'_n}^{-1}$, et donc :
\begin{eqnarray*}
 P'_1=\left[ \begin{array}{c}
           1
\end{array} \right] ;\hspace{8mm}
 P'_2=\left[ \begin{array}{cc}
0&1\\
1&-2
\end{array} \right] ;&&
 P'_3=\left[ \begin{array}{ccccc}
0&0&0&0&1\\
0&0&0&1&-2\\
0&0&1&-1&-1\\
0&1&-1&0&0\\
1&-2&-1&0&3
\end{array} \right] ;
\end{eqnarray*}
$$ P'_4=\left[ \begin{array}{cccccccccccccc}
0&0&0&0&0&0&0&0&0&0&0&0&0&1\\
0&0&0&0&0&0&0&0&0&0&0&0&1&-2\\
0&0&0&0&0&0&0&0&0&0&0&1&-1&-1\\
0&0&0&0&0&0&0&0&0&0&1&-1&0&0\\
0&0&0&0&0&0&0&0&0&1&-2&-1&0&3\\
0&0&0&0&0&0&0&0&1&0&0&-1&0&-1\\
0&0&0&0&0&0&0&1&-2&-1&0&1&-1&2\\
0&0&0&0&0&0&1&-1&0&1&-1&-1&1&0\\
0&0&0&0&0&1&-2&0&-1&-1&2&1&0&1\\
0&0&0&0&1&0&-1&1&-1&-1&0&2&-1&0\\
0&0&0&1&-2&0&0&-1&2&0&2&-2&0&0\\
0&0&1&-1&-1&-1&2&-1&1&2&-2&-2&2&0\\
0&1&-1&0&0&0&-1&1&0&-1&0&2&-1&0\\
1&-2&-1&0&3&-1&2&0&1&0&0&0&0&-4\\
\end{array}\right].$$

\end{document}